\documentclass[reqno]{amsart}
\pdfoutput=1
\usepackage{graphicx}
\usepackage{amsmath}
\usepackage{amsfonts}
\usepackage{amssymb}%
\providecommand{\U}[1]{\protect\rule{.1in}{.1in}}
\newtheorem{theorem}{Theorem}

\newtheorem{conclusion}[theorem]{Conclusion}

\newtheorem{corollary}[theorem]{Corollary}

\newtheorem{definition}[theorem]{Definition}

\newtheorem{lemma}[theorem]{Lemma}
\newtheorem{notation}[theorem]{Notation}

\newtheorem{proposition}[theorem]{Proposition}
\newtheorem{remark}[theorem]{Remark}

\numberwithin{theorem}{section}
\numberwithin{equation}{section}
\begin{document}
\title[Brown measure of Brownian motions]{The Brown measure of a family of free multiplicative Brownian motions}
\author{Brian C. Hall}
\address{Department of Mathematics, University of Notre Dame, Notre Dame, IN 46556, USA}
\email{bhall@nd.edu}
\author{Ching-Wei Ho}
\address{Institute of Mathematics, Academia Sinica, Taipei 10617, Taiwan; Department of
Mathematics, University of Notre Dame, Notre Dame, IN 46556, USA}
\email{chwho@gate.sinica.edu.tw}
\thanks{Hall's research supported in part by a grant from the Simons Foundation}

\begin{abstract}
We consider a family of free multiplicative Brownian motions $b_{s,\tau}$
parametrized by a real variance parameter $s$ and a complex covariance
parameter $\tau.$ We compute the Brown measure $\mu_{s,\tau}$ of $ub_{s,\tau
},$ where $u$ is a unitary element freely independent of $b_{s,\tau}.$ We find
that $\mu_{s,\tau}$ has a simple structure, with a density in logarithmic
coordinates that is constant in the $\tau$-direction. These results generalize
those of Driver--Hall--Kemp and Ho--Zhong for the case $\tau=s.$

We also establish a remarkable \textquotedblleft model deformation
phenomenon,\textquotedblright\ stating that all the Brown measures with $s$
fixed and $\tau$ varying are related by push-forward under a natural family of maps.

Our proofs use a first-order nonlinear PDE of Hamilton--Jacobi type satisfied
by the regularized log potential of the Brown measures. Although this approach
is inspired by the PDE method introduced by Driver--Hall--Kemp, our methods
are substantially different at both the technical and conceptual level.

\end{abstract}
\maketitle
\tableofcontents

\section{Introduction}

\subsection{Additive and multiplicative models and their large-$N$
limits\label{introAddMult.sec}}

The most basic random matrix models, such as the Gaussian unitary ensemble and
the Ginibre ensemble, are described by Gaussian measures. Thus, in light of
the central limit theorem, such a model can be thought of as arising as the
sum of a large number of small independent, identically distributed random
matrices. More explicitly, a Gaussian random matrix model can be described in
terms of Brownian motion in the appropriate vector space of matrices.

It is then natural to consider \textquotedblleft
multiplicative\textquotedblright\ random matrix models, which are constructed
as \textit{products} of independent random matrices that are close to the
identity. More explicitly, we can consider random matrix models arising from a
left-invariant Brownian motion in some group of invertible matrices, such as
the unitary group $U(N)$ or the general linear group $GL(N;\mathbb{C}).$ For
such a Brownian motion $g_{t}$, the \textit{multiplicative} increments,
computed as $g_{s}^{-1}g_{t}$ for $s<t,$ are independent.

In \cite{BianeFields}, Biane constructed the large-$N$ limit of Brownian
motions $U_{t}^{N}$ in $U(N)$, the so-called free unitary Brownian motion
$u_{t},$ as an element of an operator algebra with a trace. Biane then
computed the law of $u_{t},$ which is a probability measure $\nu_{t}$ on the
unit circle. Zhong \cite{Zhong} then extended Biane's result by computing the
law of $uu_{t}$, where $u$ is an arbitrary unitary element that is freely
independent of $u_{t}.$ Meanwhile, in \cite{BianeJFA}, Biane introduced a free
multiplicative Brownian motion $b_{t}$ (denoted $\Lambda_{t}$ in
\cite{BianeJFA}). As conjectured by Biane and proved by Kemp \cite{KempLargeN}%
, $b_{t}$ is the large-$N$ limit, in the sense of $\ast$-distribution, of a
Brownian motion $B_{t}^{N}$ in $GL(N;\mathbb{C}).$

The paper \cite{DHKBrown} computed the Brown measure of $b_{t}.$ This Brown
measure is supported in a domain $\Sigma_{t}$ introduced by Biane and has a
special structure. It is also related to the law $\nu_{t}$ of the free unitary
Brownian motion. The results of \cite{DHKBrown} were then extended by Ho and
Zhong \cite{HZ} to compute the Brown measure of $ub_{t},$ where $u$ is an
arbitrary unitary element freely independent of $b_{t}.$

It is possible to study also a two-parameter family of Brownian motions in
$GL(N;\mathbb{C}),$ where the two parameters control the diffusion rates in
the unitary and positive directions. This family was introduced in the context
of the Segal--Bargmann transform in \cite{DriverHall,HallCJM} and then used
also in the large-$N$ limit of the Segal--Bargmann transform in
\cite{DHKlargeN,KempLargeN,HoSBT}. This family of Brownian motions was also
studied from the point of view of random matrix theory (that is, looking at
the eigenvalue distribution) in \cite{LNW} and \cite{HKsupport}.

In \cite{DHKcomplex}, it was then noted that one can naturally include a third
parameter that controls the correlation between the unitary and positive
diffusions. According to \cite[Theorem 3.2]{DHKcomplex}, this three-parameter
family of Brownian motions is (up to multiplying by a scalar process) the most
general family of Brownian motions that is invariant under the left action of
$GL(N;\mathbb{C})$ and the right action of $U(N).$ See also the work of Chan
\cite{Chan} on the large-$N$ limit of the Segal--Bargmann transform associated
to this three-parameter family.

Following \cite{DHKcomplex}, it is convenient to describe these Brownian
motions by one real parameter $s$ and one complex parameter $\tau,$ where $s$
controls the overall diffusion rate, while the real and imaginary parts of
$\tau$ control the diffusion rate in the unitary direction and the
correlations between the unitary and positive diffusions. The complex
parameter $\tau$ is nonzero and satisfies%
\[
\left\vert \tau-s\right\vert \leq s.
\]
See Section \ref{freeBM.sec} for details of the parametrization. We therefore
obtain a family $B_{s,\tau}^{N}(r)$ of Brownian motions in $GL(N;\mathbb{C}),$
where $s$ and $\tau$ are parameters and $r$ is the time variable. In this
paper, we will consider a \textquotedblleft free\textquotedblright\ version
$b_{s,\tau}(r)$ of these Brownian motions. For any $r>0,$ the $\ast
$-distribution of $b_{s,\tau}(r)$ is the same as that of $b_{rs,r\tau}(1).$ It
is therefore sufficient to consider the case $r=1$ and we use the notation%
\begin{equation}
b_{s,\tau}:=b_{s,\tau}(1). \label{bstauDef}%
\end{equation}
The case $\tau=s$ gives the \textquotedblleft standard\textquotedblright\ free
multiplicative Brownian motion $b_{s},$ as defined by Biane and denoted by
$\Lambda_{s}$ in \cite[Section 4]{BianeJFA}. When $\tau$ is real, results of
Kemp \cite{KempLargeN} show that $b_{s,\tau}$ is the large-$N$ limit of
$B_{s,\tau}^{N}$ in the sense of $\ast$-distribution. We expect that this
result extends in a straightforward way to general values of $\tau,$ but this
has not been proved.

\begin{figure}[ptb]
\centering
\includegraphics[scale=0.6]{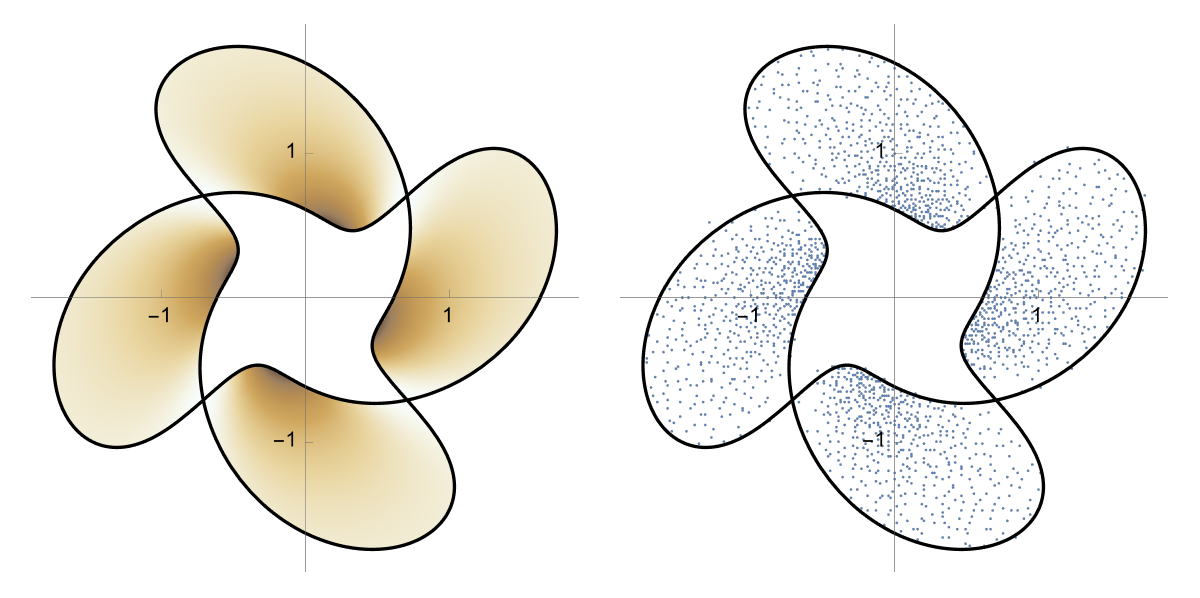}
\caption{Comparison of the Brown measure of $ub_{s,\tau}$ (left) with a
simulation of the eigenvalues of the corresponding random matrix model
(right). Shown for $s=1,$ $\tau=1+i/2,$ with the law of $u$ giving equal mass
to the four points $\pm1$ and $\pm i$.}
\label{theoryexp.fig}
\end{figure}

We then consider an arbitrary unitary element $u$ freely independent of
$b_{s,\tau}(r)$ and we consider the element%
\begin{equation}
ub_{s,\tau}. \label{bsZero}%
\end{equation}
We let $\mu_{s,\tau}$ denote the \textquotedblleft Brown
measure\textquotedblright\ of $ub_{s,\tau},$ as defined in Section
\ref{brown.sec}. We believe that the Brown measure of $ub_{s,\tau}$ coincides
with the limiting eigenvalue distribution of any random matrix model of the
form $U^{N}B_{s,\tau}^{N}\,$, where $U^{N}$ is unitary and independent of
$B_{s,\tau}^{N},$ and where and limiting eigenvalue distribution of $U^{N}$ is
the law of the unitary element $u$ in (\ref{bsZero}). Our belief is supported
by numerical simulations of $U^{N}B_{s,\tau}^{N},$ as in Figure
\ref{theoryexp.fig}.

In this paper, we will compute the Brown measure $\mu_{s,\tau}$ and find that
it has a remarkably simple structure, with a density in logarithmic
coordinates that is constant in the $\tau$-direction. (See Section
\ref{BrownFormula.sec}.) We will also establish a \textquotedblleft model
deformation phenomenon\textquotedblright\ as follows. Suppose we deform the
free random matrix model $ub_{s,\tau}$ by varying $\tau$ with $s$ and $u$
fixed. Then the Brown measures $\mu_{s,\tau}$ vary in a simple way: the Brown
measures are all related by push-forward under a natural family of maps. (See
Section \ref{Relating.sec}.) Assuming that the Brown measures are the same as
the limiting eigenvalue distributions of the associated random matrix models,
this result represents a new phenomenon in random matrix theory: the limiting
eigenvalue distribution of one model can be obtained by applying an explicitly
computable plane map to the eigenvalues of a \textit{different} random matrix model.

In \cite{HoElliptic}, the second author obtains the additive counterpart of
the results of the present paper, in the special case that $\tau$ is real.
That is to say, when $\tau$ is real, \cite{HoElliptic} gives the Brown measure
of $x_{0}+w_{s,\tau},$ where $w_{s,\tau}$ is the free (additive) Brownian
motion $w_{s,\tau}(r)$ in Section \ref{freeBM.sec} evaluated at $r=1$ and
where $x_{0}$ is a self-adjoint element freely independent of $w_{s,\tau}.$
The results of \cite{HoElliptic} are based on \cite{HHadditive}, which uses
the PDE\ method.

For general complex values of $\tau$ (with $\left\vert \tau-s\right\vert \leq
s$), the results of the present paper should have straightforward analogs for
the \textquotedblleft additive case,\textquotedblright\ that is, the Brown
measure of $x_{0}+w_{s,\tau},$ where $x_{0}$ is self-adjoint and freely
independent of $w_{s,\tau}.$ Indeed, all the proofs given here should adapt to
the additive case in a straightforward way, except that certain parts of the
argument become simpler. We plan to provide the details of this argument in a
later publication.

Meanwhile, after the first version of this paper was posted on the arXiv,
Zhong \cite{Zhong2} posted a preprint that computes the Brown measure of
$x_{0}+w_{s,\tau}$, where $x_{0}$ is assumed to be freely independent of
$w_{s,\tau}$ but is otherwise arbitrary. Zhong's results therefore go well
beyond the additive counterpart of our results here, which rely heavily on the
assumption that the element $u$ in (\ref{bsZero}) is unitary. Zhong uses free
probability in place of the PDE\ method.

\subsection{Random matrices and the Brown measure\label{brown.sec}}

Girko's method for determining the eigenvalues of a non-normal random matrix
$Z$ consists of looking at the function%
\[
s(\lambda):=\frac{1}{N}\log(\left\vert \det(Z-\lambda)\right\vert ^{2}%
)=\frac{2}{N}\sum_{j=1}^{N}\log\left\vert \lambda-\lambda_{j}\right\vert ,
\]
where $\lambda_{1},\ldots,\lambda_{N}$ are the eigenvalues of $Z.$ Since
$\frac{1}{2\pi}\log\left\vert \lambda-\lambda_{j}\right\vert $ is the Green's
function of the Laplacian in the plane, we find that $\frac{1}{4\pi}$ time the
Laplacian of $s,$ computed in the distribution sense, is just the empirical
eigenvalue distribution of $Z$:%
\[
\frac{1}{4\pi}\Delta s(\lambda)=\frac{1}{N}\sum_{j=1}^{N}\delta_{\lambda_{j}%
}(\lambda).
\]
It is then an elementary exercise to show that $s$ may also be computed as%
\begin{equation}
s(\lambda)=\frac{1}{N}\operatorname{trace}\left[  \log((Z-\lambda)^{\ast
}(Z-\lambda))\right]  , \label{sLambdaIntro}%
\end{equation}
where \textquotedblleft$\log$\textquotedblright\ denotes the self-adjoint
logarithm of a self-adjoint, positive matrix.

Meanwhile, pioneering work of Voiculescu developed an operator-algebra
formalism for understanding the large-$N$ limits of random matrix models.
Specifically, if $Z^{N}$ is a sequence of $N\times N$ random matrices (not
necessarily Hermitian or even normal) one looks for a von Neumann algebra
$\mathcal{A}$ with a \textquotedblleft trace\textquotedblright%
\ $\operatorname{tr}:\mathcal{A}\rightarrow\mathbb{C}$ and an element $z$ of
$\mathcal{A}$ such that%
\begin{equation}
\lim_{N\rightarrow\infty}\mathbb{E}\left\{  \frac{1}{N}\operatorname{trace}%
\left[  p(Z^{N},(Z^{N})^{\ast}\right]  \right\}  =\operatorname{tr}\left[
p(z,z^{\ast})\right]  \label{starLim}%
\end{equation}
for every polynomial $p$ in two noncommuting variables. (Here
$\operatorname{tr}$ is \textit{not} the trace in the sense of trace-class
operators, but is a linear functional having properties similar to the
normalized trace of matrices.) One indication of the power of this approach is
the work of Voiculescu \cite{Voi2,VoiStrengthened}, which showed that, under
quite general circumstances, independent random matrices become \textit{freely
}independent in such limits. (See also Chapter 23 of \cite{NicaSpeicher} and
Section 4.2 of \cite{MS}.)

If $z$ is an element of a tracial von Neumann algebra $(\mathcal{A}%
,\operatorname{tr}),$ one can define the \textbf{Brown measure} of $z$
\cite{Br} denoted $\mu^{z},$ by imitating the preceding construction for
matrices. (See also Chapter 11 of \cite{MS}.) One first introduces a
regularized analog of (\ref{sLambdaIntro}),%
\begin{equation}
S(\lambda,\varepsilon)=\operatorname{tr}\left[  \log((z-\lambda)^{\ast
}(z-\lambda)+\varepsilon^{2})\right]  ,\quad\varepsilon>0, \label{generalS}%
\end{equation}
and then defines
\[
S_{0}(\lambda)=\lim_{\varepsilon\rightarrow0^{+}}S(\lambda,\varepsilon).
\]
The Brown measure $\mu^{z}$ is then $\frac{1}{4\pi}$ times the distributional
Laplacian of $S_{0}$:%
\[
\mu^{z}=\frac{1}{4\pi}\Delta S_{0}.
\]

An important technical issue with this approach is that if even if $Z^{N}$
converges to $z$ in the sense described in (\ref{starLim}), the empirical
eigenvalue distribution of $Z^{N}$ need not converge almost surely to the
Brown measure of $z.$ On the other hand, tools have been developed to show
that, in many interesting cases, $\mu^{z}$ is indeed the limiting eigenvalue
distribution. We mention here work of Girko \cite{GirkoCircular}, Bai
\cite{Bai}, Tao--Vu \cite{TV}, and Guionnet--Krishnapur--Zeitouni
\cite{GKZsingleRing}, among others. (Although these works do not explicitly
use the Brown measure terminology, they all use Girko's approach to computing
the eigenvalue distribution in the large-$N$ limit.)

\subsection{Formula for the Brown measure\label{BrownFormula.sec}}

The first main result of this paper is a computation of the following measure:%
\begin{subequations}
\[
\mu_{s,\tau}=\text{ Brown measure of }ub_{s,\tau},
\]
where the answer will depend on the choice of the unitary element $u$. The
case $\tau=s$ corresponds to the Brown measure computed in \cite{DHKBrown,HZ}.%

\begin{figure}[ptb]
\centering
\includegraphics[scale=0.55]{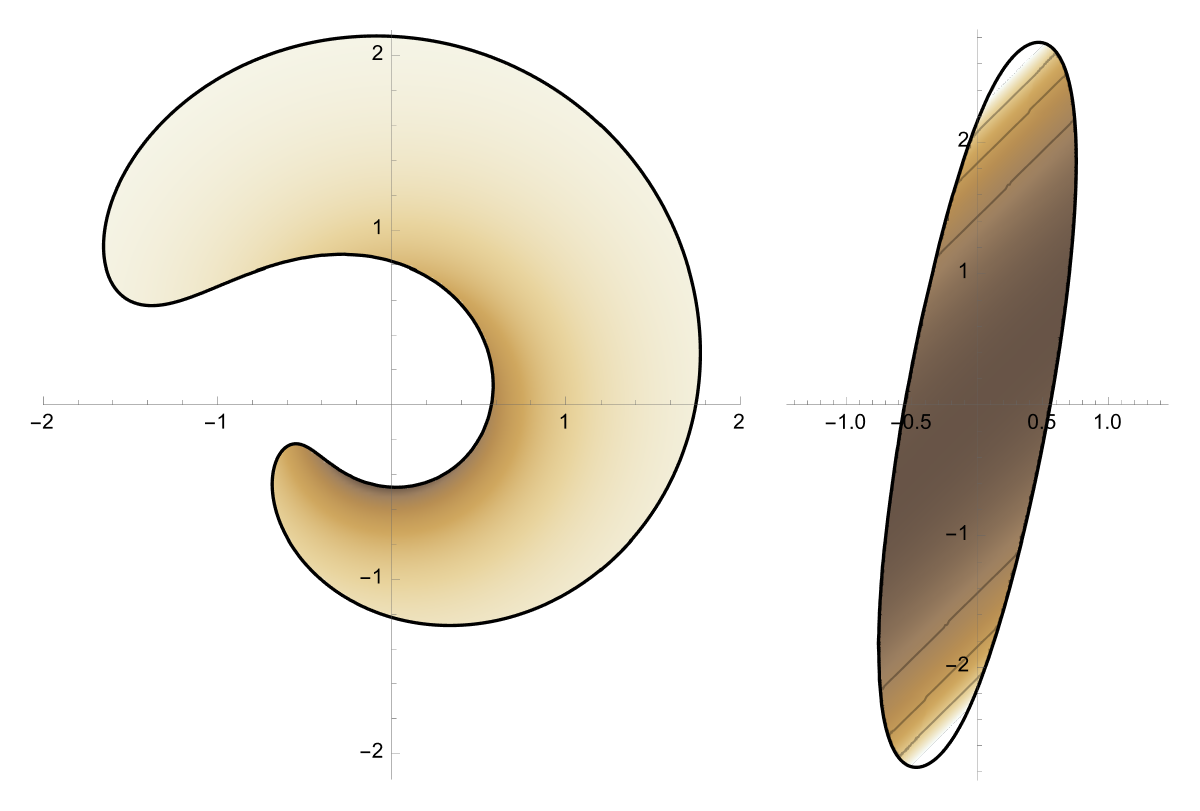}
\caption{The Brown measure $\mu_{s,\tau}$ in rectangular coordinates (left)
and in logarithmic coordinates (right). Shown with $u=1,$ $s=3,$ and
$\tau=1+i.$ The density in logarithmic coordinates is constant in the $\tau
$-direction.}
\label{brownwithlog.fig}
\end{figure}

We now briefly indicate the structure of the Brown measure, leaving a precise
statement of the results to Section \ref{results.sec}. We first identify a
certain open set $\Sigma_{s,\tau}$ containing all the mass of the Brown
measure $\mu_{s,\tau}$ and then give a formula for the density of $\mu
_{s,\tau}$ on $\Sigma_{s,\tau}.$ We use the notation%
\end{subequations}
\[
\tau=\tau_{1}+i\tau_{2},\quad\tau_{1},\tau_{2}\in\mathbb{R},
\]
and introduce twisted logarithmic coordinates of a nonzero complex number
$\lambda,$ given as
\begin{equation}
v=\frac{1}{\tau_{1}}\log\left\vert \lambda\right\vert ;\quad\delta=\arg
\lambda-\frac{\tau_{2}}{\tau_{1}}\log\left\vert \lambda\right\vert .
\label{coordinates}%
\end{equation}
These coordinates are the unique real numbers $v$ and $\delta$ such that
\[
\lambda=e^{v\tau}e^{i\delta}.
\]

We are now ready to state our first main result.

\begin{theorem}
\label{brownIntro.thm}There is a smooth, increasing function $\phi^{s,\tau}$
such that the Brown measure in $\Sigma_{s,\tau}$ is given by the formula%
\begin{equation}
d\mu_{s,\tau}(\lambda)=\frac{1}{2\pi\tau_{1}}\frac{1}{\left\vert
\lambda\right\vert ^{2}}\frac{d\phi^{s,\tau}(\delta)}{d\delta}~dx~dy,
\label{introFormula}%
\end{equation}
where $\lambda=x+iy.$ We can also write the result in logarithmic coordinates,
$\rho=\log\left\vert \lambda\right\vert $ and $\theta=\arg\lambda,$ as%
\[
d\mu_{s,\tau}=\frac{1}{2\pi\tau_{1}}\frac{d\phi^{s,\tau}(\delta)}{d\delta
}~d\rho~d\theta.
\]
Since the map $(\rho,\theta)\mapsto(\rho+\tau_{1},\theta+\tau_{2})$ does not
change the value of the $\delta$-coordinate, we see that \textit{the density
of }$\mu_{s,\tau}$\textit{ in logarithmic coordinates is constant in the
}$\tau$\textit{-direction}.
\end{theorem}

The geometric meaning of the function $\phi^{s,\tau}$ is explained in Section
\ref{mainResults.sec}; see Figures \ref{thetaphi.fig} and \ref{thetadelta.fig}
and Notation \ref{phistau.notation}. We could also write the Brown measure in
the \textit{twisted} logarithmic coordinates $v$ and $\delta,$ but this change
only affects the density by a constant, since $v$ and $\delta$ are obtained
from $\rho$ and $\theta$ by a linear change of variables.

\begin{figure}[ptb]
\centering
\includegraphics[scale=0.65]{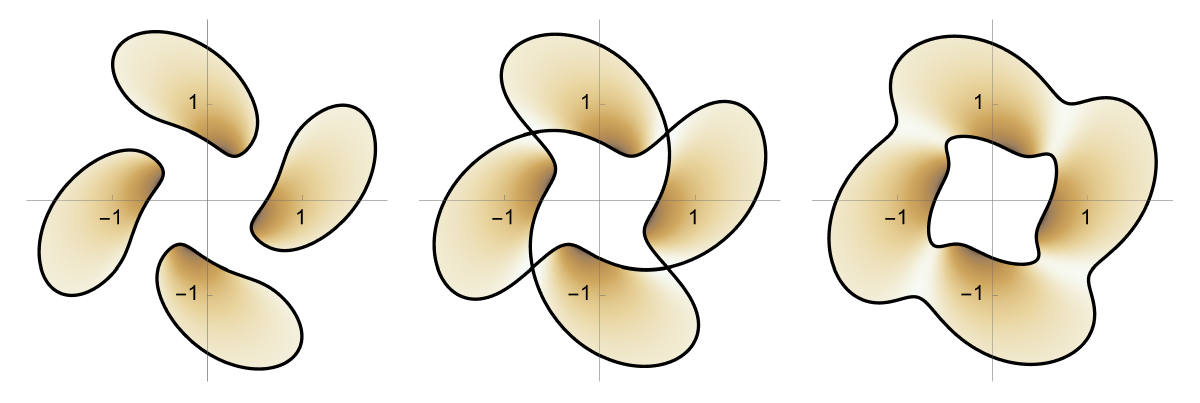}
\caption{The Brown measure $\mu_{s,\tau}$ with $\tau=1+i/2$ and the law of $u$
putting equal mass at the four points $\pm1$ and $\pm i.$ Shown for $s=0.9,$
$s=1,$ and $s=1.1.$}
\label{threebrowns.fig}
\end{figure}

Figure \ref{brownwithlog.fig} gives a density plots of an example of the Brown
measure, in both rectangular and logarithmic coordinates. See also Figure
\ref{threebrowns.fig} for additional examples of Brown measures, showing how
the Brown measure changes as $s$ changes with $\tau$ fixed. See Section
\ref{inside.sec} for details about Theorem \ref{brownIntro.thm}.%

\subsection{Relating different values of $\tau$\label{Relating.sec}}

Our second main result is that all the Brown measures of $ub_{s,\tau}$ with
$s$ and $u$ fixed and $\tau$ varying, are related.

\begin{theorem}
\label{pushforward.thm}Fix a positive real number $s$ and the unitary element
$u$ and consider a nonzero complex number $\tau$ satisfying $\left\vert
s-\tau\right\vert \leq s.$ Then there is an invertible map $\Phi_{s,\tau
}:\Sigma_{s}\rightarrow\Sigma_{s,\tau}$ with the following properties. First,
$\Phi_{s,\tau}$ maps each radial segment in $\Sigma_{s}$ to a curve with a
constant value of the coordinate $\delta$ in $\Sigma_{s,\tau}$; this curve is
a portion of an exponential spiral. Second, the associated Brown measures
$\mu_{s,s}$ and $\mu_{s,\tau}$ are related by push-forward:
\[
(\Phi_{s,\tau})_{\ast}(\mu_{s,s})=\mu_{s,\tau}.
\]

Furthermore, let $\Phi_{s}$ denote the limit of $\Phi_{s,\tau}$ as $\tau$
tends to zero. Then $\Phi_{s}$ maps $\Sigma_{s}$ into the unit circle and the
push-forward of $\mu_{s,s}$ by $\Phi_{s}$ is the law of $uu_{s},$ where
$u_{s}$ is the free unitary Brownian motion and $u$ is as in (\ref{bsZero}).
\end{theorem}

\begin{figure}[ptb]
\centering
\includegraphics[scale=0.6]{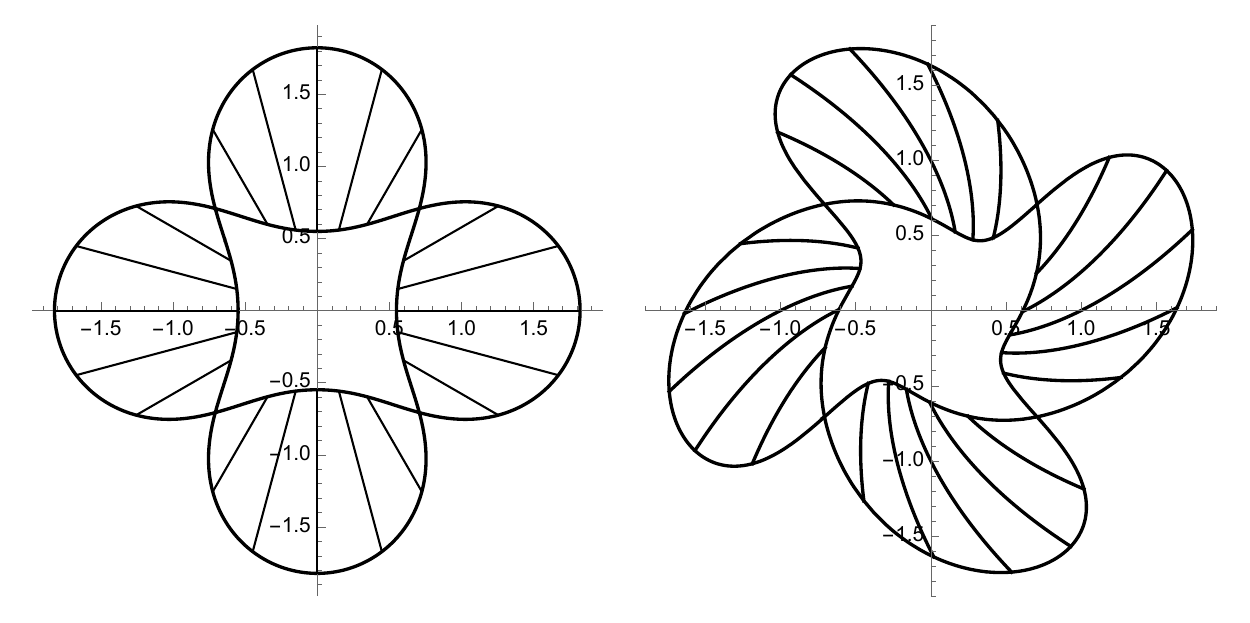}
\caption{The map $\Phi_{s,\tau}$ maps radial segments in $\Sigma_{s}$ (left)
to spiral segments in $\Sigma_{s,\tau}$ (right). Shown for $s=1$ and
$\tau=1+i/2,$ with the law of $u$ giving equal mass to the four points $\pm1$
and $\pm i$.}
\label{pushmap0.fig}
\end{figure}

\begin{figure}[ptb]
\centering
\includegraphics[scale=0.55]{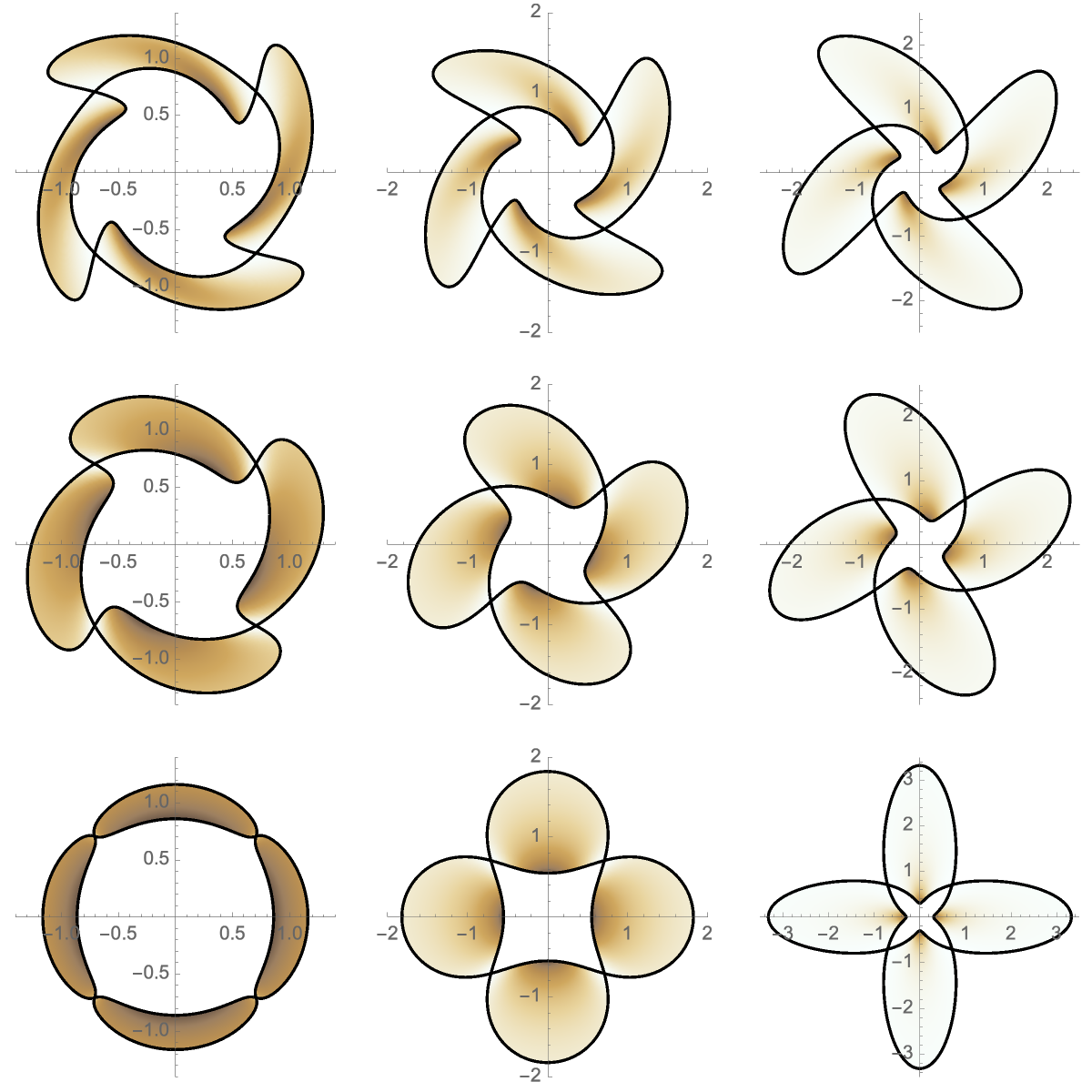}
\caption{The domains $\Sigma_{s,\tau}$ and associated Brown measures for $s=1$
and several different values of $\tau$, with the law of $u$ giving equal mass
to the four points $\pm1$ and $\pm i.$ The real part of $\tau$ increases from
left to right and the imaginary part of $\tau$ is zero on the bottom row and
increases from bottom to top. The case $\tau=s$ is in the middle of the bottom
row.}
\label{grid1.fig}
\end{figure}

See Section \ref{push.sec} for details. The map $\Phi_{s,\tau}$ is illustrated
in Figure \ref{pushmap0.fig}. Figure \ref{grid1.fig} then shows the domains
and Brown measures for a fixed value of $s$ and several different values of
$\tau.$ Although the different Brown measures in the figure look quite
different from one another, they are actually all connected by maps of the
form $\Phi_{s,\tau}\circ\Phi_{s,\tau^{\prime}}^{-1}.$ In particular, the
topology of the closed support of $\mu_{s,\tau}$ is independent of $\tau$ with
$s$ fixed.

Theorem \ref{pushforward.thm} demonstrates what we call a \textquotedblleft
model deformation phenomenon.\textquotedblright\ As we noted in Section
\ref{introAddMult.sec}, we believe that the Brown measure $\mu_{s,\tau}$ is
the limiting eigenvalue distribution of the random matrix model $U^{N}%
B_{s,\tau}^{N}.$ Assuming this is the case, Theorem \ref{pushforward.thm}
tells us that the limiting eigenvalue distribution of one random matrix model
$U^{N}B_{s,\tau}^{N}$ can be converted into the limiting eigenvalue
distribution of another random matrix model $U^{N}B_{s,\tau^{\prime}}^{N}$ by
applying the map $\Phi_{s,\tau}\circ\Phi_{s,\tau^{\prime}}^{-1}$ to the
eigenvalues of the first model. As we will explain in the next paragraph, the
limiting case of this result as $\tau\rightarrow0$ was established in
\cite{DHKBrown} and \cite{HZ}. Otherwise, it seems to be a new phenomenon in
random matrix theory. (After the first version of this paper was posted on the
arXiv, Zhong posted a preprint \cite{Zhong2} obtaining similar results in the
additive case.)

The paper \cite{DHKBrown} showed that the push-forward of the Brown measure of
$b_{s,s}$ under a certain map $\Phi_{s}$ is equal to the law of the free
unitary Brownian motion $u_{s}.$ This result was extended in \cite{HZ} to
relate the Brown measure of $ub_{s}$ to the law of $uu_{s},$ where $u$ is an
arbitrary unitary element freely independent of $b_{s}$ and of $u_{s}.$ The
map $\Phi_{s}$ in \cite{DHKBrown,HZ} is the same as the one in the last part
of Theorem \ref{pushforward.thm}, showing that the push-forward result in
those papers is simply a limiting case ($\tau\rightarrow0$) of a much more
general family of results.%

\subsection{The method of proof}

The proofs use a variant of the PDE method introduced in \cite{DHKBrown},
which studied the \textquotedblleft standard\textquotedblright\ free
multiplicative Brownian motion with trivial initial condition. (See also
\cite{PDEmethods} for a gentle introduction to the PDE method.) In the
notation of the present paper, \cite{DHKBrown} corresponds to taking $\tau=s$
and $u=1$. The method in \cite{DHKBrown} was then extended by the second
author and Zhong \cite{HZ} to analyze the standard free multiplicative
Brownian motion with arbitrary unitary initial condition ($\tau=s$ and $u$ is
a arbitrary). The method has also been used in \cite{DemniHamdi} and
\cite{HHadditive} and discussed in the physics literature in \cite{GNT}.

Although the results of the present paper include the results of
\cite{DHKBrown} and \cite{HZ} as special cases, the PDE we use here
\textit{does not} include the PDE in those papers as a special case. The PDEs
differ because in our PDE we vary $\tau$ with $s$ fixed, whereas in
\cite{DHKBrown,HZ}, $\tau$ is always equal to $s.$ This point is discussed in
detail in Section \ref{discussion.sec}. Although the approach we are using has
substantial advantages (described in Section \ref{discussion.sec}) relative to
the approach in \cite{DHKBrown,HZ}, it also requires new techniques.
Specifically, the technical details associated to analyzing the Brown measure
in the region where it is not zero are substantially different here than in
previous papers. See the beginning of Section \ref{inside.sec} for a
discussion of the relevant differences. To derive the relevant PDE, we use a
factorization result (Theorem \ref{thm:factorization}) for elements of the
form $b_{s,\tau}$ in (\ref{bstauDef}) that is of independent interest. The
proof of this result is surprisingly subtle; see Appendix
\ref{factorization.app}.

\subsection{Table of notation}

We collect the following notations.%

\[%
\begin{array}
[c]{ccl}
&  & \text{\textbf{Table of Notation}}\\
u_{s} &  & \text{free unitary Brownian motion starting at 1}\\
b_{s,\tau} &  & \text{free multiplicative Brownian motion in Notation
\ref{bstau.notation}}\\
u &  & \text{unitary element freely independent of }b_{s,\tau}\\
\mu_{0} &  & \text{law of }u\\
\mu_{s} &  & \text{law of }uu_{s}\\
m_{s} &  & \text{density of }\mu_{s}\text{ w.r.t. normalized Lebesgue measure
on }S^{1}\\
\mu_{s,\tau} &  & \text{Brown measure of }ub_{s,\tau}\\
f_{\beta} &  & \text{holomorphic function in Definition \ref{fFunct.def}}\\
\Sigma_{s}\text{ and }\Sigma_{s,\tau} &  & \text{domains in Definitions
\ref{sigmas.def} and \ref{def:Sigmastaudef}}\\
S(s,\tau,\lambda,\varepsilon) &  & \text{regularized log potential of
}ub_{s,\tau}\text{ in (\ref{Sfunction})}\\
S_{0}(s,\tau,\lambda) &  & \text{limiting value of }S(s,\tau,\lambda
,\varepsilon)\text{ as }\varepsilon\text{ tends to zero from above}\\
v\text{ and }\delta &  & \text{twisted logarithmic coordinates defined in
(\ref{coordinates})}\\
\tau_{1}\text{ and }\tau_{2} &  & \text{real and imaginary parts of }\tau
\end{array}
\]

\section{Set-up and statement of results\label{results.sec}}

\subsection{The free Brownian motions\label{freeBM.sec}}

We now consider a tracial von Neumann algebra $\mathcal{A},$ that is, a von
Neumann algebra together with a faithful, normal, tracial state. To avoid a
conflict of notation with the complex parameter $\tau,$ we denote the trace on
$\mathcal{A}$ by
\[
\operatorname{tr}:\mathcal{A}\rightarrow\mathbb{C},
\]
while emphasizing that $\operatorname{tr}\left[  \cdot\right]  $ is not the
trace in the sense of trace-class operators. Rather, $\operatorname{tr}\left[
\cdot\right]  $ is a linear functional that behaves like the normalized trace
of matrices.

In free probability, a semicircular element is the large-$N$ limit of a GUE.
More precisely, an element $X$ of a tracial von Neumann algebra $(\mathcal{A}%
,\operatorname{tr})$ is \textbf{semicircular with variance }$t$ if its law is
the semicircular measure supported on the interval $[-2\sqrt{t},2\sqrt{t}].$
There then exists a \textbf{semicircular Brownian motion} $x_{t}$,
characterized as a continuous process with $x_{0}=0$ and having freely
independent increments such that $x_{t}-x_{s}$ is semicircular with variance
$t-s$ for all $t>s.$ See Section 1.1 of \cite{BianeFields}.

We then define a \textbf{rotated elliptic element} to be one of the form%
\begin{equation}
Z=e^{i\theta}(aX+ibY), \label{rotatedElliptic}%
\end{equation}
where $X$ and $Y$ are freely independent semicircular elements, $a$, $b,$ and
$\theta$ are real numbers, and we assume $a$ and $b$ are not both zero. The
$\ast$-distribution of $Z$ is unchanged if we switch the sign of $a$ or $b.$
The $\ast$-distribution is also unchanged if we interchange $a$ and $b$ and
then replace $\theta$ by $\theta+\pi.$ Finally, if $a=b,$ then $Z$ is a
circular element with $\ast$-distribution independent of $\theta.$

We parametrize rotated elliptic elements by two parameters: a positive
variance parameter $s$ and a complex covariance parameter $\tau$ defined by%
\begin{align}
s  &  =\operatorname{tr}\left[  Z^{\ast}Z\right] \label{rot1}\\
\tau &  =\operatorname{tr}\left[  Z^{\ast}Z\right]  -\operatorname{tr}\left[
Z^{2}\right]  . \label{rot2}%
\end{align}
By the Cauchy--Schwarz inequality for the inner product $\left\langle
A,B\right\rangle :=\operatorname{tr}(A^{\ast}B),$ any rotated elliptic element
must satisfy%
\[
\left\vert \tau-s\right\vert \leq s.
\]
We can recover the parameters $a,$ $b,$ and $\theta$---up to the
transformations in the previous paragraph---by plugging (\ref{rotatedElliptic}%
) into (\ref{rot1}) and (\ref{rot2}), giving
\begin{align*}
\operatorname{tr}\left[  Z^{\ast}Z\right]   &  =a^{2}+b^{2}=s\\
\operatorname{tr}\left[  Z^{2}\right]   &  =e^{2i\theta}(a^{2}-b^{2})=s-\tau
\end{align*}
where we have used that $\operatorname{tr}\left[  X^{2}\right]
=\operatorname{tr}\left[  Y^{2}\right]  =1$ and $\operatorname{tr}\left[
XY\right]  =\operatorname{tr}\left[  X\right]  \operatorname{tr}\left[
Y\right]  =0.$ If we assume $a\geq b\geq0,$ then we easily find%
\begin{align}
\theta &  =\frac{1}{2}\arg(s-\tau)\nonumber\\
a  &  =\sqrt{\frac{1}{2}(s+\left\vert \tau-s\right\vert )}\nonumber\\
b  &  =\sqrt{\frac{1}{2}(s-\left\vert \tau-s\right\vert )}. \label{parameters}%
\end{align}

In the special case $\tau=s,$ we may take $a=b=\sqrt{s/2}$ and choose $\theta$
is arbitrarily. If $\left\vert \tau-s\right\vert =s,$ we have $a=\sqrt{s}$ and
$b=0,$ so that $Z=e^{i\theta}\sqrt{s}X.$ In particular, if $\tau=0,$ then $Z$
is just a semicircular element of variance $s.$

We then define a \textbf{free additive }$(s,\tau)$\textbf{-Brownian motion} as
a continuous process $w_{s,\tau}(r)$ with $w_{s,\tau}(0)=0$ having freely
independent increments such that for all $r_{2}>r_{1},$
\begin{equation}
\frac{w_{s,\tau}(r_{2})-w_{s,\tau}(r_{1})}{\sqrt{r_{2}-r_{1}}} \label{wIncr}%
\end{equation}
is a rotated elliptic element with parameters $s$ and $\tau.$ Explicitly, we
can construct $w_{s,\tau}(r)$ as%
\[
w_{s,\tau}(r)=e^{i\theta}(aX_{r}+ibY_{r})
\]
where $a$ and $b$ are chosen as above and where $X_{r}$ and $Y_{r}$ are freely
independent semicircular Brownian motions.

We then construct a \textbf{free multiplicative }$(s,\tau)$\textbf{-Brownian
motion} $b_{s,\tau}(r)$ as the solution to the free stochastic differential
equation%
\begin{equation}
db_{s,\tau}(r)=b_{s,\tau}(r)\left(  i~~dw_{s,\tau}(r)-\frac{1}{2}%
(s-\tau)~dr\right)  , \label{bstauSDE}%
\end{equation}
with
\[
b_{s,\tau}(0)=1.
\]
The $dr$ term on the right-hand side of (\ref{bstauSDE}) is an
\textquotedblleft It\^{o} term\textquotedblright\ equal to half the square of
$i~dw_{s,\tau}(r)$, computed using the It\^{o} rules in Section
\ref{itoRules.sec}. Existence and uniqueness of the solution to
(\ref{bstauSDE}) follows from a standard Picard iteration argument, as in
Proposition A.1 of \cite{CapDon}. (The key making the Picard iteration work is
the Burkholder--Davis--Gundy estimate for the norm of a free stochastic
integral, as in \cite[Theorem 3.2.1]{BS1} or \cite[Theorem 3.1.12]{Nik}.) It
follows from (\ref{wIncr}) that $w_{s,\tau}(r)\cong w_{rs,r\tau}(1),$ where
$\cong$ indicates that the two elements have the same $\ast$-distribution. It
is then not hard to see that
\[
b_{s,\tau}(r)\cong b_{rs,r\tau}(1),
\]
We will therefore usually assume that $r=1,$ without loss of generality.

We also fix a unitary element $u$ that is freely independent of $b_{s,\tau
}(r).$

\begin{notation}
\label{bstau.notation}For $s>0$ and $\left\vert \tau-s\right\vert \leq s,$ we
use $b_{s,\tau}$ to denote the value of $b_{s,\tau}(r)$ at $r=1$:%
\[
b_{s,\tau}:=b_{s,\tau}(1).
\]
Then we define%
\[
\mu_{s,\tau}=\text{ Brown measure of }ub_{s,\tau}.
\]

\end{notation}

When $\tau=0,$ we may take $\theta=0,$ $a=\sqrt{s},$ and $b=0$ in
(\ref{parameters}) so that $w_{s,\tau}(r)=\sqrt{s}X_{r}.$ In that case, the
SDE (\ref{bstauSDE}) becomes%
\[
db_{s,\tau}(r)=b_{s,\tau}(r)\left(  i\sqrt{s}~dX_{r}-\frac{1}{2}s~dr\right)
,
\]
and the solution is a free unitary Brownian motion, with time-parameter scaled
by $s$. We express this result as%
\begin{equation}
\left.  b_{s,\tau}(r)\right\vert _{\tau=0}\cong u_{sr}, \label{AtTauZero}%
\end{equation}
where $u_{s}$ is a free unitary Brownian motion with initial condition $1,$ as
introduced by Biane in \cite{BianeFields}.

\subsection{Main results\label{mainResults.sec}}

Throughout the paper, we will assume that $\tau$ is a nonzero complex number
satisfying $\left\vert \tau-s\right\vert \leq s$ . We always write $\tau$ in
the form%
\[
\tau=\tau_{1}+i\tau_{2},\quad\tau_{1},\tau_{2}\in\mathbb{R}.
\]
We also assume that the unitary element $u$ is fixed and we define a
probability measure $\mu_{0}$ on the unit circle by%
\begin{equation}
\mu_{0}=\text{ law of }u. \label{mu0law}%
\end{equation}
That is to say, $\mu_{0}$ is the unique probability measure on $S^{1}$ such
that%
\[
\int_{S^{1}}\xi^{k}~d\mu_{0}(\xi)=\operatorname{tr}\left[  u^{k}\right]
\]
for every integer $k.$ Although most of the objects defined below depend on
the measure $\mu_{0},$ we do not indicated this dependence in the notation. We
also consider the free unitary Brownian motion $u_{s}$ introduced in
\cite{BianeFields} and let
\[
\mu_{s}=\text{law of }uu_{s}.
\]
The measure $\mu_{s}$ was computed by Zhong in \cite{Zhong}. In particular,
Proposition 3.6 of Zhong states that $\mu_{s}$ is absolutely continuous with
respect to the Lebesgue measure on the circle.

\begin{definition}
\label{fFunct.def}For any complex number $\beta$, we define a function
$f_{\beta}$ by
\[
f_{\beta}(z)=z\exp\left\{  \frac{\beta}{2}\int_{S^{1}}\frac{\xi+z}{\xi
-z}\,d\mu_{0}(\xi)\right\}
\]
for any $z$ for which the integral is absolutely convergent.
\end{definition}

In the case that $\mu_{0}$ is the $\delta$-measure at $1$ and $\beta=s>0$, the
function $f_{s}$ reduces to the one defined by Biane in \cite{BianeJFA} and
used in \cite{DHKBrown} to compute the Brown measure of $b_{s,s}$. The general
definition of $f_{s}$ was given by Zhong in \cite{Zhong} and was used in
\cite{HZ} to compute the Brown measure of $ub_{s,s}$. In general, the function
$f_{\beta}$ is certainly defined and holomorphic on the open set
$\mathbb{C}\setminus\operatorname{supp}(\mu_{0}),$ where $\operatorname{supp}%
(\mu_{0})$ is the closed support (in $S^{1}$) of the measure $\mu_{0}.$ In
this paper, the motivation for Definition \ref{fFunct.def} is that the
characteristic curves of a certain family of PDEs will be expressed in terms
of functions of the form $f_{\beta}$; see Remark \ref{domainMotivation.remark}.

It is easily verified that for any positive real number $s,$ the function
$f_{s}$ satisfies%
\begin{equation}
f_{s}(z)=\left(  \overline{f(1/\bar{z})}\right)  ^{-1}. \label{fsZbar}%
\end{equation}
By Proposition 2.3 of \cite{Zhong} (drawing on results of \cite{BelBer}),
$f_{s}$ has a right inverse function $\chi_{s}$ defined on the unit disk. We
may then extend $\chi_{s}$ holomorphically to the set of $z$ such that
$\left\vert z\right\vert \neq1$ by defining
\begin{equation}
\chi_{s}(z)=\left(  \overline{\chi_{s}\left(  1/\bar{z}\right)  }\right)
^{-1}, \label{chisZbar}%
\end{equation}
for $\left\vert z\right\vert >1.$ By (\ref{fsZbar}), this extended $\chi_{s}$
is still a right inverse to $f_{s}.$ Then by \cite[Theorem 3.8]{Zhong},
$\chi_{s}$ extends by continuity to a neighborhood of any point $\xi$ in the
unit circle outside the closed support of the measure $\mu_{s}.$ Then
$\chi_{s}$ is holomorphic on $\mathbb{C}\setminus\operatorname{supp}(\mu
_{s}).$ By applying the Schwarz reflection principle to the function
$w\mapsto\frac{1}{i}\log\chi_{s}(e^{iw}),$ we see that $\chi_{s}$ is
holomorphic on $\mathbb{C}\setminus\mathrm{supp}(\mu_{s}).$ It also satisfies%
\[
f_{s}(\chi_{s}(z))=z,\quad z\in\mathbb{C}\setminus\operatorname{supp}(\mu
_{s}).
\]

We then define a domain $\Sigma_{s}$ as the interior of the complement of the
image of $\chi_{s},$ as follows.

\begin{definition}
\label{sigmas.def}For any $s>0,$ we define sets $\tilde{\Sigma}_{s}$ and
$\Sigma_{s}$ by%
\begin{align*}
\tilde{\Sigma}_{s}  &  =\left[  \chi_{s}(\mathbb{C}\setminus
\operatorname{supp}(\mu_{s}))\right]  ^{c}\\
\Sigma_{s}  &  =\text{interior of }\tilde{\Sigma}_{s}.
\end{align*}

\end{definition}

\begin{figure}[ptb]
\centering
\includegraphics[scale=0.8]{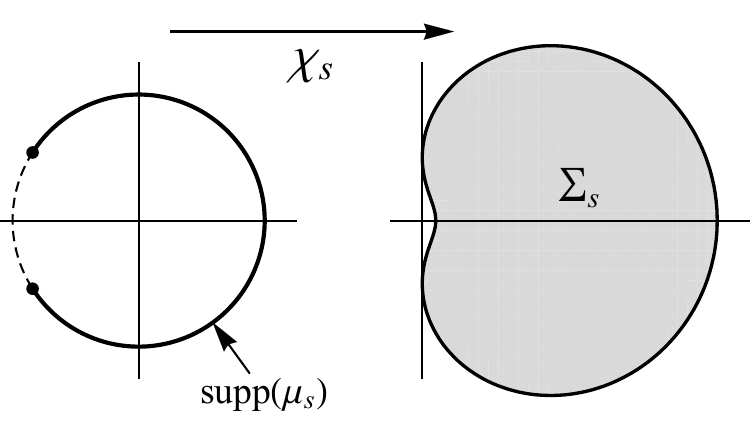}
\caption{The map $\chi_{s}=f_{s}^{-1}$ maps the complement of the support of
$\mu_{s}$ (left) to the complement of the closure of $\Sigma_{s}$ (right).
Shown for $s=2$ and $\mu_{0}=\delta_{1}$.}
\label{first2domains.fig}
\end{figure}

Since Proposition \ref{SigmasChar.prop} will show that $\tilde{\Sigma}_{s}$ is
the closure of $\Sigma_{s},$ we will subsequently use the notation
$\overline{\Sigma}_{s}$ in place of $\tilde{\Sigma}_{s}.$ See Figure
\ref{first2domains.fig}.%

In Section \ref{domSigmas.sec}, we will give the following characterization of
the domain $\Sigma_{s}.$

\begin{proposition}
\label{SigmasChar.prop}There is a unique continuous function $r_{s}$ on the
unit circle satisfying $0<r_{s}(\theta)\leq1$ such that%
\[
\left\vert f_{s}(r_{s}(\theta)e^{i\theta})\right\vert =1
\]
and such that the domain $\Sigma_{s}$ may be computed as
\[
\Sigma_{s}=\{\left.  re^{i\theta}\right\vert r_{s}(\theta)<r<1/r_{s}%
(\theta)\}.
\]
Furthermore, we have
\[
\tilde{\Sigma}_{s}=\overline{\Sigma}_{s}.
\]
We also have that%
\[
f_{s}(r_{s}(\theta)e^{i\theta})=f_{s}(r_{s}(\theta)^{-1}e^{i\theta}).
\]

\end{proposition}

We now introduce a domain $\Sigma_{s,\tau},$ which we define by specifying its complement.

\begin{definition}
\label{def:Sigmastaudef}Fix a positive real number $s$ and a nonzero complex
number $\tau$ with $\left\vert s-\tau\right\vert \leq s.$ Let $\Sigma_{s,\tau
}$ be the domain defined by
\begin{align}
(\tilde{\Sigma}_{s,\tau})^{c}  &  =f_{s-\tau}(\tilde{\Sigma}_{s}%
^{c})=f_{s-\tau}(\chi_{s}(\mathbb{C\setminus}\operatorname{supp}(\mu
_{s})))\label{stauTilde}\\
\Sigma_{s,\tau}  &  =\text{interior of }\tilde{\Sigma}_{s,\tau}.
\label{stauDef}%
\end{align}

\end{definition}

\begin{figure}[ptb]
\centering
\includegraphics[scale=0.8]{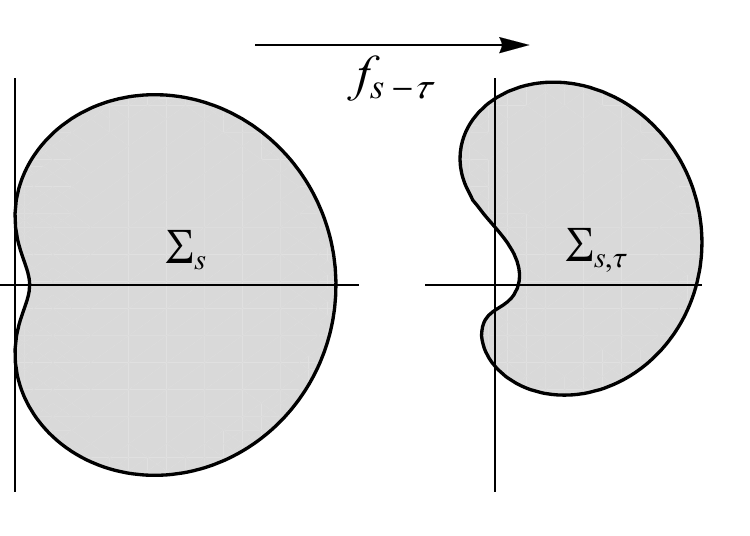}
\caption{The map $f_{s-\tau}$ maps the complement of the closure of
$\Sigma_{s}$ (left) to the complement of the closure of $\Sigma_{s,\tau}$
(right). Shown for $s=2$, $\tau=3/2+i/2,$ and $\mu_{0}=\delta_{1}.$}%
\label{last2domains.fig}
\end{figure}

We will see in Section \ref{domains.sec} that $\tilde{\Sigma}_{s,\tau}$ is the
closure of $\Sigma_{s,\tau},$ so we will subsequently use the notation
$\bar{\Sigma}_{s,\tau}$ in place of $\tilde{\Sigma}_{s,\tau}.$ See Figure
\ref{last2domains.fig}. In the case that $\tau$ is real and $u=1,$ the domain
$\Sigma_{s,\tau}$ and the map $f_{s-\tau}\circ\chi_{s}$ were studied in the
paper \cite{HoSBT} of the second author.%

\begin{remark}
\label{domainMotivation.remark}The way (the complement of) the domain
$\Sigma_{s,\tau}$ arises in this paper is the following. To compute the Brown
measure $\mu_{s,\tau}$ of $ub_{s,\tau},$ we consider the regularized
logarithmic potential $S$ of $ub_{s,\tau}$ (as in (\ref{generalS})). We will
derive (Section \ref{pde.sec}) a partial differential equation for $S$ and
then analyze this equation by the Hamilton--Jacobi method. In this method
(Section \ref{hj.sec}), we compute $S$ along certain \textquotedblleft
characteristic curves\textquotedblright\ $\lambda(\tau)$ and $\varepsilon
(\tau),$ which depend on a choice of $\lambda_{0}$ and $\varepsilon_{0}.$

Since the definition of the Brown measure involves letting $\varepsilon$ tend
to zero, we will try to choose $\lambda_{0}$ and $\varepsilon_{0}$ so that
$\varepsilon(\tau)=0.$ The condition $\varepsilon(\tau)=0$ can be achieved by
taking $\varepsilon_{0}=0$---but with the restriction that $\lambda_{0}$
cannot be in the closed support $\operatorname{supp}(\mu_{s})$ of the measure
$\mu_{s}$ in (\ref{mu0law}). (This restriction arises because the initial
momenta in the problem become undefined if $\lambda_{0}\in\operatorname{supp}%
(\mu_{s}).$) In Section \ref{outside.sec}, we will see that \emph{the Brown
measure }$\mu_{s,\tau}$\emph{ is zero} in a neighborhood of any point of the
form $\lambda(\tau)$ computed with $\varepsilon_{0}=0$ and $\lambda_{0}%
\in\mathbb{C}\setminus\operatorname{supp}(\mu_{s}).$

Now, if we compute $\lambda(\tau)$ with $\varepsilon_{0}=0$ and $\lambda_{0}$
outside $\operatorname{supp}(\mu_{s}),$ we will find (Proposition
\ref{solutionse0.prop}) that%
\[
\lambda(\tau)=f_{s-\tau}(\chi_{s}(\lambda_{0})).
\]
Thus, by (\ref{stauTilde}), we find that $\Sigma_{s,\tau}$ is the interior of
the complement of the set of points of the form $\lambda(\tau)=f_{s-\tau}%
(\chi_{s}(\lambda_{0})).$ Since the Brown measure is zero near each point of
the form $\lambda(\tau),$ we see that the closed support of $\mu_{s,\tau}$ is
contained in $\tilde{\Sigma}_{s,\tau}=\overline{\Sigma}_{s,\tau}.$
\end{remark}

We now introduce \textquotedblleft twisted logarithmic
coordinates\textquotedblright\ $v$ and $\delta$ for a nonzero complex number
$\lambda,$ which are real numbers satisfying%
\[
\lambda=e^{v\tau}e^{i\delta}.
\]
We can compute $v$ and $\delta$ explicitly as
\begin{equation}
v=\frac{1}{\tau_{1}}\log\left\vert \lambda\right\vert ;\quad\delta=\arg
\lambda-\frac{\tau_{2}}{\tau_{1}}\log\left\vert \lambda\right\vert ,
\label{vAndDeltaDef}%
\end{equation}
where $\arg\lambda$ is the argument of $\lambda$ and where $\tau=\tau
_{1}+i\tau_{2}.$ If we fix $\delta$ and let $v$ vary, $\log\lambda$ varies
along a straight line in the $\tau$-direction in the plane, while $\lambda$
itself varies along an exponential spiral that intersects the unit circle at
angle $\delta.$

\begin{figure}[ptb]
\centering
\includegraphics[scale=0.6]{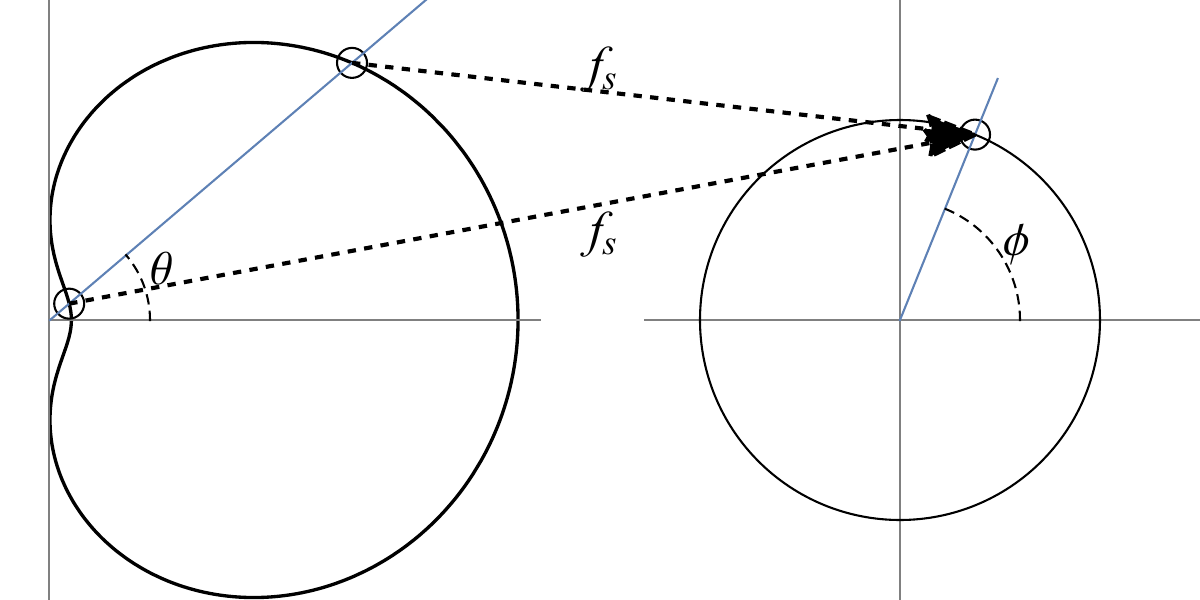}
\caption{The dependence of $\phi$ on $\theta.$ Shown for $s=2$ and $\mu
_{0}=\delta_{1}$.}
\label{thetaphi.fig}
\end{figure}

\begin{figure}[ptb]
\centering
\includegraphics[scale=0.6]{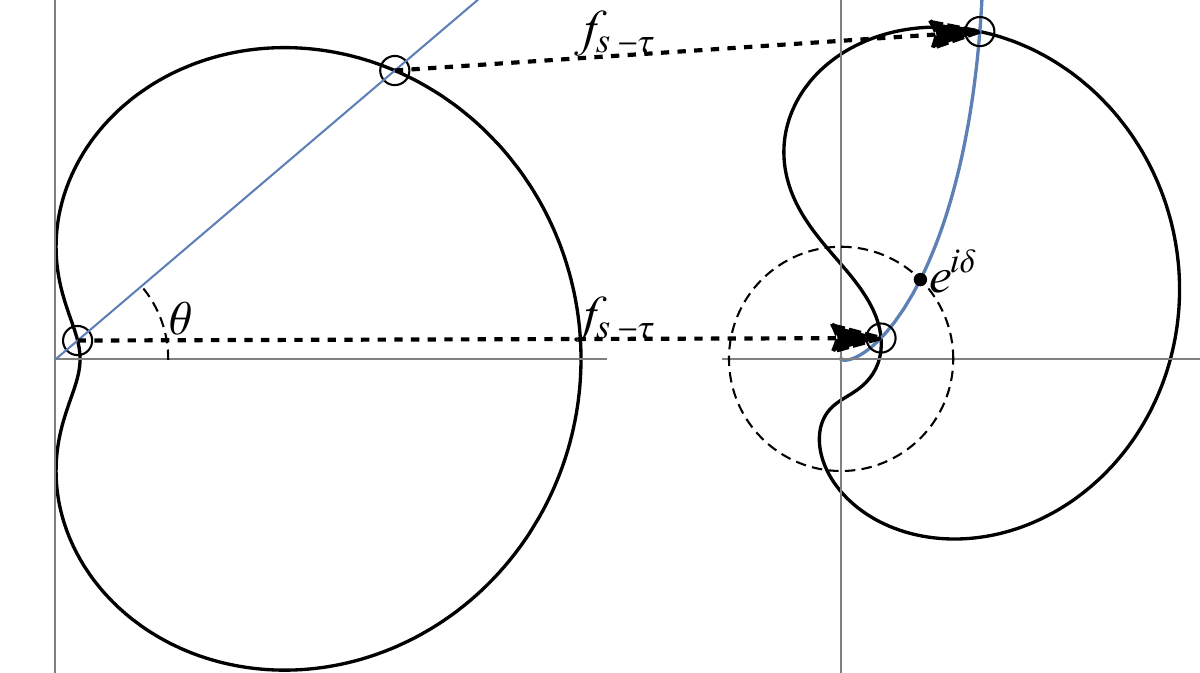}
\caption{The dependence of $\delta$ on $\theta.$ Shown for $s=2,$ $\tau
=\frac{3}{2}+\frac{1}{2}i,$ and $\mu_{0}=\delta_{1}$.}
\label{thetadelta.fig}
\end{figure}

We will show that for each $\delta,$ if an exponential spiral $\{e^{v\tau
}e^{i\delta}\}_{v\in\mathbb{R}}$ intersects $\Sigma_{s,\tau},$ this spiral
intersects the boundary of $\Sigma_{s,\tau}$ in exactly two points,
corresponding to $v=v_{1}(\delta)$ and $v=v_{2}(\delta),$ with $v_{1}%
(\delta)<v_{2}(\delta).$ For any two such boundary points, we will show that
there is a unique pair $r_{s}(\theta)e^{i\theta}$ and $r_{s}(\theta
)^{-1}e^{i\theta}$ of points on the boundary of $\Sigma_{s}$ such that
\begin{equation}
f_{s-\tau}(r_{s}(\theta)e^{i\theta})=e^{v_{1}(\delta)\tau}e^{i\delta};\quad
f_{s-\tau}(r_{s}(\theta)^{-1}e^{i\theta})=e^{v_{2}(\delta)\tau}e^{i\delta}.
\label{boundaryToBoundary}%
\end{equation}
Proposition \ref{SigmasChar.prop} then tells us that $f_{s}$ maps the points
$r_{s}(\theta)e^{i\theta}$ and $r_{s}(\theta)^{-1}e^{i\theta}$ to the same
point $e^{i\phi}$ on the unit circle. Thus, to each such $\delta,$ we
associate first an angle $\theta$ and then an angle $\phi.$ The relationship
among $\phi,$ $\theta$, and $\delta$ is shown in Figures \ref{thetaphi.fig}
and \ref{thetadelta.fig}.

\begin{notation}
\label{phistau.notation}Fix a $\delta$ such that the exponential spiral
$e^{v\tau}e^{i\delta}$ intersects $\Sigma_{s,\tau}.$ Let $\theta$ be the angle
such that (\ref{boundaryToBoundary}) holds. Then define $\phi^{s,\tau}%
(\delta)$ by%
\begin{align*}
\phi^{s,\tau}(\delta)  &  =\arg f_{s}(r_{s}(\theta)e^{i\theta})\\
&  =\arg f_{s}(r_{s}(\theta)^{-1}e^{i\theta}).
\end{align*}

\end{notation}

We are now ready to state our main result.

\begin{theorem}
\label{introMain.thm}All the mass of the Brown measure $\mu_{s,\tau}$ is
concentrated on the open set $\Sigma_{s,\tau}.$ The measure $\mu_{s,\tau}$ has
a density on $\Sigma_{s,\tau}$ given by%
\begin{equation}
d\mu_{s,\tau}(\lambda)=\frac{1}{2\pi\tau_{1}}\frac{1}{\left\vert
\lambda\right\vert ^{2}}\frac{d\phi^{s,\tau}(\delta)}{d\delta}~dx~dy,
\label{introMainFormula}%
\end{equation}
where $\lambda=x+iy.$ In logarithmic coordinates $\rho=\log\left\vert
\lambda\right\vert $ and $\theta=\arg\lambda,$ (\ref{introMainFormula})
becomes%
\[
d\mu_{s,\tau}(\lambda)=\frac{1}{2\pi\tau_{1}}\frac{d\phi^{s,\tau}(\delta
)}{d\delta}~d\rho~d\theta,
\]
so that the density in logarithmic coordinates is constant in the $\tau$-direction.
\end{theorem}

Here, when we say \textquotedblleft constant in the $\tau$%
-direction,\textquotedblright\ we are thinking of the complex number $\tau$ as
a vector in $\mathbb{C}\cong\mathbb{R}^{2}.$ See the right-hand side of Figure
\ref{brownwithlog.fig}. The density of $\mu_{s,\tau}$ in the twisted
logarithmic coordinates $v$ and $\delta$ differs from the density in ordinary
logarithmic coordinates by a constant, since $v$ and $\delta$ are linear
functions of $\rho$ and $\theta.$

Let us consider what happens to (\ref{introMainFormula}) in the case $\tau=s.$
In that case, $f_{s-\tau}(z)=f_{0}(z)=z,$ so that $\Sigma_{s,\tau}=\Sigma_{s}$
and in Figure \ref{thetadelta.fig}, $\delta$ is equal to $\theta.$ Thus,
writing $\theta$ in place of $\delta$ and $s$ in place of $\tau,$ we have%
\[
d\mu_{s,s}(\lambda)=\frac{1}{2\pi s}\frac{1}{\left\vert \lambda\right\vert
^{2}}\frac{d\phi^{s,s}(\theta)}{d\theta}~dx~dy.
\]
This formula agrees with the result previously obtained in Theorem 2.2 and
Proposition 8.5 of \cite{DHKBrown} (in the case $\mu_{0}=\delta_{1}$) and
Proposition 4.24 and Theorem 4.28 of \cite{HZ} (in general).

\section{The domains $\Sigma_{s}$ and $\Sigma_{s,\tau}$\label{domains.sec}}

\subsection{The domain $\Sigma_{s}$\label{domSigmas.sec}}

In this section, we study the domain $\Sigma_{s}$ in Definition
\ref{sigmas.def}.

As in \cite[Equation (4.1)]{DHKBrown} and \cite[Equation (4.32)]{HZ}, we
define a function $T:\mathbb{C}\rightarrow\mathbb{R}$ by%
\begin{equation}
T(\lambda)=\frac{\log(\left\vert \lambda\right\vert ^{2})}{\left\vert
\lambda\right\vert ^{2}-1}\left(  \int_{S^{1}}\frac{1}{\left\vert \lambda
-\xi\right\vert ^{2}}~d\mu_{0}(\xi)\right)  ^{-1}, \label{Tdef}%
\end{equation}
where we give $T$ the value 0 if the integral in the definition has the value
$+\infty.$ The function $T$ is positive for $\lambda$ outside the closed
support of the measure $\mu_{0}.$ If $\left\vert \lambda\right\vert =1,$ we
interpret $\log(\left\vert \lambda\right\vert ^{2})/(\left\vert \lambda
\right\vert ^{2}-1)$ as having the value 1 in accordance with the limit
$\lim_{x\rightarrow1}\log(x)/(x-1)=1.$ The motivation for this definition is
this: If $\left\vert z\right\vert \neq1,$ then for any real number $s,$ we can
easily compute that
\begin{equation}
\left\vert f_{s}(z)\right\vert =1\Longleftrightarrow T(z)=s.
\label{equivalence}%
\end{equation}
(In \cite{DHKBrown,HZ}, the function $T$ also gives the small-$\varepsilon
_{0}$ lifetime to a solution of a certain system of ODEs.) By \cite[Lemma
4.15]{HZ}, $T$ satisfies $T(re^{i\theta})=T(r^{-1}e^{i\theta})$ and
\thinspace$T(re^{i\theta})$ is a decreasing function of $r$ for $0<r<1$ with
$\theta$ fixed.

As in \cite[Equation (2.16)]{HZ}, we then define a function $r_{s}$ by
\[
r_{s}(\theta)=\left\{
\begin{array}
[c]{cl}%
\text{the unique }r<1\text{ such that }T(re^{i\theta})=s & \text{if such an
}r\text{ exists}\\
1 & \text{if no such }r\text{ exists}%
\end{array}
\right.  .
\]
In light of (\ref{equivalence}), we may equivalently say that $r_{s}(\theta)$
is the unique $r<1$ for which $\left\vert f_{s}(re^{i\theta})\right\vert =1,$
or 1 if no such $r$ exists. Since $\left\vert f_{s}(z)\right\vert =1$ when
$\left\vert z\right\vert =1,$ we see that
\begin{equation}
\left\vert f_{s}(r_{s}e^{i\theta})\right\vert =1 \label{frsOne}%
\end{equation}
for all $\theta,$ whether $r_{s}(\theta)<1$ or $r_{s}(\theta)=1.$ The function
$r_{s}$ is continuous everywhere and analytic when $r_{s}(\theta)<1$
\cite[Remark 2.8]{HZ}. As a consequence of (\ref{fsZbar}) and (\ref{frsOne}),
we also have the relation%
\begin{equation}
f_{s}(r_{s}(\theta)e^{i\theta})=f_{s}(r_{s}(\theta)^{-1}e^{i\theta}).
\label{frsID}%
\end{equation}

The following result is a strengthening of Proposition \ref{SigmasChar.prop}.

\begin{proposition}
\label{moreSigmasChars.prop}The domain $\Sigma_{s}$ in Definition
\ref{sigmas.def} can be computed as%
\begin{equation}
\Sigma_{s}=\left\{  re^{i\theta}\left\vert r_{s}(\theta)<r<\frac{1}%
{r_{s}(\theta)}\right.  \right\}  \label{rsChar}%
\end{equation}
and also as%
\begin{equation}
\Sigma_{s}=\left\{  \left.  z\in\mathbb{C}\right\vert T(z)<s\right\}  .
\label{TlessS}%
\end{equation}
Finally, if $z$ is outside of $\overline{\Sigma}_{s}$ then $T(z)>s$ and $z$ is
outside the closed support of $\mu_{0}.$
\end{proposition}

\begin{proof}
The set $\mathbb{C}\setminus\operatorname{supp}(\mu_{s})$ consists of three
parts: the set $\left\{  \left\vert z\right\vert <1\right\}  ,$ the set
$\left\{  \left\vert z\right\vert >1\right\}  ,$ and $S^{1}\setminus
\operatorname{supp}(\mu_{s}).$ According to \cite[Theorem 3.2]{Zhong}, we have%
\[
\chi_{s}(\{\left\vert z\right\vert <1\})=\{\left.  re^{i\theta}\right\vert
r<r_{s}(\theta)\}.
\]
Then by (\ref{chisZbar}),%
\[
\chi_{s}(\left\{  \left\vert z\right\vert >1\right\}  )=\left\{  \left.
re^{i\theta}\right\vert r>r_{s}(\theta)^{-1}\right\}  .
\]
Then by Proposition 3.7 and Theorem 3.8 of \cite{Zhong}, the image of
$S^{1}\setminus\operatorname{supp}(\mu_{s})$ is the set of points $e^{i\theta
}$ in the unit circle for which $r_{s}\equiv1$ near $\theta.$ Thus, in light
of Definition \ref{sigmas.def}, $\tilde{\Sigma}_{s}$ is the complement of the
union of these three images. It is easy to see that $\tilde{\Sigma}_{s}$ is
the closure of the set on the right-hand side of (\ref{rsChar}). Then
$\Sigma_{s},$ defined as the interior of $\tilde{\Sigma}_{s},$ is simply the
set on the right-hand side of (\ref{rsChar}), so that $\tilde{\Sigma}%
_{s}=\overline{\Sigma}_{s}.$ Theorem 4.10 of \cite{HZ} then gives the
alternative characterization of $\Sigma_{s}$ in (\ref{TlessS}) and shows that
$T>s$ outside the closure of $\Sigma_{s}.$ That points outside the closure of
$\Sigma_{s}$ are outside the closed support of $\mu_{0}$ follows by the same
argument as in the proof of Lemma 6.3 of \cite{HHadditive}.
\end{proof}

We will use the following notation.

\begin{definition}
\label{phiS.def}We define a continuous function $\phi^{s}$ from $S^{1}$ to
$S^{1}$ by the relation%
\[
e^{i\phi^{s}(\theta)}=f_{s}(r_{s}(\theta)e^{i\theta})=f_{s}(r_{s}(\theta
)^{-1}e^{i\theta}).
\]

\end{definition}

Here we have used the identity (\ref{frsID}). The notation is illustrated in
Figure \ref{thetaphi.fig}. According to \cite[Proposition 3.7]{Zhong}, the map
$\phi_{s}$ defines a homeomorphism of the unit circle to itself.

\subsection{The domain $\Sigma_{s,\tau}$}

We study the domain $\Sigma_{s,\tau}$ in Definition \ref{def:Sigmastaudef}. We
will use the twisted logarithmic coordinates introduced in (\ref{coordinates}%
). Recall that the $(v,\delta)$-coordinates of a nonzero complex number $z$
may be obtained by writing $z$ as%
\[
z=e^{v\tau}e^{i\delta},
\]
with $v$ and $\delta$ in $\mathbb{R}.$

We let $J$ denote the Herglotz integral of $\mu_{0}$, given by
\begin{equation}
J(z)=\int_{S^{1}}\frac{\xi+z}{\xi-z}\,d\mu_{0}(\xi), \label{herglotzJ}%
\end{equation}
which appears in Definition \ref{fFunct.def} defining the function $f_{\beta
}.$ Using the Cauchy--Schwarz inequality, we compute that%
\begin{align*}
\left\vert J(z)\right\vert  &  =\left\vert 1+2z\int_{S^{1}}\frac{1}{\xi
-z}~d\mu_{0}(\xi)\right\vert \\
&  \leq1+2\left\vert z\right\vert \left(  \int_{S^{1}}\frac{1}{\left\vert
\xi-z\right\vert ^{2}}~d\mu_{0}(\xi)\right)  ^{1/2}.
\end{align*}
Thus, $J(z)$ is defined and finite whenever the integral in the definition
(\ref{Tdef}) of the function $T$ is finite.

If there is some $s>0$ for which $z$ is outside $\Sigma_{s},$ then by
(\ref{TlessS}), $T(z)\neq0$ and $J(z)$ will be defined and finite.
Furthermore,
\begin{align}
&  \left\vert J(z_{1})-J(z_{2})\right\vert \nonumber\\
&  =2\left\vert z_{1}-z_{2}\right\vert \left\vert \int_{S^{1}}\frac{\xi}%
{(\xi-z_{1})(\xi-z_{2})}~d\mu_{0}(\xi)\right\vert \nonumber\\
&  \leq2\left\vert z_{1}-z_{2}\right\vert \left(  \int_{S^{1}}\frac
{1}{\left\vert \xi-z_{1}\right\vert ^{2}}~d\mu_{0}(\xi)\int_{S^{1}}\frac
{1}{\left\vert \xi-z_{2}\right\vert ^{2}}~d\mu_{0}(\xi)\right)  ^{1/2}%
.\label{Jdiff}%
\end{align}
From this and (\ref{TlessS}), we can easily verify that $J$ is continuous on
the complement of $\Sigma_{s},$ for every $s>0.$

We will prove the following results. One of the results concerns the $\delta
$-coordinate of $f_{s-\tau}(r_{s}(\theta)e^{i\theta})$; the other gives a
characterization of $\Sigma_{s,\tau}$.

\begin{theorem}
\label{thm:Sigmastau} Fix a positive real number $s$ and a nonzero complex
number $\tau$ satisfying $\left\vert \tau-s\right\vert \leq s.$

\begin{enumerate}
\item \label{DomainDelta.point}Define a function $\delta^{s,\tau}%
:\mathbb{R}\rightarrow\mathbb{R}$ by
\begin{equation}
\delta^{s,\tau}(\theta)=\theta+\frac{s-\left\vert \tau\right\vert ^{2}%
/\tau_{1}}{2}\operatorname{Im}\left[  \int_{S^{1}}\left(  \frac{\xi
+r_{s}(\theta)e^{i\theta}}{\xi-r_{s}(\theta)e^{i\theta}}\right)  \,d\mu
_{0}(\xi)\right]  . \label{eq:deltadef}%
\end{equation}
Then $\delta^{s,\tau}$ is continuous and strictly increasing and satisfies
\[
\delta^{s,\tau}(\theta+2\pi)=\delta^{s,\tau}(\theta)+2\pi.
\]
Furthermore, the $\delta$-coordinate of $f_{s-\tau}(r_{s}(\theta)e^{i\theta})$
may be computed as%
\[
\delta(f_{s-\tau}(r_{s}(\theta)e^{i\theta}))=\delta^{s,\tau}(\theta).
\]

\item \label{Domainu1u2.point}For each $\delta\in\mathbb{R},$ choose
$\theta\in\mathbb{R}$ so that $\delta^{s,\tau}(\theta)=\delta,$ which is
possible by Point 1. Then define $v_{1}^{s,\tau}(\delta)$ and $v_{2}^{s,\tau
}(\delta)$ to be the $v$-coordinates of $f_{s-\tau}(r_{s}(\theta)e^{i\theta})$
and $f_{s-\tau}(r_{s}(\theta)^{-1}e^{i\theta}),$ respectively, so that
\begin{equation}
f_{s-\tau}(r_{s}(\theta)e^{i\theta})=e^{\tau v_{1}^{s,\tau}(\delta)}%
e^{i\delta};\quad f_{s-\tau}(r_{s}(\theta)^{-1}e^{i\theta})=e^{\tau
v_{2}^{s,\tau}(\delta)}e^{i\delta}. \label{eq:u1u2}%
\end{equation}
Then $v_{1}^{s,\tau}(\delta)\leq v_{2}^{s,\tau}(\delta)$ for all $\delta$ and
the domain $\Sigma_{s,\tau}$ may be computed as
\begin{equation}
\{\left.  e^{v\tau}e^{i\delta}\right\vert v_{1}^{s,\tau}(\delta)<v<v_{2}%
^{s,\tau}(\delta)\} \label{u1u2Char}%
\end{equation}
and $\tilde{\Sigma}_{s,\tau}=\overline{\Sigma}_{s,\tau}.$
\end{enumerate}
\end{theorem}

Note that if $v_{1}^{s,\tau}(\delta)=v_{2}^{s,\tau}(\delta),$ then there are
no points $\Sigma_{s,\tau}$ of the form $e^{v\tau}e^{i\delta}.$%

\begin{figure}[ptb]
\centering
\includegraphics[scale=0.6]{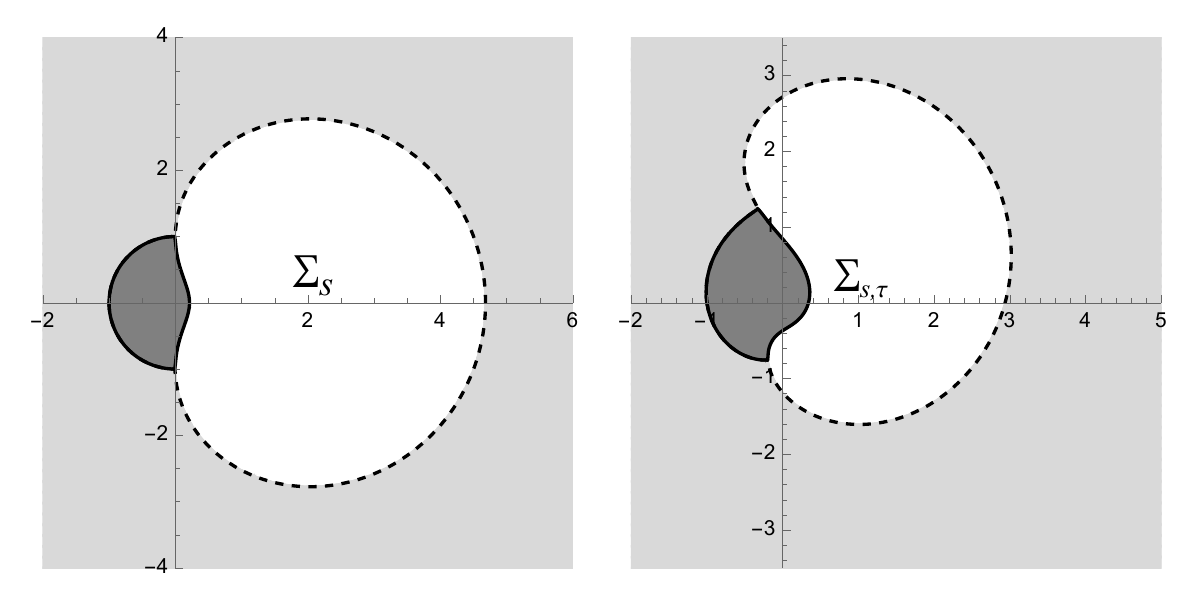}
\caption{The domains $\Sigma_{s}$ (left) and $\Sigma_{s,\tau}$ (right). The
complement of each domain in the Riemann sphere consists of two topological
disks (dark and light gray), which may overlap on the boundary. The map
$f_{s-\tau}$ takes each disk on the left homeomorphically to the corresponding
disk on the right. Shown for $\mu_{0}=\delta_{1},$ $s=2,$ and $\tau=\frac
{3}{2}+\frac{1}{2}i$. }
\label{twodomains.fig}
\end{figure}

Point \ref{DomainDelta.point} of the theorem is proved in Section
\ref{sect:fInjBdry}. Section \ref{fInjOutside.sec} then shows that $f_{s-\tau
}$ is injective on the complement of $\Sigma_{s}.$ The characterization of the
domain in (\ref{u1u2Char}) is then given in Section \ref{sect:charSigmastau}.
We briefly outline the argument for this characterization. The complement of
$\Sigma_{s}$ in the Riemann sphere can be computed as the union of two closed
topological disks, $\{re^{i\theta}|\,r\leq r_{s}(\theta)\}$ and $\{re^{i\theta
}|\,r\geq r_{s}(\theta)^{-1}\}\cup\{\infty\},$ which intersect at the points
(if any) where $r_{s}(\theta)=1.$ These two disks are shown in dark gray and
light gray, respectively, on the left-hand side of Figure \ref{twodomains.fig}%
. We will show that $f_{s-\tau}$ maps these disks to two closed topological
disks, namely%
\begin{equation}
\{e^{\tau v}e^{i\delta}|v\leq v_{1}^{s,\tau}(\delta)\}\cup\{0\} \label{disk1}%
\end{equation}
and%
\begin{equation}
\{e^{\tau v}e^{i\delta}|v\geq v_{2}^{s,\tau}(\delta)\}\cup\{\infty\}.
\label{disk2}%
\end{equation}
These disks are shown in dark gray and light gray, respectively, on the
right-hand side of Figure \ref{twodomains.fig}. Since the disks (\ref{disk1})
and (\ref{disk2}) cover $v\leq v_{1}^{s,\tau}(\delta)$ and $v\geq
v_{2}^{s,\tau}(\delta),$ the complement of their union is the set in
(\ref{u1u2Char}) where $v_{1}^{s,\tau}(\delta)<v<v_{2}^{s,\tau}(\delta).$ Once
this result is established, we can easily compute $\tilde{\Sigma}_{s,\tau
}:=[f_{s-\tau}((\overline{\Sigma}_{s})^{c})]^{c}$ and verify that
$\Sigma_{s,\tau}:=\operatorname{int}(\tilde{\Sigma}_{s})$ is the set in
(\ref{u1u2Char}).

\subsection{Behavior of $f_{s-\tau}$ on the boundary of $\Sigma_{s}%
$\label{sect:fInjBdry}}

Throughout the section, we assume that $s$ is a positive real number and
$\tau$ is a nonzero complex number satisfying $\left\vert \tau-s\right\vert
\leq s.$ As usual, we let $\tau_{1}$ and $\tau_{2}$ denote the real and
imaginary parts of $\tau,$ respectively.

\subsubsection{Proof of Point 1 of Theorem \ref{thm:Sigmastau}.}

We now work our way toward the proof of the first part of Theorem
\ref{thm:Sigmastau}. We define $R_{s}(\theta)$ and $I_{s}(\theta)$ as the real
and imaginary parts of the Herglotz function $J$ in (\ref{herglotzJ}), on the
curve $r_{s}(\theta)e^{i\theta}$:%
\begin{equation}
R_{s}(\theta)=\operatorname{Re}[J(r_{s}(\theta)e^{i\theta})];\quad
I_{s}(\theta)=\operatorname{Im}[J(r_{s}(\theta)e^{i\theta})]. \label{RsIs}%
\end{equation}

\begin{lemma}
\label{lem:spiralBoundary} The $\delta$-coordinate of the point $f_{s-\tau
}(r_{s}(\theta)e^{i\theta})$ may be computed as in (\ref{eq:deltadef}) and we
have the inequality
\begin{equation}
0<\frac{|\tau|^{2}}{s\tau_{1}}\leq2. \label{tauSquaredIneq}%
\end{equation}

\end{lemma}

\begin{proof}
We may easily compute that $f_{s-\tau}(z)=f_{s}(z)e^{-\frac{\tau}{2}J(z)}.$
Since $\left\vert f_{s}(r_{s}(\theta)e^{i\theta})\right\vert $ equals 1, we
must have%
\begin{equation}
f_{s}(r_{s}(\theta)e^{i\theta})=e^{i\left(  \theta+\frac{s}{2}I_{s}%
(\theta)\right)  }. \label{IsAndPhi}%
\end{equation}
Thus,%
\[
f_{s-\tau}(r_{s}(\theta)e^{i\theta})=\exp\left\{  i\theta+i\frac{s}{2}%
I_{s}(\theta)\right\}  \exp\left\{  -\frac{\tau}{2}R_{s}(\theta)\right\}
\exp\left\{  -\frac{i\tau}{2}I_{s}(\theta)\right\}  .
\]
Now, multiplying a complex number $z$ by $e^{i\alpha},$ $\alpha\in\mathbb{R},$
just adds $\alpha$ to the $\delta$-coordinate of $z$, while multiplying $z$ by
$e^{v\tau},$ $v\in\mathbb{R},$ does not change the $\delta$-coordinate of $z.$
Thus,%
\begin{align*}
\delta\left(  f_{s-\tau}(r_{s}(\theta)e^{i\theta})\right)   &  =\theta
+\frac{s}{2}I_{s}(\theta)+\delta\left(  \exp\left\{  -\frac{i\tau}{2}%
I_{s}(\theta)\right\}  \right) \\
&  =\theta+\frac{s}{2}I_{s}(\theta)-\frac{\left\vert \tau\right\vert ^{2}%
}{2\tau_{1}}I_{s}(\theta),
\end{align*}
where the second equality is a direct computation using the definition
(\ref{coordinates}) of $\delta.$ This expression is the claimed result
(\ref{eq:deltadef}).

Since $\tau$ is nonzero and satisfies $\left\vert \tau-s\right\vert \leq s,$
both $\left\vert \tau\right\vert ^{2}$ and $\tau_{1}$ are positive, so that
$0<\left\vert \tau\right\vert ^{2}/\tau_{1}.$ The inequality $\left\vert
\tau\right\vert ^{2}/(2\tau_{1})\leq s$ follows from the assumption
$\left\vert \tau-s\right\vert ^{2}\leq s^{2},$ after writing $\left\vert
\tau-s\right\vert ^{2}$ as $\left\vert \tau\right\vert ^{2}-2s\tau_{1}+s^{2}.$
\end{proof}

Next, we estimate the derivative of $\delta^{s,\tau}(\theta)$ at points
$\theta$ where $r_{s}(\theta)<1$.

\begin{proposition}
\label{prop:dpsidtheta} The function $\delta^{s,\tau}$ in (\ref{eq:deltadef})
is differentiable at every point $\theta$ for which $r_{s}(\theta)<1,$ and at
every such point, we have
\[
0<\frac{d\delta^{s,\tau}}{d\theta}<2.
\]

\end{proposition}

\begin{proof}
If $I_{s}$ is as in (\ref{RsIs}), then (\ref{eq:deltadef}) says that
\begin{equation}
\delta^{s,\tau}(\theta)=\theta+\frac{s-\left\vert \tau\right\vert ^{2}%
/\tau_{1}}{2}I_{s}(\theta). \label{deltaFromIs}%
\end{equation}
Meanwhile, the relation (\ref{IsAndPhi}) shows that the function $\phi^{s}$ in
(\ref{phiS.def}) can be computed as
\begin{equation}
\phi^{s}(\theta)=\theta+\frac{s}{2}I_{s}(\theta). \label{phisFormula}%
\end{equation}
If we solve for $I_{s}$ in (\ref{phisFormula}), substitute into
(\ref{deltaFromIs}), and differentiate, we get
\[
\frac{d\delta^{s,\tau}(\theta)}{d\theta}=\frac{\left\vert \tau\right\vert
^{2}}{s\tau_{1}}+\left(  1-\frac{\left\vert \tau\right\vert ^{2}}{s\tau_{1}%
}\right)  \frac{d\phi^{s}}{d\theta}.
\]

Now, \cite[Lemma 4.20]{HZ} shows that for every $\theta$ with $r_{s}%
(\theta)<1,$ the function $\phi^{s}$ is differentiable at $\theta$ with the
estimate $0<d\phi^{s}/d\theta<2.$ Then if $\left\vert \tau\right\vert
^{2}/(s\tau_{1})<1,$ we get%
\[
\frac{\left\vert \tau\right\vert ^{2}}{s\tau_{1}}<\frac{d\delta^{s,\tau}%
}{d\theta}<2-\frac{\left\vert \tau\right\vert ^{2}}{s\tau_{1}},
\]
and if $\left\vert \tau\right\vert ^{2}/(s\tau_{1})>1,$ we get%
\[
2-\frac{\left\vert \tau\right\vert ^{2}}{s\tau_{1}}<\frac{d\delta^{s,\tau}%
}{d\theta}<\frac{\left\vert \tau\right\vert ^{2}}{s\tau_{1}},
\]
while if $\left\vert \tau\right\vert ^{2}/(s\tau_{1})=1,$ we get%
\[
\frac{d\delta^{s,\tau}(\theta)}{d\theta}=\frac{\left\vert \tau\right\vert
^{2}}{s\tau_{1}}=1.
\]
In all cases, the inequality (\ref{tauSquaredIneq}) shows that the lower bound
on $d\delta^{s,\tau}/d\theta$ is greater than or equal to zero and the upper
bound is less than or equal to 2.
\end{proof}

\begin{proposition}
\label{prop:psiincreasing} The function $\delta^{s,\tau}$ is strictly
increasing on $\{\theta\in\mathbb{R}\vert r_{s}(\theta)=1\}$. That is, if
$\theta_{2}>\theta_{1}$ and $r_{s}(\theta_{1})=r_{s}(\theta_{2})=1$, then
$\delta^{s,\tau}(\theta_{2})>\delta^{s,\tau}(\theta_{1})$.
\end{proposition}

\begin{proof}
Denote $\left\vert \tau\right\vert ^{2}/\tau_{1}$ by $t$ and let $\theta
_{2}>\theta_{1}$. Note that $(\xi+e^{i\theta_{j}})/(\xi-e^{i\theta_{j}})$,
$j=1,2$, are purely imaginary. If $r_{s}(\theta_{j})=1,$ we compute that
\begin{align}
&  \delta^{s,\tau}(\theta_{2})-\delta^{s,\tau}(\theta_{1})\nonumber\\
=  &  \left(  \theta_{2}-\frac{s-t}{2}I_{s}(\theta)\right)  -\left(
\theta_{1}-\frac{s-t}{2}I_{s}(\theta)\right) \nonumber\\
=  &  (\theta_{2}-\theta_{1})\left(  1-(s-t)\frac{e^{i\theta_{2}}%
-e^{i\theta_{1}}}{\theta_{2}-\theta_{1}}\int_{S^{1}}\frac{\xi\,d\mu_{0}(\xi
)}{(\xi-e^{i\theta_{2}})(\xi-e^{i\theta_{1}})}\right)  . \label{delta12}%
\end{align}
Since $\theta_{1}\neq\theta_{2}$, we have
\begin{equation}
\left\vert \frac{e^{i\theta_{2}}-e^{i\theta_{1}}}{\theta_{2}-\theta_{1}%
}\right\vert =\left\vert \frac{\sin[(\theta_{2}-\theta_{1})/2]}{(\theta
_{2}-\theta_{1})/2}\right\vert <1. \label{eq:strictineq}%
\end{equation}

Now, since $r_{s}(\theta_{j})=1,$ $j=1,2,$ Proposition \ref{SigmasChar.prop}
tells us that the point $e^{i\theta_{j}}$ is not in $\Sigma_{s}.$ Thus, by
Proposition \ref{moreSigmasChars.prop}, $T(e^{i\theta_{j}})\geq s,$ or
\[
\int_{S^{1}}\frac{1}{\left\vert \xi-e^{i\theta_{j}}\right\vert ^{2}}\,d\mu
_{0}(\xi)\leq\frac{1}{s}\frac{\log r_{s}(\theta_{j})^{2}}{r_{s}(\theta
_{j})^{2}-1}=\frac{1}{s}.
\]
We thus obtain, using the Cauchy--Schwarz inequality, that
\begin{equation}
\left\vert \int_{S^{1}}\frac{\xi\,d\mu_{0}(\xi)}{(\xi-e^{i\theta_{2}}%
)(\xi-e^{i\theta_{1}})}\right\vert \leq\frac{1}{s}. \label{CSs}%
\end{equation}
By (\ref{tauSquaredIneq}), we have $0<t\leq2s$ or $-s<t-s<s.$ It then follows
from (\ref{eq:strictineq}) and (\ref{CSs}) that
\[
\left\vert (s-t)\frac{e^{i\theta_{2}}-e^{i\theta_{1}}}{\theta_{2}-\theta_{1}%
}\int_{S^{1}}\frac{\xi\,d\mu_{0}(\xi)}{(\xi-e^{i\theta_{2}})(\xi
-e^{i\theta_{1}})}\right\vert <1.
\]
Plugging this estimate into (\ref{delta12}) shows that $\delta^{s,\tau}%
(\theta_{2})-\delta^{s,\tau}(\theta_{1})>0.$
\end{proof}

\begin{proof}
[Proof of Point 1 of Theorem \ref{thm:Sigmastau}]We have already shown (Lemma
\ref{lem:spiralBoundary}) that $\delta^{s,\tau}$ is a local continuous version
of the $\delta$-coordinate of $f_{s-\tau}(r_{s}(\theta)e^{i\theta})$ for
$\theta$ in a neighborhood of $\theta_{0}$. By (\ref{Jdiff}), $f_{s-\tau}$ is
continuous on $\Sigma_{s}^{c},$ so $\delta^{s,\tau}$ is continuous.

Next, we show that $\delta^{s,\tau}$ is strictly increasing on $\mathbb{R}$.
Let $\theta_{1}<\theta_{2}$. We consider four cases, corresponding to whether
$r_{s}(\theta_{1})$ and $r_{s}(\theta_{2})$ are $1$ or less than $1$. If both
$r_{s}(\theta_{1})$ and $r_{s}(\theta_{2})$ are both $1$,
Proposition~\ref{prop:psiincreasing} shows that $\delta^{s,\tau}(\theta
_{1})<\delta^{s,\tau}(\theta_{2})$. If $r_{s}(\theta_{1})=1$ but $r_{s}%
(\theta_{2})<1$, then let $\alpha$ be the infimum of the interval $I$ around
$\theta_{2}$ on which $r_{s}$ is positive, so that $r_{s}(\alpha)=1$ and
$\theta_{1}\leq\alpha$. Then $\delta^{s,\tau}(\theta_{1})\leq\delta^{s,\tau
}(\alpha)$ by Proposition~\ref{prop:psiincreasing} and $\delta^{s,\tau}%
(\alpha)<\delta^{s,\tau}(\theta_{2})$ by the positivity of $d\delta/d\theta$
on $I$ in Proposition~\ref{prop:dpsidtheta}. The remaining cases are similar;
the case where both $r_{s}(\theta_{1})$ and $r_{s}(\theta_{2})$ are less than
$1$ can be subdivided into two cases depending on whether or not $\theta_{1}$
and $\theta_{2}$ are in the same interval where $r_{s}$ is less than $1$.

Meanwhile, from (\ref{eq:deltadef}), we can see that $\delta^{s,\tau}%
(\theta)-\theta$ is $2\pi$-periodic and continuous on $\mathbb{R}$. Since also
$\delta^{s,\tau}$ is strictly increasing, it defines a homeomorphism of
$S^{1}$ to $S^{1}.$
\end{proof}

\subsubsection{The functions $v_{1}$ and $v_{2}$}

In the following lemma, we show that the function $f_{s-\tau}$ maps the points
$r_{s}(\theta)e^{i\theta}$ and $r_{s}(\theta)^{-1}e^{i\theta}$ to a pair of
points lying on the same exponential spiral; the point at which this spiral
crosses the unit circle is $e^{i\delta^{s,\tau}(\theta)}$. Recall the
definition of the functions $v_{1}^{s,\tau}$ and $v_{2}^{s,\tau}$ in
(\ref{eq:u1u2}).

\begin{lemma}
\label{lem:u1u2} The points $f_{s-\tau}(r_{s}(\theta)e^{i\theta})$ and
$f_{s-\tau}(r_{s}(\theta)^{-1}e^{i\theta})$ have the same $\delta$-coordinate,
which we denote as $\delta^{s,\tau}(\theta).$ Furthermore, we have that
$v_{1}^{s,\tau}(\delta^{s,\tau}(\theta))<v_{2}^{s,\tau}(\delta^{s,\tau}%
(\theta))$ if and only if $r_{s}(\theta)<1.$ In fact,%
\begin{equation}
v_{2}^{s,\tau}(\delta^{s,\tau}(\theta))-v_{1}^{s,\tau}(\delta^{s,\tau}%
(\theta))=-\frac{2}{s}\log r_{s}(\theta). \label{uMinusU}%
\end{equation}

\end{lemma}

\begin{proof}
The Herglotz integral $J(z)$ in (\ref{herglotzJ}) is easily seen to satisfy
\begin{equation}
J(1/\bar{z})=-\overline{J(z)}, \label{J1zbar}%
\end{equation}
so that%
\begin{align*}
f_{s-\tau}(z)  &  =f_{s}(z)e^{-\frac{\tau}{2}J(z)}\\
f_{s-\tau}(1/\bar{z})  &  =f_{s}(1/\bar{z})e^{-\frac{\tau}{2}J(1/\bar{z}%
)}=\left(  \overline{f_{s}(z)}\right)  ^{-1}e^{\frac{\tau}{2}\overline{J(z)}}.
\end{align*}
If we evaluate at $z=r_{s}(\theta)e^{i\theta}$ where $\left\vert
f_{s}(z)\right\vert =1$ and use the notation $R_{s}(\theta)$ and $I_{s}%
(\theta)$ defined in (\ref{RsIs}), we obtain
\begin{align*}
f_{s-\tau}(r_{s}(\theta)e^{i\theta})  &  =f_{s}(r_{s}(\theta)e^{i\theta}%
)\exp\left\{  -\frac{\tau}{2}R_{s}(\theta)\right\}  \exp\left\{  -i\frac{\tau
}{2}I_{s}(\theta)\right\} \\
f_{s-\tau}(r_{s}(\theta)^{-1}e^{i\theta})  &  =f_{s}(r_{s}(\theta)e^{i\theta
})\exp\left\{  \frac{\tau}{2}R_{s}(\theta)\right\}  \exp\left\{  -i\frac{\tau
}{2}I_{s}(\theta)\right\}  .
\end{align*}
Since multiplying a complex number by $e^{v\tau}$ for $v\in\mathbb{R}$ does
not change the $\delta$-coordinate, $f_{s-\tau}(r_{s}(\theta)e^{i\theta})$ and
$f_{s-\tau}(r_{s}(\theta)^{-1}e^{i\theta})$ have the same $\delta$-coordinates.

Meanwhile, recalling that the $v$-coordinate of a point $\lambda$ is computed
as $\log\left\vert \lambda\right\vert /\tau_{1},$ we find that the difference
of the $v$-coordinates is $R_{s}(\theta)$. But using the identity $\left\vert
f_{s}(r_{s}(\theta)e^{i\theta}\right\vert =1,$ we may compute that%
\[
R_{s}(\theta)=-\frac{2}{s}\log r_{s}(\theta),
\]
giving the claimed formula (\ref{uMinusU}).
\end{proof}

\subsection{Injectivity of $f_{s-\tau}$ on the complement of $\Sigma_{s}%
$\label{fInjOutside.sec}}

In this section, we prove the following result. We continue to assume $s$ is a
positive real number and $\tau$ is a nonzero complex number with $\left\vert
\tau-s\right\vert \leq s.$

\begin{theorem}
\label{fInjective.thm} The map $f_{s-\tau}$ is injective on the complement of
$\Sigma_{s},$ including on the boundary.
\end{theorem}

We begin with the following lemma, whose proof is deferred until the end of
this section.

\begin{lemma}
\label{diffQuotient.lem}For all complex numbers $w_{1}$ and $w_{2},$ we have
\begin{equation}
\left\vert \frac{e^{w_{1}}-e^{w_{2}}}{w_{1}-w_{2}}\right\vert ^{2}\leq
\frac{(e^{2\operatorname{Re}w_{1}}-1)(e^{2\operatorname{Re}w_{2}}%
-1)}{(2\operatorname{Re}w_{1})(2\operatorname{Re}w_{2})}. \label{reverse}%
\end{equation}
Here all fractions are interpreted as having their limiting values when the
denominator is zero. Equality holds in (\ref{reverse}) only if $w_{2}=-\bar
{w}_{1}.$
\end{lemma}

\begin{proof}
[Proof of Theorem \ref{fInjective.thm}]Suppose, towards a contradiction, that
$f_{s-\tau}(z_{1})=f_{s-\tau}(z_{2})$ for two distinct points $z_{1}$ and
$z_{2}$ outside of $\Sigma_{s}$. Since $f_{s-\tau}(z)=0$ if and only if $z=0,$
we may assume $z_{1}$ and $z_{2}$ are nonzero. Using the definition of
$f_{s-\tau}$ and a bit of algebraic manipulation, we obtain
\[
\frac{z_{1}}{z_{2}}=\exp\left\{  \left(  s-\tau\right)  (z_{2}-z_{1}%
)\int_{S^{1}}\frac{\xi}{(\xi-z_{1})(\xi-z_{2})}~d\mu_{0}(\xi)\right\}  .
\]
We can then find some choice of a logarithm of $z_{1}$ and a logarithm of
$z_{2}$ such that taking the log of both sides gives
\begin{equation}
\frac{\log z_{1}-\log z_{2}}{z_{1}-z_{2}}=-(s-\tau)\int_{S^{1}}\frac{\xi}%
{(\xi-z_{1})(\xi-z_{2})}~d\mu_{0}(\xi). \label{firstLogs}%
\end{equation}

Applying the Cauchy--Schwarz inequality, we find that, since $\xi\in S^{1},$
we have
\begin{align}
&  \left\vert \int_{S^{1}}\frac{\xi}{(\xi-z_{1})(\xi-z_{2})}~d\mu_{0}%
(\xi)\right\vert \nonumber\\
&  \leq\left(  \int_{S^{1}}\frac{1}{\left\vert \xi-z_{1}\right\vert ^{2}}%
~d\mu_{0}(\xi)~\int_{S^{1}}\frac{1}{\left\vert \xi-z_{1}\right\vert ^{2}}%
~d\mu_{0}(\xi)\right)  ^{1/2}. \label{CSint}%
\end{align}
Now, since $z_{1}$ and $z_{2}$ are outside the domain $\Sigma_{s}$,
Proposition \ref{moreSigmasChars.prop} tells us that $T(z_{j})\geq s,$
$j=1,2$, and $T(z_{j})>s$ unless $z_{j}$ is in the boundary of $\Sigma_{s}.$
Thus,
\begin{equation}
\int_{S^{1}}\frac{1}{\left\vert \xi-z_{j}\right\vert ^{2}}~d\mu_{0}(\xi
)\leq\frac{1}{s}\frac{\log(\left\vert z_{j}\right\vert ^{2})}{\left\vert
z_{j}\right\vert ^{2}-1}, \label{IntZj}%
\end{equation}
and the inequality is strict unless $z_{j}$ is in the boundary of $\Sigma_{s}$.

Assume at least one of $z_{1}$ and $z_{2}$ is not in the boundary of
$\Sigma_{s}$ and use (\ref{CSint}) and (\ref{IntZj}) to get
\begin{equation}
\left\vert \int\frac{\xi}{(\xi-z_{1})(\xi-z_{2})}~d\mu_{0}(\xi)\right\vert
<\frac{2}{s}\left(  \frac{\log\left\vert z_{1}\right\vert \log\left\vert
z_{2}\right\vert }{(\left\vert z_{1}\right\vert ^{2}-1)(\left\vert
z_{2}\right\vert ^{2}-1)}\right)  ^{1/2}. \label{CSestimate}%
\end{equation}
By taking the absolute value of (\ref{firstLogs}) and using (\ref{CSestimate})
and the assumption $\left\vert \tau-s\right\vert \leq s,$ we get
\begin{equation}
\left\vert \frac{\log z_{1}-\log z_{2}}{z_{1}-z_{2}}\right\vert <2\left(
\frac{\log\left\vert z_{1}\right\vert \log\left\vert z_{2}\right\vert
}{(\left\vert z_{1}\right\vert ^{2}-1)(\left\vert z_{2}\right\vert ^{2}%
-1)}\right)  ^{1/2}. \label{inequality}%
\end{equation}
We now let $w_{1}=\log z_{1}$ and $w_{2}=\log z_{2},$ square both sides of
(\ref{inequality}), and take reciprocals, giving
\begin{equation}
\left\vert \frac{e^{w_{1}}-e^{w_{2}}}{w_{1}-w_{2}}\right\vert ^{2}%
>\frac{(e^{2\operatorname{Re}w_{1}}-1)(e^{2\operatorname{Re}w_{2}}%
-1)}{(2\operatorname{Re}w_{1})(2\operatorname{Re}w_{2})}. \label{inequality2}%
\end{equation}

What we have shown is that if $f_{s-\tau}(z_{1})=f_{s-\tau}(z_{2})$ for two
distinct points $z_{1}$ and $z_{2}$ outside $\Sigma_{s}$ with at least one of
$z_{1}$ and $z_{2}$ not in $\partial\Sigma_{s}$, then there exist \textit{some
choices} $w_{1}$ and $w_{2}$ of the logarithms of $z_{1}$ and $z_{2},$
respectively, for which (\ref{inequality2}) holds. But (\ref{inequality2})
contradicts Lemma \ref{diffQuotient.lem}.

It remains, then, to address the case in which both $z_{1}$ and $z_{2}$ are on
the boundary of $\Sigma_{s}.$ By Proposition \ref{moreSigmasChars.prop}, any
point $z$ on the boundary of $\Sigma_{s}$ must have the form $z=r_{s}%
(\theta)e^{i\theta}$ or $r_{s}(\theta)^{-1}e^{i\theta}.$ If $z_{1}$ and
$z_{2}$ have the same value of $\theta,$ then they must have the form
$z=r_{s}(\theta)e^{i\theta}$ and $r_{s}(\theta)^{-1}e^{i\theta}$ with
$r_{s}(\theta)<1,$ in which case Lemma \ref{lem:u1u2} tells us that the
$v$-coordinates of $f_{s-\tau}(z_{1})$ and $f_{s-\tau}(z_{2})$ are different.
If $z_{1}$ and $z_{2}$ have different values of $\theta,$ then Point 1 of
Theorem \ref{thm:Sigmastau} tells us that the $\delta$-coordinates of
$f_{s-\tau}(z_{1})$ and $f_{s-\tau}(z_{2})$ are different.
\end{proof}

We note that if $\tau=0,$ Lemma \ref{lem:u1u2} does not make sense as written,
because the coordinate $v$ is defined as $\log\left\vert \lambda\right\vert
/\tau_{1}.$ But if we rewrite the lemma using the coordinate $\rho
=\log\left\vert \lambda\right\vert ,$ we find that $\rho_{2}^{s,\tau}%
(\delta^{s,\tau}(\theta))=\rho_{1}^{s,\tau}(\delta^{s,\tau}(\theta))$ when
$\tau=0.$ Thus, when $\tau=0,$ injectivity of $f_{s-\tau}=f_{s}$ fails on the
boundary of $\Sigma_{s},$ as we have already observed in (\ref{frsID}).

We conclude this section by supplying the proof of Lemma
\ref{diffQuotient.lem}.

\begin{proof}
If we let $a=(w_{1}+w_{2})/2$ and $b=(w_{1}-w_{2})/2$, we can compute that
\begin{equation}
\left\vert \frac{e^{w_{1}}-e^{w_{2}}}{w_{1}-w_{2}}\right\vert ^{2}%
=e^{2\operatorname{Re}(a)}\left\vert \frac{\sinh(b)}{b}\right\vert ^{2}.
\label{diffQuotient}%
\end{equation}
We now claim that
\begin{equation}
\left\vert \frac{\sinh(b)}{b}\right\vert ^{2}\leq\left(  \frac{\sinh
(\operatorname{Re}b)}{\operatorname{Re}b}\right)  ^{2}, \label{sinhIneq}%
\end{equation}
with equality only if $b$ is real. To verify (\ref{sinhIneq}), we first
compute that with $b=x+iy,$ we have
\[
\left\vert \sinh(b)\right\vert ^{2}=\sinh^{2}(x)+\sin^{2}(y).
\]
We then use the elementary inequalities
\[
\sinh^{2}(x)\geq x^{2};\quad\sin^{2}(y)\leq y^{2},
\]
which hold with equality only at $x=0$ and $y=0,$ respectively. These
inequalities give
\begin{equation}
\left\vert \frac{\sinh(b)}{b}\right\vert ^{2}\leq\frac{\sinh^{2}(x)+y^{2}%
}{x^{2}+y^{2}}\leq\frac{\sinh^{2}(x)}{x^{2}}, \label{sinhIneq2}%
\end{equation}
with equality only if $x=0$ or $y=0.$ This result establishes the inequality
(\ref{sinhIneq}). If equality holds in (\ref{sinhIneq}), the equality holds in
(\ref{sinhIneq2}), so that either $y=0$---meaning that $b$ is real---or $x=0.$
But if $x=0,$ (\ref{sinhIneq}) says that $\left\vert \sin y/y\right\vert
\leq1,$ and equality holds there only if $y=0.$

Combining (\ref{diffQuotient}) and (\ref{sinhIneq}), we see that
\begin{align}
\left\vert \frac{e^{w_{1}}-e^{w_{2}}}{w_{1}-w_{2}}\right\vert ^{2}  &  \leq
e^{2\operatorname{Re}\,a}\left(  \frac{\sinh(\operatorname{\operatorname{Re}%
}\,b)}{\operatorname{Re}b}\right)  ^{2}\nonumber\\
&  =e^{\operatorname{\operatorname{Re}}w_{1}+\operatorname{Re}w_{2}}\left(
\frac{\sinh((\operatorname{Re}w_{1}-\operatorname{Re}w_{2})/2)}%
{(\operatorname{Re}w_{1}-\operatorname{Re}w_{2})/2}\right)  ^{2},
\label{stage}%
\end{align}
with equality only if $\operatorname{Im}\,b=0,$ that is, if $\operatorname{Im}%
w_{1}=\operatorname{Im}w_{2}.$ Now, the function $\log(\sinh(x)/x)$ is
strictly convex on the real line, which we can prove by computing its second
derivative as $1/x^{2}-1/\sinh^{2}x,$ which is positive for all $x.$ Thus,
\[
\frac{\sinh^{2}\left(  \frac{u+v}{2}\right)  }{\left(  \frac{u+v}{2}\right)
^{2}}\leq\frac{\sinh u}{u}\frac{\sinh v}{v}%
\]
with equality only when $u=v.$ Applying this with $u=\operatorname{Re}w_{1}$
and $v=-\operatorname{Re}w_{2}$ gives
\begin{equation}
\left(  \frac{\sinh((\operatorname{Re}w_{1}-\operatorname{Re}w_{2}%
)/2)}{(\operatorname{Re}w_{1}-\operatorname{\operatorname{Re}}w_{2}%
)/2}\right)  ^{2}\leq\frac{\sinh(\operatorname{Re}w_{1})}{\operatorname{Re}%
w_{1}}\frac{\sinh(\operatorname{Re}w_{2})}{\operatorname{Re}w_{2}},
\label{logConvex2}%
\end{equation}
with equality only if $\operatorname{Re}w_{2}=-\operatorname{Re}w_{1}.$
Plugging (\ref{logConvex2}) back into (\ref{stage}) and simplifying give the
claimed formula.
\end{proof}

\subsection{The characterization of $\Sigma_{s,\tau}$%
\label{sect:charSigmastau}}

Finally, we are ready to prove Point 2 of Theorem~\ref{thm:Sigmastau}, using
the strategy indicated in Figure \ref{twodomains.fig}. We use the following
elementary topological result.

\begin{lemma}
\label{disks.lem}Suppose $D_{1}$ and $D_{2}$ are closed topological disks in
$\mathbb{C}$ and $f$ is an injective, continuous map of $D_{1}$ into
$\mathbb{C}.$ Suppose $f$ maps at least one point in the interior of $D_{1}$
into the interior of $D_{2}$ and that $f$ maps the boundary of $D_{1}$
homeomorphically onto the boundary of $D_{2}.$ Then the image of $D_{1}$ is
contained in $D_{2}$ and $f$ is a homeomorphism of $D_{1}$ onto $D_{2}.$
\end{lemma}

\begin{proof}
Pick $z_{1}$ in the interior of $D_{1}$ mapping to the interior of $D_{2}.$ If
some $z_{2}$ in the interior of $D_{1}$ maps outside of $D_{2},$ connect
$z_{1}$ to $z_{2}$ by a continuous path in the interior of $D_{1}.$ The image
of this path then travels from the interior of $D_{2}$ to the complement of
$D_{2}$ and must cross $\partial D_{2}.$ Thus, there is some $z_{3}$ in the
interior of $D_{1}$ mapping to $\partial D_{2}.$ But $f$ maps $\partial D_{1}$
\textit{onto} $\partial D_{2},$ so there is $z_{4}$ in $\partial D_{1}$ with
$f(z_{4})=f(z_{3}),$ contradicting the injectivity of $f.$

Meanwhile, $f$ maps $\partial D_{1}$ to a curve that winds exactly once around
each point in the interior of $D_{2}.$ It then follows easily that $f$ must
map $D_{1}$ onto $D_{2}.$
\end{proof}

\begin{proof}
[Proof of Point 2 of Theorem~\ref{thm:Sigmastau}]The claim that $v_{1}%
(\delta)\leq v_{2}(\delta)$ follows from Lemma \ref{lem:u1u2}. To establish
the characterization of $\Sigma_{s,\tau},$ we apply the lemma with
$f=f_{s-\tau}$, with $D_{1}=\{\left.  re^{i\theta}\right\vert 0\leq
r<r_{s}(\theta)\},$ and with $D_{2}=\{e^{\tau v}e^{i\delta}|\,-\infty\leq
v\leq v_{1}^{s,\tau}(\delta)\}.$ (We define the $v$-coordinate of $0$ to be
$-\infty$.) Theorem \ref{fInjective.thm} tells us that $f_{s-\tau}$ is
injective on $D_{1}$. Furthermore, $f_{s-\tau}$ maps the point $0$ in $D_{1}$
to the point $0$ in $D_{2}.$ Finally, Point \ref{DomainDelta.point} of Theorem
\ref{thm:Sigmastau}, together with the definition (\ref{eq:u1u2}), tells us
that $f_{s-\tau}$ maps $\partial D_{1}$ homeomorphically onto $\partial
D_{2}.$ Thus, Lemma \ref{disks.lem} tells us that $f_{s-\tau}$ maps $D_{1}$
into $D_{2}$ and is a homeomorphism. A similar argument shows that $f_{s-\tau
}$ maps the disk%
\[
\{\left.  re^{i\theta}\right\vert r\geq r_{s}(\theta)^{-1}\}\cup\{\infty\}
\]
in the Riemann sphere onto the disk
\[
\{\left.  e^{v\tau}e^{i\delta}\right\vert v\geq v_{2}^{s,\tau}(\delta
))\}\cup\{\infty\}.
\]
We therefore conclude that
\[
f_{s-\tau}(\Sigma_{s}^{c})=\{e^{\tau v}e^{i\delta}|\,-\infty\leq v\leq
v_{1}^{s,\tau}(\delta)\text{ or }v\geq v_{2}^{s,\tau}(\delta)\}.
\]
Now, the boundary of $\Sigma_{s}$ consists of points of the form $r_{s}%
(\theta)e^{i\theta}$ and $r_{s}(\theta)^{-1}e^{i\theta}$ where $r_{s}%
(\theta)<1,$ together with limits of such points. By (\ref{eq:u1u2})\ and
Lemma \ref{lem:u1u2}, $f_{s-\tau}$ will map $\partial\Sigma_{s}$ to points
with $v$-coordinates $v_{1}^{s,\tau}(\delta)$ or $v_{2}^{s,\tau}(\delta)$
where $v_{1}^{s,\tau}(\delta)<v_{2}^{s,\tau}(\delta),$ together with limits of
such points. Then as in the proof of Proposition \ref{moreSigmasChars.prop},
we can easily see that $\tilde{\Sigma}_{s}$ is the closure of the set in
(\ref{u1u2Char}), so that $\Sigma_{s}$ is the set in (\ref{u1u2Char}) and
$\tilde{\Sigma}_{s}=\overline{\Sigma}_{s}.$
\end{proof}

\section{The PDE for $S$\label{pde.sec}}

\subsection{The main result}

Recall the definition of $b_{s,\tau}$ in (\ref{bstauDef}) and recall that $u$
is a unitary element freely independent of $b_{s,\tau}.$

\begin{definition}
\label{s.def} Define a function $S$ by
\begin{equation}
S(s,\tau,\lambda,\varepsilon)=\operatorname{tr}\left[  \log((ub_{s,\tau
}-\lambda)^{\ast}(ub_{s,\tau}-\lambda)+\varepsilon^{2})\right]
\label{Sfunction}%
\end{equation}
for all $s>0$, $\tau\in\mathbb{C}$ with $\left\vert \tau-s\right\vert \leq s$,
$\lambda\in\mathbb{C},$ and $\varepsilon>0.$
\end{definition}

We call $S$ the regularized log potential of $ub_{s,\tau}.$ Note that we
regularize the logarithm on the right-hand side of (\ref{Sfunction}) by
$\varepsilon^{2}.$ By contrast, \cite{DHKBrown}\ and \cite{HZ} use
$\varepsilon.$

\begin{theorem}
\label{thm:dSdtau}The function $S$ in Definition \ref{Sfunction} satisfies the
PDE
\begin{equation}
\frac{\partial S}{\partial\tau}=\frac{1}{8}\left[  1-\left(  1-\varepsilon
\frac{\partial S}{\partial\varepsilon}-2\lambda\frac{\partial S}%
{\partial\lambda}\right)  ^{2}\right]  . \label{thePDE}%
\end{equation}
for all $\lambda\in\mathbb{C},$ all $\varepsilon>0,$ and all $\tau>0$
satisfying $\left\vert \tau-s\right\vert <s,$ with the initial condition at
$\tau=0$ given by
\[
S(s,0,\lambda,\varepsilon)=\operatorname{tr}\left[  \log((uu_{s}%
-\lambda)^{\ast}(uu_{s}-\lambda)+\varepsilon^{2})\right]  ,
\]
where $u_{s}$ is the free unitary Brownian motion. More explicitly, we have%
\begin{equation}
S(s,0,\lambda,\varepsilon)=\int_{S^{1}}\log(\left\vert \xi-\lambda\right\vert
^{2}+\varepsilon^{2})~d\mu_{s}(\xi), \label{theInitialCondIntegral}%
\end{equation}
where $\mu_{s}$ is the law of the unitary element $uu_{s}.$
\end{theorem}

This result is proved in Section \ref{PDEpf}. Note that although in Definition
\ref{s.def} we allow $\left\vert \tau-s\right\vert \leq s,$ Theorem
\ref{thm:dSdtau} only asserts that the PDE (\ref{thePDE}) holds for $\tau$ in
the open set $\left\vert \tau-s\right\vert <s.$

Note that $S$ takes positive real values and therefore cannot be holomorphic
in $\tau$ or $\lambda$ (unless it is independent of that variable). Thus, the
derivative with respect to $\tau$ in (\ref{thePDE}) should be interpreted as
the usual complex partial derivative%
\[
\frac{\partial}{\partial\tau}=\frac{1}{2}\left(  \frac{\partial}{\partial
\tau_{1}}-i\frac{\partial}{\partial\tau_{2}}\right)  ,
\]
where $\tau=\tau_{1}+i\tau_{2},$ and similarly for $\partial/\partial\lambda.$

\subsection{A discussion of our approach\label{discussion.sec}}

Recall that the three-parameter free multiplicative Brownian motion
$b_{s,\tau}$ is defined as the solution $b_{s,\tau}(r)$ to a the free SDE
(\ref{bstauSDE}) evaluated at $r=1$ (Notation \ref{bstau.notation}). The
reader might then naturally expect that we would analyze the Brown measure of
$ub_{s,\tau}$ by deriving a PDE for the regularized logarithmic potential $S$
of $ub_{s,\tau}(r)$, with $r$ playing the role of the time variable. Indeed,
this is the approach taken in \cite{DHKBrown} and \cite{HZ} in the case
$\tau=s.$ Although we will, in fact, derive such a PDE---see (\ref{rPDE})
below---it is not the PDE we will actually use to compute the Brown measure.
Instead, to compute the Brown measure, we will use the PDE (\ref{thePDE}) in
Theorem \ref{thm:dSdtau}, in which we differentiate the regularized log
potential $S$ of $ub_{s,\tau}$ \textit{with respect to }$\tau$\textit{ while
keeping }$s$\textit{ fixed}.

This approach can be motivated by the work of the first author with Driver and
Kemp on the complex-time Segal--Bargmann transform \cite{DHKcomplex}. In that
paper, the space of $L^{2}$ functions on $U(N)$ with respect to a heat kernel
measure $\rho_{s}$ are transformed unitarily into holomorphic $L^{2}$
functions on $GL(N;\mathbb{C})$ with respect to a heat kernel measure
$\nu_{s,\tau}$ by means of the heat operator with complex time $\tau.$ The
measure $\nu_{s,\tau}$ is just the law of the finite-$N$ Brownian motion
$B_{s,\tau}^{N}(1)$ and the appearance of the heat operator with time $\tau$
makes it natural to differentiate in $\tau$ with $s$ fixed.

Regardless of that motivation, there are several advantages to differentiating
with respect to $\tau$ with $s$ fixed. First, the PDE\ in (\ref{thePDE}) is
simpler than the one obtained by differentiating in $r,$ which reads
\begin{align}
\frac{\partial S}{\partial r}  &  =\operatorname{Re}\left[  \frac{\tau}%
{4}\left(  1-\left(  1-\varepsilon\frac{\partial S}{\partial\varepsilon
}-2\lambda\frac{\partial S}{\partial\lambda}\right)  ^{2}\right)  \right]
\nonumber\\
&  \qquad+s\operatorname{Re}\left[  \lambda^{2}\left(  \frac{\partial
S}{\partial\lambda}\right)  ^{2}-\lambda\frac{\partial S}{\partial\lambda
}+\frac{\left\vert \lambda\right\vert ^{2}}{4}\left(  \frac{\partial
S}{\partial\varepsilon}\right)  ^{2}\right]  . \label{rPDE}%
\end{align}
(Apply Theorem \ref{thm:rPDE} below with $(s,\tau)$ replaced by $(0,0)$ and
$(s^{\prime},\tau^{\prime})$ replaced by $(s,\tau).$) Although the first term
on the right-hand side of (\ref{rPDE}) is closely related to (\ref{thePDE}),
there is a complicated second term. Second, in the PDE (\ref{thePDE}), it is
possible (at least formally) to set $\varepsilon=0$ and obtain a
self-contained PDE\ for
\[
S_{0}(s,\tau,\lambda):=S(s,\tau,\lambda,\varepsilon).
\]
That is to say, if we set $\varepsilon=0$ in (\ref{thePDE}), all derivatives
with respect to $\varepsilon$ also disappear, which is not the case for the
second term in the PDE (\ref{rPDE}). Since we are ultimately interested in
setting $\varepsilon=0,$ this property of the PDE\ will prove useful. Third,
related to the second, differentiating with respect to $\tau$ with $s$ fixed
will help us understand a crucial property of the Brown measures $\mu_{s,\tau
},$ namely that there is close relationship between two Brown measures with
the same value of $s$ but different values of $\tau.$ To understand this
property from the PDE\ perspective, we first set $\varepsilon=0,$ after
justifying that this is allowed. As we have noted, putting $\varepsilon=0$
gives a self-contained PDE\ for the \textit{un}-regularized log potential
$S_{0}$ of $ub_{s,\tau}.$ We will then apply the Laplacian to both sides of
this equation, obtaining a formula for how the Brown measure $\mu_{s,\tau}$
changes as $\tau$ changes with $s$ fixed. See Section \ref{push.sec}.

One disadvantage of differentiating in $\tau$ with $s$ fixed is that the
initial condition---computed from $ub_{s,\tau}$ at $\tau=0$---is not
\textquotedblleft trivial.\textquotedblright\ Specifically, when $\tau=0,$ we
have%
\[
ub_{s,0}=uu_{s},
\]
where $u_{s}$ is Biane's free unitary Brownian motion and $u$ is our unitary
element, assumed to be freely independent of $u_{s}.$ By contrast, if we
differentiated $b_{s,\tau}(r)$ with respect to $r,$ then the initial condition
would be computed from the value of $ub_{s,\tau}(r)$ at $r=0,$ which is simply
$u.$ Fortunately, Zhong \cite{Zhong} has computed the law of $uu_{s}$ in terms
of the law of $u$ using Biane's computation of the law of $u_{s}$ and the
technique of free multiplicative convolution. We will use Zhong's results in
Section \ref{e0ode.sec} (and then indirectly in Sections \ref{outside.sec} and
\ref{inside.sec}).

Although our PDE is in many ways simpler than the one in \cite{DHKBrown} and
\cite{HZ}, the analysis of it requires new techniques. Thus, the technical
details involved in computing the Brown measure are completely different here
as compared to those earlier papers. See the beginning of Section
\ref{inside.sec} for more information.

\subsection{It\^{o} rules\label{itoRules.sec}}

The theory of stochastic integration with respect to a semicircular Brownian
motion was developed by Biane and Speicher \cite{BS1}. This paper is
foundational to the whole theory of free stochastic calculus and gave the
first version of a free It\^{o} formula \cite[Proposition 4.3.4]{BS1}. We will
need an extension of that result to a pair of freely independent semicircular
Brownian motions. If $x_{r}$ and $\tilde{x}_{r}$ are two freely independent
semicircular Brownian motions, they satisfy the following free It\^{o} rules
for a continuous adapted process $A_{r}$:
\begin{align}
\operatorname{tr}[A_{r}\,dx_{r}]  &  =\operatorname{tr}[A_{r}\,d\tilde{x}%
_{r}]=0\nonumber\\
dx_{r}\,A_{r}\,dx_{r}  &  =d\tilde{x}_{r}\,A_{r}\,d\tilde{x}_{r}%
=\operatorname{tr}[A_{r}]\,dr\nonumber\\
dx_{r}\,A_{r}\,d\tilde{x}_{r}  &  =d\tilde{x}_{r}\,A_{r}\,dx_{r}=0\nonumber\\
dx_{r}~A_{r}~dr  &  =dr~A_{r}~dx_{r}=0, \label{eq:scIto}%
\end{align}
together with a stochastic product rule for processes $m_{r}^{1}$ and
$m_{r}^{2}$ satisfying SDEs involving $x_{r}$ and $\tilde{x}_{r}$
\begin{equation}
d(m^{1}m^{2})_{r}=dm_{r}^{1}\,m_{r}^{2}+m_{1}^{1}\,dm_{r}^{2}+(dm_{r}%
^{1})(dm_{r}^{2}). \label{eq:Itoprod}%
\end{equation}
These rules are widely used in the literature and were established by
K\"{u}mmerer and Speicher in the setting of the Cuntz algebra \cite{KumSpeich}%
. We will use a general form of these rules developed by Nikitopoulos
\cite[Theorem 3.2.5]{Nik}, stated in an equivalent form for a circular
Brownian motion and its adjoint.

The following It\^{o} rules for the Brownian motions $w_{s,\tau}$ introduced
in Section \ref{freeBM.sec} can be obtained easily from (\ref{eq:scIto}):
\begin{equation}
\operatorname{tr}[A_{r}\,dw_{s,\tau}(r)]=\operatorname{tr}[A_{r}\,dw_{s,\tau
}^{\ast}(r)]=0, \label{eq:ItoTr}%
\end{equation}
and
\begin{align}
dw_{s,\tau}(r)\,A_{r}\,dw_{s,\tau}^{\ast}(r)  &  =dw_{s,\tau}^{\ast}%
(r)\,A_{r}\,dw_{s,\tau}(r)=s\operatorname{tr}[A_{r}]\,dr\nonumber\\
dw_{s,\tau}(r)\,A_{r}\,dw_{s,\tau}(r)  &  =(s-\tau)\operatorname{tr}%
[A_{r}]\,dr\nonumber\\
dw_{s,\tau}^{\ast}(r)\,A_{r}\,dw_{s,\tau}^{\ast}(r)  &  =(s-\bar{\tau
})\operatorname{tr}[A_{r}]\,dr \label{eq:3Ito}%
\end{align}
and%
\begin{align}
dw_{s,\tau}(r)~A_{r}~\,dr  &  =dr~A_{r}\,~dw_{s,\tau}(r)=0\nonumber\\
dw_{s,\tau}^{\ast}(r)~A_{r}~\,dr  &  =dr~A_{r}\,~dw_{s,\tau}^{\ast}(r)=0.
\label{dwdrZero}%
\end{align}

\subsection{The PDE with respect to $r$}

We will make use of the following result.

\begin{theorem}
[Factorization Theorem]\label{thm:factorization} Choose $s,s^{\prime}\geq0$
and $\tau,\tau^{\prime}\in\mathbb{C}$ so that $\left\vert \tau-s\right\vert
\leq s$ and $\left\vert \tau^{\prime}-s^{\prime}\right\vert \leq s^{\prime}$.
Let $b_{s,\tau}$ and $b_{s^{\prime},\tau^{\prime}}^{\prime}$ be two elements
constructed as in Notation \ref{bstau.notation} but chosen to be freely
independent. Then $b_{s+s^{\prime},\tau+\tau^{\prime}}$ and $b_{s,\tau
}b_{s^{\prime},\tau^{\prime}}^{\prime}$ have the same $\ast$-distribution.
\end{theorem}

The proof of Theorem \ref{thm:factorization} is subtle and is deferred to
Appendix \ref{factorization.app}.

Define%
\[
g=ub_{s,\tau}b_{s^{\prime},\tau^{\prime}}^{\prime}(r)
\]
and then set%
\begin{equation}
P(r,\lambda,\varepsilon)=\operatorname{tr}[\log((g-\lambda)^{\ast}%
(g-\lambda)+\varepsilon^{2})], \label{eq:Tdef}%
\end{equation}
where $b_{s^{\prime},\tau^{\prime}}^{\prime}(r)$ is freely independent of
$b_{s,\tau}$ and $u$ and where $s,$ $s^{\prime},$ $\tau,$ and $\tau^{\prime}$
are fixed. By Theorem \ref{thm:factorization} and the fact that $b_{s^{\prime
},\tau^{\prime}}(r)$ has the same $\ast$-distribution as $b_{rs^{\prime}%
,r\tau^{\prime}}(1),$ we see that
\begin{equation}
P(r,\lambda,\varepsilon)=S(s+rs^{\prime},\tau+r\tau^{\prime},\lambda
,\varepsilon). \label{TandS}%
\end{equation}
We note that the process $ub_{s,\tau}b_{s^{\prime},\tau^{\prime}}^{\prime}(r)$
satisfies the same SDE as $b_{s^{\prime},\tau^{\prime}}^{\prime}(r),$ but with
a different initial condition.

\begin{theorem}
\label{thm:rPDE} The function $P$ satisfies the PDE
\begin{align}
\frac{\partial P}{\partial r}  &  =\operatorname{Re}\left[  \frac{\tau
^{\prime}}{4}\left(  1-\left(  1-\varepsilon\frac{\partial P}{\partial
\varepsilon}-2\lambda\frac{\partial P}{\partial\lambda}\right)  ^{2}\right)
\right] \nonumber\\
&  +s^{\prime}\operatorname{Re}\left[  \lambda^{2}\left(  \frac{\partial
P}{\partial\lambda}\right)  ^{2}-\lambda\frac{\partial P}{\partial\lambda
}+\frac{\left\vert \lambda\right\vert ^{2}}{4}\left(  \frac{\partial
P}{\partial\varepsilon}\right)  ^{2}\right]  . \label{eq:rPDE}%
\end{align}

\end{theorem}

\begin{remark}
If we take $s=\tau=0$ and $s^{\prime}=\tau^{\prime}=1,$ then the function $P$
in (\ref{eq:Tdef}) becomes $S(r,r,\lambda,\varepsilon),$ which is the function
considered in \cite{DHKBrown} and \cite{HZ}---except that those papers
regularize with $\varepsilon$ instead of $\varepsilon^{2}.$ In this case, the
PDE (\ref{eq:rPDE}) simplifies to
\[
\frac{\partial P}{\partial r}=\frac{\partial P}{\partial\varepsilon}\left(
\frac{\varepsilon}{2}+\frac{1}{4}\left(  \left\vert \lambda\right\vert
^{2}-\varepsilon^{2}\right)  \frac{\partial P}{\partial\varepsilon}%
-\frac{\varepsilon}{2}\left(  x\frac{\partial P}{\partial x}+y\frac{\partial
P}{\partial y}\right)  \right)
\]
where $\lambda=x+iy$. This is just the PDE in \cite{DHKBrown} and \cite{HZ},
after making the change of variable $\varepsilon\mapsto\varepsilon^{2}$ to
account for the difference in regularization ($\varepsilon$ in \cite{DHKBrown}
and \cite{HZ} and $\varepsilon^{2}$ here).
\end{remark}

We will prove Theorem \ref{thm:rPDE} after establishing some preliminary
results. We use the notation%
\begin{align*}
g  &  =ub_{s,\tau}b_{s^{\prime},\tau^{\prime}}^{\prime}(r)\\
g_{\lambda}  &  =g-\lambda\\
m_{r}  &  =g_{\lambda}^{\ast}g_{\lambda}\\
R  &  =(m_{r}+\varepsilon^{2})^{-1}.
\end{align*}
Note that
\begin{equation}
P(r,\lambda,\varepsilon)=\operatorname{tr}[\log(m_{r}+\varepsilon^{2})].
\label{PfromMr}%
\end{equation}

\begin{proposition}
\label{rDeriv.prop}We have
\begin{equation}
\frac{\partial}{\partial r}P(r,\lambda,\varepsilon)=\frac{\operatorname{tr}%
\left[  R~dm_{r}\right]  }{dr}-\frac{1}{2}\frac{\operatorname{tr}\left[
R~dm_{r}~R~dm_{r}\right]  }{dr} \label{functionalIto}%
\end{equation}
where
\begin{equation}
\frac{\operatorname{tr}\left[  R~dm_{r}\right]  }{dr}=-\frac{s^{\prime}%
-\bar{\tau}^{\prime}}{2}\operatorname{tr}\left[  R\,g^{\ast}g_{\lambda
}\right]  -\frac{s^{\prime}-\tau^{\prime}}{2}\operatorname{tr}\left[
R\,g_{\lambda}^{\ast}g\right]  +s^{\prime}\operatorname{tr}\left[  g^{\ast
}g\right]  \operatorname{tr}\left[  R\right]  \label{dSfirst}%
\end{equation}
and
\begin{align}
\frac{\operatorname{tr}\left[  R~dm_{r}~R~dm_{r}\right]  }{dr}  &
=-2\operatorname{Re}\left\{  (s^{\prime}-\tau^{\prime})\operatorname{tr}%
\left[  Rg_{\lambda}^{\ast}g\right]  \operatorname{tr}\left[  Rg_{\lambda
}^{\ast}g\right]  \right\} \nonumber\\
&  +2s^{\prime}\operatorname{tr}\left[  R\right]  \operatorname{tr}\left[
Rg_{\lambda}^{\ast}gg^{\ast}g_{\lambda}\right]  . \label{dSsecond}%
\end{align}

\end{proposition}

Observe that the right-hand sides of (\ref{dSfirst}) and (\ref{dSsecond})
involve $g$ and $g^{\ast}$ in addition to $g_{\lambda}=g-\lambda$ and
$g_{\lambda}^{\ast}=(g-\lambda)^{\ast}.$ Since the derivatives of $P$ with
respect to $\lambda$ and $\varepsilon$ (Lemma \ref{lem:CollDerivatives})
involve only $g_{\lambda}$ and $g_{\lambda}^{\ast},$ we will eventually need
to rewrite $g$ as $g_{\lambda}+\lambda$ and simplify in (\ref{dSfirst}) and
(\ref{dSsecond}).

\begin{proof}
Using (\ref{PfromMr}), we can see that the expression (\ref{functionalIto}) is
a consequence of Eq. (33) in Example 3.5.5 of \cite{Nik}. To compute more
explicitly, we use the stochastic product rule (\ref{eq:Itoprod}) to obtain%
\begin{equation}
dm_{r}=(dg^{\ast})g_{\lambda}+g_{\lambda}^{\ast}(dg)+(dg^{\ast})(dg),
\label{dmr}%
\end{equation}
and then compute $dg$ and $dg^{\ast}$ (thinking of $r$ as the time variable)
using the SDE (\ref{bstauSDE}) for $b_{s^{\prime},\tau^{\prime}}(r).$ When
computing $\operatorname{tr}\left[  R~dm_{r}\right]  /dr,$ we use
(\ref{eq:ItoTr}) and (\ref{dwdrZero}), which tell us that we keep only the
$dr$ terms in the formulas for $dg$ and $dg^{\ast}$ in the first two terms on
the right-hand side of (\ref{dmr}), but that we can drop the $dr$ terms in the
last term on the right-hand side of (\ref{dmr}). Three terms will remain,
which correspond to the three terms in (\ref{dSfirst}).

We then compute that
\begin{equation}
\operatorname{tr}\left[  R\,dm_{r}\,R\,dm_{r}\right]  =\operatorname{tr}%
\left\{  R\left[  (dg^{\ast})g_{\lambda}+g_{\lambda}^{\ast}(dg)\right]
R\left[  (dg^{\ast})g_{\lambda}+g_{\lambda}^{\ast}(dg)\right]  \right\}  ,
\label{RdmRdm}%
\end{equation}
where we note that the $(dg^{\ast})(dg)$ term in (\ref{dmr}) does not
contribute here. Furthermore, the $dr$ terms in the formula for $dg$ and
$dg^{\ast}$ do not contribute. We are then left with four terms, which we
compute using (\ref{eq:3Ito}) to obtain (\ref{dSsecond}).
\end{proof}

\begin{lemma}
\label{lem:CollDerivatives}We have
\[
\frac{\partial P}{\partial\varepsilon}=2\varepsilon\operatorname{tr}\left[
R\right]  ;\quad\frac{\partial P}{\partial\lambda}=-\operatorname{tr}\left[
Rg_{\lambda}^{\ast}\right]  ;\quad\frac{\partial P}{\partial\bar{\lambda}%
}=-\operatorname{tr}\left[  Rg_{\lambda}\right]  .
\]

\end{lemma}

\begin{proof}
Direct calculation using Lemma 1.1 in Brown's paper \cite{Br}, which states
that
\[
\frac{d}{du}\operatorname{tr}\left[  \log(f(u))\right]  =\operatorname{tr}%
\left[  f(u)^{-1}\frac{df}{du}\right]
\]
for any smooth function $f(u)$ taking values in the set of strictly positive
elements of $\mathcal{A}.$
\end{proof}

We now provide the proof of Theorem \ref{thm:rPDE}.

\begin{proof}
[Proof of Theorem \ref{thm:rPDE}]Our strategy is (1) to simplify the results
in (\ref{dSfirst}) and (\ref{dSsecond}) using the relations $g=g_{\lambda
}+\lambda,$ $g^{\ast}=g_{\lambda}^{\ast}+\bar{\lambda},$ and
\[
Rg_{\lambda}^{\ast}g_{\lambda}=R(g_{\lambda}^{\ast}g_{\lambda}+\varepsilon
-\varepsilon)=1-\varepsilon R,
\]
and (2) to express the resulting expressions in terms of derivatives of $P$
with respect to $\lambda$ and $\varepsilon$ using Lemma
\ref{lem:CollDerivatives}. There is, however, a difficulty with this approach,
namely that the last term in (\ref{dSfirst}) and the last term in
(\ref{dSsecond}) will give rise to \textquotedblleft bad\textquotedblright%
\ terms (such as $\operatorname{tr}\left[  g_{\lambda}^{\ast}g_{\lambda
}\right]  $) that do not show up in Lemma \ref{lem:CollDerivatives}.
Fortunately, the bad terms arising from (\ref{dSfirst}) cancel with the bad
terms arising from (\ref{dSsecond}). We omit the details of the argument,
which is similar to the proof of Theorem 2.8 in Section 5 of \cite{DHKBrown}.
\end{proof}

\subsection{The PDE with respect to $\tau$\label{PDEpf}}

We now prove Theorem~\ref{thm:dSdtau}, establishing the PDE (\ref{thePDE}).

\begin{proof}
[Proof of Theorem~\ref{thm:dSdtau}]Recall the definition of $P$ in
(\ref{eq:Tdef}) that
\[
P(r,\lambda,\varepsilon)=S(s+rs^{\prime},\tau+r\tau^{\prime},\lambda
,\varepsilon).
\]
If we differentiate both sides of this relation with respect to $r$ at $r=0,$
using the chain rule, we obtain
\[
\left.  \frac{\partial P}{\partial r}\right\vert _{r=0}=s^{\prime}%
\frac{\partial S}{\partial s}(s,\tau,\lambda,\varepsilon)+\tau^{\prime}%
\frac{\partial S}{\partial\tau}(s,\tau,\lambda,\varepsilon)+\bar{\tau}%
^{\prime}\frac{\partial S}{\partial\bar{\tau}^{\prime}}(s,\tau,\lambda
,\varepsilon).
\]
We can then compute $\partial S/\partial\tau$ as%
\[
\frac{\partial S}{\partial\tau}=\frac{\partial}{\partial\tau^{\prime}}\left(
\left.  \frac{\partial P}{\partial r}\right\vert _{r=0}\right)  .
\]
If we evaluate the PDE (\ref{eq:rPDE}) at $r=0$ and then differentiate the
right-hand side with respect to $\tau^{\prime},$ we obtain (\ref{thePDE}).
\end{proof}

\section{The Hamilton--Jacobi analysis\label{hj.sec}}

\subsection{ The complex-time Hamilton--Jacobi formulas}

We now analyze solutions to the PDE (\ref{thePDE}). We introduce the
complex-valued Hamiltonian function by replacing $\partial S/\partial\lambda$
by $p_{\lambda}$ and $\partial S/\partial\varepsilon$ by $p_{\varepsilon}$ on
the right-hand side of (\ref{thePDE}), with an overall minus sign:%
\begin{equation}
H(\lambda,\varepsilon,p_{\lambda},p_{\varepsilon})=-\frac{1}{8}\left[
1-\left(  1-\varepsilon p_{\varepsilon}-2\lambda p_{\lambda}\right)
^{2}\right]  . \label{compexHamiltonian}%
\end{equation}
Here the variables $\lambda$ and $p_{\lambda}$ are complex valued and the
variables $\varepsilon$ and $p_{\varepsilon}$ are real valued. Define
\begin{align}
p_{\lambda,0}  &  =\frac{\partial}{\partial\lambda_{0}}S(s,0,\lambda
_{0},\varepsilon_{0})\label{InitialMomentum1}\\
p_{\varepsilon,0}  &  =\frac{\partial}{\partial\varepsilon_{0}}S(s,0,\lambda
_{0},\varepsilon_{0}). \label{InitialMomentum2}%
\end{align}
Writing $S(s,0,\lambda_{0},\varepsilon_{0})$ as in
(\ref{theInitialCondIntegral}), we obtain
\begin{equation}
p_{\lambda,0}=-\int_{S^{1}}\frac{\bar{\xi}-\bar{\lambda}}{\left\vert
\xi-\lambda_{0}\right\vert ^{2}+\varepsilon_{0}^{2}}~d\mu_{s}(\xi)
\label{pLambda0}%
\end{equation}
and
\begin{equation}
p_{\varepsilon,0}=\int_{S^{1}}\frac{2\varepsilon_{0}}{\left\vert \xi
-\lambda_{0}\right\vert ^{2}+\varepsilon_{0}^{2}}~d\mu_{s}(\xi). \label{p0}%
\end{equation}

Then define curves $\lambda(\tau),$ $\varepsilon(\tau),$ $p_{\lambda}(\tau)$,
and $p_{\varepsilon}(\tau)$ by%
\begin{align}
\lambda(\tau)  &  =\lambda_{0}\exp\left\{  \frac{\tau}{2}(\varepsilon
_{0}p_{\varepsilon,0}+2\lambda_{0}p_{\lambda,0}-1)\right\}
\label{lambdaOfTau}\\
\varepsilon(\tau)  &  =\varepsilon_{0}\exp\left\{  \operatorname{Re}\left[
\frac{\tau}{2}(\varepsilon_{0}p_{\varepsilon,0}+2\lambda_{0}p_{\lambda
,0}-1)\right]  \right\} \label{epsilonOfTau}\\
p_{\lambda}(\tau)  &  =p_{\lambda,0}\exp\left\{  -\frac{\tau}{2}%
(\varepsilon_{0}p_{\varepsilon,0}+2\lambda_{0}p_{\lambda,0}-1)\right\}
\label{pLambdaOfTau}\\
p_{\varepsilon}(\tau)  &  =p_{\varepsilon,0}\exp\left\{  -\operatorname{Re}%
\left[  \frac{\tau}{2}(\varepsilon_{0}p_{\varepsilon,0}+2\lambda_{0}%
p_{\lambda,0}-1)\right]  \right\}  \label{pEpsilonOfTau}%
\end{align}

We now derive Hamilton--Jacobi formulas for solutions to the PDE (\ref{thePDE}).

\begin{theorem}
\label{HJ3param.thm} For all $\tau$ with $\left\vert \tau-s\right\vert \leq
s,$ we have the first Hamilton--Jacobi formula
\begin{equation}
S(s,\tau,\lambda(\tau),\varepsilon(\tau))=S(s,0,\lambda_{0},\varepsilon
_{0})+2\operatorname{Re}[\tau H_{0}]+\frac{1}{2}\operatorname{Re}\left[
\tau(\varepsilon_{0}p_{\varepsilon,0}+2\lambda_{0}p_{\lambda,0})\right]  ,
\label{HJfirst}%
\end{equation}
where $H_{0}=H(\lambda_{0},\varepsilon_{0},p_{\lambda,0},p_{\varepsilon,0})$,
and the second Hamilton--Jacobi formulas
\begin{align}
\frac{\partial S}{\partial\lambda}(s,\tau,\lambda(\tau),\varepsilon(\tau))  &
=p_{\lambda}(\tau)\label{HJsecond1}\\
\frac{\partial S}{\partial\varepsilon}(s,\tau,\lambda(\tau),\varepsilon
(\tau))  &  =p_{\varepsilon}(\tau). \label{HJsecond2}%
\end{align}

\end{theorem}

Note that the ordinary Hamilton--Jacobi method is not directly applicable to
the PDE (\ref{thePDE}), because the \textquotedblleft time\textquotedblright%
\ variable $\tau$ in that equation is complex. Our method for proving Theorem
\ref{HJ3param.thm} will then be to reduce the problem to an ordinary
Hamilton--Jacobi PDE\ with real time variable.

\subsection{Ordinary Hamilton--Jacobi method}

Our strategy for proving Theorem \ref{HJ3param.thm} is to consider the
function
\begin{equation}
S^{\tau}(s,t,\lambda,\varepsilon):=S(s,t\tau,\lambda,\varepsilon),
\label{SsTauDef}%
\end{equation}
where $t$ is a positive real number. Note that on the right-hand side of
(\ref{SsTauDef}), we are scaling the complex number $\tau$ by the positive
real number $t.$ With $\tau$ fixed, $S^{\tau}(s,t,\lambda,\varepsilon)$ is
defined for all $t$ such that $\left\vert t\tau-s\right\vert \leq s.$

\begin{proposition}
\label{dStauDt.thm}The function $S^{\tau}$ satisfies the PDE
\begin{equation}
\frac{\partial S^{\tau}}{\partial t}=\operatorname{Re}\left[  \frac{\tau}%
{4}\left(  1-\left(  1-\varepsilon\frac{\partial S^{\tau}}{\partial
\varepsilon}-2\lambda\frac{\partial S^{\tau}}{\partial\lambda}\right)
^{2}\right)  \right]  \label{SPDEt}%
\end{equation}
for all $t$ such that $\left\vert t\tau-s\right\vert <s,$ with the initial
condition at $t=0$ given by
\begin{equation}
S^{\tau}(s,0,\lambda,\varepsilon)=\operatorname{tr}[\log((uu_{s}%
-\lambda)^{\ast}(uu_{s}-\lambda)+\varepsilon^{2})] \label{sPDEinitial}%
\end{equation}

\end{proposition}

The PDE\ for $S^{\tau}$ is of ordinary Hamilton--Jacobi type with respect to
the real time variable $t.$ Note that the coefficients of the PDE in
(\ref{SPDEt}) do not depend on $s$; the $s$-dependence in the problem is only
in the initial condition (\ref{sPDEinitial}).

\begin{proof}
Using the chain rule in the definition (\ref{SsTauDef}) of $S^{\tau},$ we find
that
\[
\frac{\partial S^{\tau}}{\partial t}=\tau\frac{\partial S}{\partial\tau}%
+\bar{\tau}\frac{\partial S}{\partial\bar{\tau}}.
\]
The result now follows from the PDE (\ref{thePDE}) for $S.$
\end{proof}

We will analyze the function $S^{\tau}$ using the Hamilton--Jacobi method. The
following proposition recalls the general Hamilton--Jacobi method.

\begin{proposition}
\label{HJgeneral.prop}Fix an open set $U\subset\mathbb{R}^{n},$ a
time-interval $[0,T],$ and a function $H(\mathbf{x},\mathbf{p}):U\times
\mathbb{R}^{n}\rightarrow\mathbb{R}.$ Consider a function $S(t,\mathbf{x})$
satisfying
\[
\frac{\partial S}{\partial t}=-H(\mathbf{x},\nabla_{\mathbf{x}}S),\quad
\mathbf{x}\in U,~t\in\lbrack0,T].
\]
Consider a curve $(\mathbf{x}(t),\mathbf{p}(t))$ with $\mathbf{x}(t)\in U,$
$\mathbf{p}(t)\in\mathbb{R}^{n},$ and $t$ ranging over an interval $[0,T_{1}]$
with $T_{1}\leq T.$ Assume this curve satisfies Hamilton's equations:
\begin{equation}
\frac{dx_{j}}{dt}=\frac{\partial H}{\partial p_{j}}(\mathbf{x}(t),\mathbf{p}%
(t));\quad\frac{dp_{j}}{dt}=-\frac{\partial H}{\partial x_{j}}(\mathbf{x}%
(t),\mathbf{p}(t)) \label{HamEq}%
\end{equation}
with initial conditions of the special form
\begin{equation}
\mathbf{x}(0)=\mathbf{x}_{0};\quad\mathbf{p}(0)=(\nabla_{\mathbf{x}%
}S)(0,\mathbf{x}_{0}). \label{HJinit}%
\end{equation}
Then we have the first Hamilton--Jacobi formula
\begin{equation}
S(t,\mathbf{x}(t))=S(0,\mathbf{x}_{0})-t~H(\mathbf{x}_{0},\mathbf{p}_{0}%
)+\int_{0}^{t}\mathbf{p}(s)\cdot\frac{d\mathbf{x}}{ds}~ds \label{HJgen1}%
\end{equation}
and the second Hamilton--Jacobi formula
\begin{equation}
(\nabla_{\mathbf{x}}S)(t,\mathbf{x}(t))=\mathbf{p}(t). \label{HJgen2}%
\end{equation}

\end{proposition}

See the proof of Proposition~\ref{HJgeneral.prop} in \cite{DHKBrown}\ for a
self-contained proof and see also Section 3.3 in the book \cite{Evans} of Evans.

We now record how this general result looks in the case of the function
$S^{\tau}(s,t,\lambda,\varepsilon)$, taking the open set $U$ to be the set of
$(s,\lambda,\varepsilon)$ with $\varepsilon>0.$ Recall that we write
$\lambda=x+iy$. We use the notation%
\begin{align}
p_{\lambda}  &  =\frac{1}{2}(p_{x}-ip_{y})\nonumber\\
p_{\bar{\lambda}}  &  =\frac{1}{2}(p_{x}+ip_{y}). \label{pLambdaDef}%
\end{align}
This is a convenient notation, so that $p_{\lambda}$ and $p_{\bar{\lambda}}$
correspond to $\partial S^{\tau}/\partial\lambda$ and $\partial S^{\tau
}/\partial\bar{\lambda},$ respectively, in the PDE and in the second
Hamilton--Jacobi formula. We then compute that
\begin{align}
\left(  \frac{\partial}{\partial p_{x}}+i\frac{\partial}{\partial p_{y}%
}\right)  p_{\lambda}  &  =1;\quad\left(  \frac{\partial}{\partial p_{x}%
}-i\frac{\partial}{\partial p_{y}}\right)  p_{\lambda}=0;\nonumber\\
\left(  \frac{\partial}{\partial p_{x}}+i\frac{\partial}{\partial p_{y}%
}\right)  p_{\bar{\lambda}}  &  =0;\quad\left(  \frac{\partial}{\partial
p_{x}}-i\frac{\partial}{\partial p_{y}}\right)  p_{\bar{\lambda}}=1.
\label{derivPlambda}%
\end{align}
The Hamiltonian, read off (with a minus sign) from the PDE (\ref{SPDEt}), is
\begin{equation}
H^{\tau}=-\operatorname{Re}\left[  \frac{\tau}{4}\left(  1-\left(
1-\varepsilon p_{\varepsilon}-2\lambda p_{\lambda}\right)  ^{2}\right)
\right]  . \label{Ham}%
\end{equation}

We now use the notation $\lambda^{\tau},$ $\varepsilon^{\tau},$ $p_{\lambda
}^{\tau},$ and $p_{\varepsilon}^{\tau}$ for solutions to Hamilton's equations
with Hamiltonian $H^{\tau}.$

\begin{proposition}
\label{HJ3param.prop}Suppose $(\lambda^{\tau}(t),\varepsilon^{\tau
}(t),p_{\lambda}^{\tau}(t),p_{\varepsilon}^{\tau}(t))$ satisfies the
Hamilton's equations (\ref{HamEq}) with Hamiltonian $H^{\tau},$ with
$\varepsilon^{\tau}(t)>0$ and initial conditions computed as in (\ref{HJinit}%
). Then for all $t>0$ such that $\left\vert t\tau-s\right\vert <s,$ we have
\begin{equation}
S^{\tau}(s,t,\lambda^{\tau}(t),\varepsilon^{\tau}(t))=S^{\tau}(s,0,\lambda
_{0},\varepsilon_{0})+tH_{0}^{\tau}+t\operatorname{Re}\left[  \frac{\tau}%
{2}(\varepsilon_{0}p_{\varepsilon,0}+2\lambda_{0}p_{\lambda,0})\right]  ,
\label{HJ1st}%
\end{equation}
where $H_{0}^{\tau}$ is the value of the Hamiltonian at $t=0.$ Furthermore,%
\begin{align}
\frac{\partial S^{\tau}}{\partial\lambda}(s,t,\lambda^{\tau}(t),\varepsilon
^{\tau}(t))  &  =p_{\lambda}^{\tau}(t)\nonumber\\
\frac{\partial S^{\tau}}{\partial\bar{\lambda}}(s,t,\lambda^{\tau
}(t),\varepsilon^{\tau}(t))  &  =p_{\bar{\lambda}}^{\tau}(t)\nonumber\\
\frac{\partial S^{\tau}}{\partial\varepsilon}(s,t,\lambda^{\tau}%
(t),\varepsilon^{\tau}(t))  &  =p_{\varepsilon}^{\tau}(t). \label{HJ2nd}%
\end{align}

\end{proposition}

\begin{proof}
In the general Hamilton--Jacobi method, we note that $dx_{j}/ds=\partial
H/\partial p_{j},$ so that $\mathbf{p}(s)\cdot d\mathbf{x}/ds=\mathbf{p}%
(s)\cdot\nabla_{\mathbf{p}}H(\mathbf{x}(s),\mathbf{p}(s)).$ We then note that
the operator $\mathbf{p}\cdot\nabla_{\mathbf{p}}$ is the homogeneity operator
in the $\mathbf{p}$ variables, acting as $k$ times the identity on a function
that is homogeneous of degree $k$ in $\mathbf{p}.$ We then write the
Hamiltonian (\ref{Ham}) as $H^{\tau}=H_{2}+H_{1},$ with%
\[
H_{2}=\operatorname{Re}\left[  \frac{\tau}{4}(\varepsilon p_{\varepsilon
}+2\lambda p_{\lambda})^{2}\right]  ;\quad H_{1}=-\operatorname{Re}\left[
\frac{\tau}{2}\left(  \varepsilon p_{\varepsilon}+2\lambda p_{\lambda}\right)
\right]  ,
\]
where $H_{2}$ and $H_{1}$ are, respectively, homogeneous of degrees 2 and 1 in
$\mathbf{p}.$ (The terms of degree zero in $\mathbf{p}$ cancel.) Thus,
\[
\mathbf{p}\cdot\nabla_{\mathbf{p}}H^{\tau}=2H_{2}+H_{1}=2H^{\tau}-H_{1}.
\]

It is then easy to check that $H^{\tau}$ and $H_{1}$ are both constants of
motion, so that $\int_{0}^{t}\mathbf{p}(s)\cdot d\mathbf{x}/ds~ds$ equals $t$
times the value of $2H^{\tau}-H_{1}$ at $t=0.$ The claimed formula
(\ref{HJ1st}) then follows easily. Meanwhile, the formulas in (\ref{HJ2nd})
follow immediately from the general formula (\ref{HJgen2}), together with the
definitions (\ref{pLambdaDef}) of $p_{\lambda}$ and $p_{\bar{\lambda}}.$
\end{proof}

\begin{proposition}
\label{limitHJ.prop}Suppose $\tau$ is a nonzero complex number such that
$\left\vert \tau-s\right\vert =s.$ Then the Hamilton--Jacobi formulas
(\ref{HJ1st}) and (\ref{HJ2nd}) continue to hold at $t=1.$
\end{proposition}

Note that in the borderline case $\left\vert \tau-s\right\vert =s,$ we have
only established that the PDE (\ref{SPDEt}) holds for $t<1.$ (See Theorem
\ref{dStauDt.thm}.) The Proposition says that the Hamilton--Jacobi formulas
nevertheless remain valid at $t=1.$

\begin{proof}
We will show in Appendix \ref{continuousDependence.sec} that the element
$b_{s,\tau}$ in Notation \ref{bstau.notation} depends continuously (in the
operator norm topology) on the parameter $\tau,$ for all $\tau$ satisfying
$\left\vert \tau-s\right\vert \leq s.$ (See Theorem \ref{continuousDep.thm}.)
Let us then fix a nonzero $\tau$ with $\left\vert \tau-s\right\vert =s.$ We
will use the continuity of the operator logarithm on strictly positive
self-adjoint operators, which may be established by expanding the function
$\log x$ in a Taylor series based at $a,$ for some large positive $a.$ It is
then easily verified that the function $S^{\tau}(s,t,\lambda,\varepsilon)$ in
and its derivatives with respect to $\lambda$ and $\varepsilon$ (computed in
Lemma \ref{lem:CollDerivatives}) depend continuously on the operator
$b_{s,\tau}$ and therefore also on $\tau,$ even up to the boundary. We may
therefore let $t$ tend to 1 on both sides of the Hamilton--Jacobi formulas
(\ref{HJ1st}) and (\ref{HJ2nd}).
\end{proof}

Following the general Hamilton--Jacobi method (as in (\ref{HJinit})), the
initial momentum $p_{\lambda,0}$ is computed as%
\[
p_{\lambda,0}=\frac{\partial}{\partial\lambda_{0}}S^{\tau}(s,0,\lambda
_{0},\varepsilon_{0}).
\]
Since $S^{\tau}(s,0,\lambda_{0},\varepsilon_{0})=S(s,0,\lambda_{0}%
,\varepsilon_{0}),$ we see that the initial momentum $p_{\lambda,0}$
associated to $S^{\tau}$ agrees with the initial momentum $p_{\lambda,0}$
associated to $S,$ as in (\ref{InitialMomentum1}). A similar statement applies
to the initial momentum $p_{\varepsilon,0}.$

\subsection{Proof of the main theorem}

We now prove Theorem \ref{HJ3param.thm}.

\begin{proof}
[Proof of Theorem \ref{HJ3param.thm}]From the definition (\ref{SsTauDef}) of
$S^{\tau},$ we have
\begin{equation}
S(s,\tau,\lambda,\varepsilon)=\left.  S^{\tau}(s,t,\lambda,\varepsilon
)\right\vert _{t=1}. \label{StauAt1}%
\end{equation}
Thus, to compute $S(s,\tau,\lambda,\varepsilon),$ we will apply the
Hamilton--Jacobi analysis in the previous subsection and then evaluate at
$t=1.$ It is therefore natural to define%
\begin{align}
\lambda(\tau)  &  =\left.  \lambda^{\tau}(t)\right\vert _{t=1};\quad
\varepsilon(\tau)=\left.  \varepsilon^{\tau}(t)\right\vert _{t=1};\nonumber\\
p_{\lambda}(\tau)  &  =\left.  p_{\lambda}^{\tau}(t)\right\vert _{t=1}%
;\quad\left.  p_{\varepsilon}(\tau)=p_{\varepsilon}^{\tau}(t)\right\vert
_{t=1}. \label{curvesTau}%
\end{align}
As noted at the end of the previous subsection, the initial momenta
$p_{\lambda,0}$ and $p_{\varepsilon,0}$ associated to $p_{\lambda}^{\tau}$ and
$p_{\varepsilon}^{\tau}$ agree with the initial momenta in the statement of
the complex-time Hamilton--Jacobi formulas, as in (\ref{InitialMomentum1}) and
(\ref{InitialMomentum2}), because $S^{\tau}(s,0,\lambda,\varepsilon)$
coincides with $S(s,0,\lambda,\varepsilon).$

We now show that the curves in (\ref{curvesTau}) are given by the formulas in
(\ref{lambdaOfTau})--(\ref{pEpsilonOfTau}). We have, for example,
\[
\frac{d\lambda^{\tau}}{dt}=\frac{dx^{\tau}}{dt}+i\frac{dy^{\tau}}{dt}=\left(
\frac{\partial}{\partial p_{x}^{\tau}}+i\frac{\partial}{\partial p_{y}^{\tau}%
}\right)  H^{\tau}.
\]
Then in (\ref{Ham}), we write the real part as half the sum of the indicated
expression and its conjugate. Using (\ref{derivPlambda}), this gives
\[
\frac{d\lambda^{\tau}}{dt}=\frac{\tau}{2}(\varepsilon^{\tau}p_{\varepsilon
}^{\tau}+2\lambda^{\tau}p_{\lambda}^{\tau}-1)\lambda^{\tau}.
\]
Since $\lambda^{\tau}p_{\lambda}^{\tau}$ and $\varepsilon^{\tau}%
p_{\varepsilon}^{\tau}$ are constants of motion, as the reader may easily
verify, we find%
\[
\lambda^{\tau}(t)=\lambda_{0}\exp\left\{  \frac{t\tau}{2}(\varepsilon
_{0}p_{\varepsilon,0}+2\lambda_{0}p_{\lambda,0}-1)\right\}  .
\]
Setting $t=1$ gives the curve in (\ref{lambdaOfTau}). The verifications for
the remaining curves are extremely similar.

Using (\ref{StauAt1}), we can easily see that the Hamilton--Jacobi formulas
(\ref{HJfirst}), (\ref{HJsecond1}), and (\ref{HJsecond2}) are just
restatements of (\ref{HJ1st}) and (\ref{HJ2nd}) using the curves defined in
(\ref{curvesTau}). By Proposition \ref{limitHJ.prop}, the formulas continue to
hold in the borderline case $\left\vert \tau-s\right\vert =s$.
\end{proof}

\subsection{The $\varepsilon_{0}=0$ case\label{e0ode.sec}}

We now wish to look at what happens to the curves $\lambda(\tau)$ and
$\varepsilon(\tau)$ in the limit as $\varepsilon_{0}$ approaches zero from
above. (We are not currently claiming that the Hamilton--Jacobi formulas
remain valid up to $\varepsilon_{0}=0$; we are merely looking at the limiting
behavior of the curves.) In the formulas (\ref{lambdaOfTau}) and
(\ref{epsilonOfTau}), the limit is obtained by simply setting $\varepsilon
_{0}=0,$ \textit{provided} that the initial momenta $p_{\varepsilon,0}$ and
$p_{\lambda,0}$ remain defined in this limit. But looking at the formulas in
(\ref{pLambda0}) and (\ref{p0}), we see that if $\varepsilon_{0}=0$ and
$\lambda_{0}$ belongs to the closed support $\operatorname{supp}(\mu_{s})$ of
$\mu_{s},$ then the integrals defining $p_{\varepsilon,0}$ and $p_{\lambda,0}$
may be divergent. We therefore assume $\lambda_{0}\notin\operatorname{supp}%
(\mu_{s})$ in the following result.

\begin{proposition}
\label{solutionse0.prop}Suppose $\lambda_{0}\notin\operatorname{supp}(\mu
_{s}).$ Then with $\varepsilon_{0}=0,$ we have%
\begin{equation}
2\lambda_{0}p_{\lambda,0}-1=-\int_{S^{1}}\frac{\xi+\lambda_{0}}{\xi
-\lambda_{0}}~d\mu_{s}(\xi),\quad(\varepsilon_{0}=0). \label{momentum1}%
\end{equation}
We also have the alternative formula
\begin{equation}
2\lambda_{0}p_{\lambda,0}-1=-\int_{S^{1}}\frac{\xi+\chi_{s}(\lambda_{0})}%
{\xi-\chi_{s}(\lambda_{0})}~d\mu_{0}(\xi),\quad(\varepsilon_{0}=0),
\label{momentum2}%
\end{equation}
where $\chi_{s}$ is the inverse function to $f_{s}.$ Note that the integral in
(\ref{momentum1}) is with respect to $\mu_{s}$ (the law of $uu_{s}$) while the
integral in (\ref{momentum2}) is with respect to $\mu_{0}$ (the law of $u$).

When $\varepsilon_{0}=0,$ the formula for $\lambda(\tau)$ in
(\ref{lambdaOfTau}) becomes
\begin{equation}
\lambda(\tau)=f_{s-\tau}(\chi_{s}(\lambda_{0})), \label{lambdaWithE0}%
\end{equation}
while the formula for $\varepsilon(\tau)$ in (\ref{epsilonOfTau}) becomes%
\[
\varepsilon(\tau)\equiv0.
\]

\end{proposition}

\begin{corollary}
\label{lambda0lambda.cor}For each fixed $\tau,$ consider the map%
\[
\lambda_{0}\mapsto\lambda(\tau)
\]
where $\lambda(\tau)$ is computed with $\lambda(0)=\lambda_{0}$ and
$\varepsilon_{0}=0.$ Then this map is a biholomorphism of $\mathbb{C}%
\setminus\operatorname{supp}(\mu_{s})$ onto $(\bar{\Sigma}_{s,\tau})^{c}.$
\end{corollary}

\begin{proof}
By Definition \ref{sigmas.def}, $\chi_{s}$ maps $\mathbb{C}\setminus
\operatorname{supp}(\mu_{s})$ to the complement of $\tilde{\Sigma}%
_{s}=\overline{\Sigma}_{s}$ and by the last point in Proposition
\ref{moreSigmasChars.prop}, the complement of $\overline{\Sigma}_{s}$ does not
intersect $\operatorname{supp}(\mu_{s}).$ Thus, the Herglotz integral $J$ in
(\ref{herglotzJ}) is defined and holomorphic on the complement of
$\operatorname{supp}(\mu_{s})$ and also on the complement of $\overline
{\Sigma}_{s}.$ Thus, the right-hand side of (\ref{lambdaWithE0}) is a
holomorphic function of $\lambda_{0}.$ The claimed bijectivity of the map then
follows from the definitions of the domains in Definitions \ref{sigmas.def}
and \ref{def:Sigmastaudef}, the injectivity of $f_{s-\tau}$ in Theorem
\ref{fInjective.thm}, and the fact that $\tilde{\Sigma}_{s,\tau}$ is the
closure of its interior. (This last point follows from the proof of Theorem
\ref{thm:Sigmastau} in Section \ref{sect:charSigmastau}.)
\end{proof}

\begin{proof}
[Proof of Proposition \ref{solutionse0.prop}]From the definition
(\ref{pLambdaDef}) of $p_{\lambda},$ we have
\[
p_{\lambda,0}=\frac{\partial}{\partial\lambda_{0}}S(s,0,\lambda_{0}%
,\varepsilon_{0}).
\]
If we then specialize to the case $\varepsilon_{0}=0,$ we have%
\begin{align*}
2\lambda_{0}p_{\lambda,0}-1  &  =2\lambda_{0}\frac{\partial}{\partial
\lambda_{0}}\int_{S^{1}}\log(\left\vert \lambda_{0}-\xi\right\vert ^{2}%
)~d\mu_{s}(\xi)-1\\
&  =2\int_{S^{1}}\frac{\lambda_{0}}{\lambda_{0}-\xi}~d\mu_{s}(\xi)-1\\
&  =-\int_{S^{1}}\frac{\xi+\lambda_{0}}{\xi-\lambda_{0}}~d\mu_{s}(\xi),
\end{align*}
provided that $\lambda_{0}$ is outside the closed support of $\mu_{s},$
establishing (\ref{momentum1}).

We now apply results from Zhong's paper \cite{Zhong}, specifically, the
relation $\eta_{\mu_{t}}(z)=\eta_{\mu}(\eta_{t}(z))$ on p. 1357, after Theorem
2.1, along with Eq. (2.1) and the relation $\Phi_{t,\mu}(\eta_{\rho_{t}%
}(z))=z$ following Lemma 2.2. The relation $\eta_{\mu_{t}}(z)=\eta_{\mu}%
(\eta_{t}(z))$ implies also that $\psi_{\mu_{t}}(z)=\psi_{\mu}(\eta_{t}(z)),$
where $\psi_{\mu}$ is defined on p. 1356 of \cite{Zhong}. If we apply this
last result with Zhong's $\mu$ corresponding to our $\overline{\mu}_{0}$ (the
push-forward of $\mu_{0}$ under the complex-conjugation map) and Zhong's
$\mu_{t}$ corresponding to our $\overline{\mu}_{s},$ we find that Zhong's
function $\eta_{t}$ corresponds to our $\chi_{s}$ and we get the relation
\begin{equation}
\int_{S^{1}}\frac{\xi+\lambda_{0}}{\xi-\lambda_{0}}~d\mu_{s}(\xi)=\int_{S^{1}%
}\frac{\xi+\chi_{s}(\lambda_{0})}{\xi-\chi_{s}(\lambda_{0})}~d\mu_{0}(\xi),
\label{pingsID}%
\end{equation}
establishing (\ref{momentum2}). The identity (\ref{pingsID}) initially holds
for $\left\vert \lambda_{0}\right\vert <1,$ but extends to $\left\vert
\lambda_{0}\right\vert >1$ by the identities (\ref{J1zbar}) and (\ref{fsZbar})
and then to $\lambda_{0}\notin\operatorname{supp}(\mu_{s})$ by continuity.

Meanwhile, when $\varepsilon_{0}=0,$ we have $\varepsilon_{0}p_{\varepsilon
,0}=0$ and (\ref{momentum2})\ applies, so that the formula for $\lambda(\tau)$
in (\ref{lambdaOfTau}) becomes%
\[
\lambda(\tau)=\lambda_{0}\exp\left\{  -\frac{\tau}{2}\int_{S^{1}}\frac
{\xi+\chi_{s}(\lambda_{0})}{\xi-\chi_{s}(\lambda_{0})}~d\mu_{0}(\xi)\right\}
.
\]
Now, we can compute that%
\[
f_{s-\tau}(z)=f_{s}(z)\exp\left\{  -\frac{\tau}{2}\int_{S^{1}}\frac{\xi+z}%
{\xi-z}~d\mu_{0}(\xi)\right\}  .
\]
Applying this result with $z=\chi_{s}(\lambda_{0})$ and recalling that
$\chi_{s}$ is a right inverse of $f_{s}$ gives the claimed formula
(\ref{lambdaWithE0}).
\end{proof}

\section{Outside the domain\label{outside.sec}}

\subsection{Outline\label{outsideOutline.sec}}

Our goal is to compute the function
\begin{equation}
S_{0}(s,\tau,\lambda):=\lim_{\varepsilon\rightarrow0^{+}}S(s,\tau
,\lambda,\varepsilon), \label{S0def}%
\end{equation}
and its Laplacian with respect to $\lambda,$ using the Hamilton--Jacobi
formulas (\ref{HJfirst}), (\ref{HJsecond1}), and (\ref{HJsecond2}). We wish to
choose the initial conditions $\lambda_{0}$ and $\varepsilon_{0}$---with the
initial momenta given by (\ref{pLambda0}) and (\ref{p0})---so that
$\lambda(\tau)$ is very close to $\lambda$ and $\varepsilon(\tau)$ is very
close to zero. From the formula (\ref{epsilonOfTau}), we see that if
$\varepsilon_{0}$ approaches zero, then $\varepsilon(\tau)$ will also approach zero.

The obstruction to letting $\varepsilon_{0}$ approach 0 is that if
$\varepsilon_{0}=0$ and $\lambda_{0}$ is in the closed support of the measure
$\mu_{s}$ (inside the unit circle), then the integrals (\ref{pLambda0}) and
(\ref{p0}) defining the initial momenta may be singular. Thus, in this
section, we always assume that $\lambda_{0}$ is outside $\operatorname{supp}%
(\mu_{s}).$

For $\lambda_{0}\in\operatorname{supp}(\mu_{s})^{c},$ we have already analyzed
in Section \ref{e0ode.sec} the behavior the solution $\lambda(\tau)$ in the
limit as $\varepsilon_{0}$ approaches zero. Let us then assume that we can
simply set $\varepsilon_{0}=0$ in the second Hamilton--Jacobi formulas
(\ref{HJsecond1}) and (\ref{HJsecond2}). Note that we are not really allowed
to do this, because $S(s,\tau,\lambda,\varepsilon)$ is not known ahead of time
to be smooth up to $\varepsilon=0$; this is a point we will need to return to
in Section \ref{outsideDetails.sec}. If we multiply (\ref{HJsecond1}) by
$\lambda(\tau)$ and set $\varepsilon_{0}=0$, then $\varepsilon(\tau)=0$ and we
obtain%
\[
\lambda(\tau)\frac{\partial S}{\partial\lambda}(s,\tau,\lambda(\tau
),0)=\lambda(\tau)p_{\lambda}(\tau)=\lambda_{0}p_{\lambda,0},
\]
where the second equality is evident from the formulas (\ref{lambdaOfTau}) and
(\ref{pLambdaOfTau}) for $\lambda(\tau)$ and $p_{\lambda}(\tau).$

Now, since $\varepsilon_{0}=0,$ the formula (\ref{lambdaWithE0}) in
Proposition \ref{solutionse0.prop} applies:%
\[
\lambda(\tau)=f_{s-\tau}(\chi_{s}(\lambda_{0})).
\]
Furthermore, Corollary \ref{lambda0lambda.cor} tell us that $f_{s-\tau}%
\circ\chi_{s}$ maps $\operatorname{supp}(\mu_{s})^{c}$ injectively onto
$(\overline{\Sigma}_{s,\tau})^{c}.$ Thus, for $\lambda\in(\overline{\Sigma
}_{s,\tau})^{c},$ we may choose
\begin{equation}
\lambda_{0}=(f_{s-\tau}\circ\chi_{s})^{-1}(\lambda), \label{chooseLambda0}%
\end{equation}
so that $\lambda_{0}$ depends holomorphically on $\lambda.$ Thus, we obtain%
\begin{equation}
\lambda\frac{\partial S_{0}}{\partial\lambda}=\left.  \lambda_{0}p_{\lambda
,0}\right\vert _{\lambda_{0}=(f_{s-\tau}\circ\chi_{s})^{-1}(\lambda)}%
,\quad\lambda\in(\overline{\Sigma}_{s,\tau})^{c}. \label{lambdaDSdLambda}%
\end{equation}

Finally, we note from (\ref{momentum1}) in Proposition \ref{solutionse0.prop}
that $\lambda_{0}p_{\lambda,0}$ depends holomorphically on $\lambda_{0}$, for
$\lambda_{0}$ outside $\operatorname{supp}(\mu_{s}).$ Thus, we conclude from
(\ref{lambdaDSdLambda}) that $\lambda\partial S_{0}/\partial\lambda$ is a
holomorphic function of $\lambda,$ for $\lambda$ outside $\overline{\Sigma
}_{s,\tau}.$ Thus, $\partial S_{0}/\partial\lambda$ is also holomorphic
outside $\overline{\Sigma}_{s,\tau},$ except possibly at the origin. But we
will see, once we compute $\lambda\partial S_{0}/\partial\lambda$ more
explicitly in Theorem \ref{outside.thm1}, that the singularity at the origin
is removable, so $\partial S/\partial\lambda$ is holomorphic on all of
$(\overline{\Sigma}_{s,\tau})^{c}.$ Thus,%
\begin{equation}
\Delta_{\lambda}S_{0}(s,\tau,\lambda)=4\frac{\partial}{\partial\bar{\lambda}%
}\frac{\partial S_{0}}{\partial\lambda}=0 \label{BrownZero}%
\end{equation}
for $\lambda$ outside of $\overline{\Sigma}_{s,\tau},$ showing that the Brown
measure is zero there.

In the next section, we will make the argument sketched here rigorous.

\subsection{Details\label{outsideDetails.sec}}

We use Definitions \ref{sigmas.def} and \ref{def:Sigmastaudef} along with
Theorem \ref{thm:Sigmastau} and the injectivity of $f_{s-\tau}$ (Theorem
\ref{fInjective.thm}) and the injectivity of $\chi_{s}$ (from the relation
$f_{s}(\chi_{s}(z))=z$). These results tell us that $f_{s-\tau}\circ\chi_{s}$
maps $\operatorname{supp}(\mu_{s})^{c}$ injectively onto $(\bar{\Sigma
}_{s,\tau}^{{}})^{c}.$ Our main objective in this section is to rigorize the
argument presented in Section \ref{outsideOutline.sec}, showing that the Brown
measure $\mu_{s,\tau}$ is zero outside $\overline{\Sigma}_{s,\tau}.$

Our first main result is a refinement of (\ref{lambdaDSdLambda}). Recall the
definition (\ref{S0def}) of $S_{0}.$

\begin{theorem}
\label{outside.thm1}For all $\lambda$ outside $\overline{\Sigma}_{s,\tau},$ we
have%
\begin{equation}
\lambda\frac{\partial S_{0}}{\partial\lambda}(s,\tau,\lambda)=\chi_{s-\tau
}(\lambda)\int_{S^{1}}\frac{1}{\chi_{s-\tau}(\lambda)-\xi}~d\mu_{0}(\xi)
\label{SderivOut}%
\end{equation}
and the Brown measure $\mu_{s,\tau}$ is zero outside $\overline{\Sigma
}_{s,\tau}.$ Here, $\chi_{s-\tau}$ is the inverse function to $f_{s-\tau},$
which exists by Definition \ref{def:Sigmastaudef} and Theorem
\ref{fInjective.thm}.
\end{theorem}

The reason that the argument in Section \ref{outsideOutline.sec} is not
rigorous is that we are not allowed to put $\varepsilon_{0}=0$---which results
in $\varepsilon(\tau)=0$ for all $\tau$---in the Hamilton--Jacobi formulas
(\ref{HJfirst}), (\ref{HJsecond1}), and (\ref{HJsecond2}). After all,
$S(s,\tau,\lambda,\varepsilon)$ is initially defined only for $\varepsilon>0$
and is not known ahead of time to have a smooth extension up to $\varepsilon
=0.$

Our strategy will then be to prove that $S(s,\tau,\lambda,\varepsilon)$
\textit{does} have a smooth extension to $\varepsilon=0,$ for $\lambda
\in(\overline{\Sigma}_{s,\tau})^{c}.$ To do this, we look at the map
$(\lambda_{0},\varepsilon_{0})\mapsto(\lambda(\tau),\varepsilon(\tau))$ for a
fixed $\tau.$ This map is defined by solving the Hamilton's equations with
Hamiltonian (\ref{Ham}) and, as we will see, it makes perfect sense to put
$\varepsilon_{0}=0,$ provided $\lambda_{0}$ is outside the closed support of
the measure $\mu_{s}.$ We will show that this map has a smooth inverse at
$(\lambda_{0},0)$ and then use the inverse function theorem and the first
Hamilton--Jacobi formula to construct the smooth extension of $S.$

\begin{proposition}
\label{PhiT.prop}Consider the map $\Psi_{\tau}$ from $\mathbb{C}%
\times(0,\infty)$ into $\mathbb{C}\times(0,\infty)$ given by%
\[
\Psi_{\tau}(\lambda_{0},\varepsilon_{0})=(\lambda(\tau),\varepsilon(\tau)),
\]
where $\lambda(\tau)$ and $\varepsilon(\tau)$ are computed with $\lambda
(0)=\lambda_{0},$ $\varepsilon(0)=\varepsilon_{0},$ and the initial momenta
given by (\ref{pLambda0}) and (\ref{p0}). Then for all $\lambda_{0}%
\in\operatorname{supp}(\mu_{s})^{c},$ the map $\Psi_{\tau}$ extends
analytically to a neighborhood of $\varepsilon_{0}=0,$ and the Jacobian of
$\Psi_{\tau}$ at $(\lambda_{0},0)$ is invertible. Furthermore, the
$\varepsilon$-component of the extended map is an odd function of
$\varepsilon_{0}.$
\end{proposition}

\begin{proof}
As long as $\lambda_{0}$ is outside the closed support of $\mu_{s},$ the
formulas (\ref{pLambda0}) and (\ref{p0}) defining the initial momenta
$p_{\lambda,0}$ and $p_{\varepsilon,0}$ remain well defined and analytic, even
in a neighborhood of $\varepsilon_{0}=0.$ The formulas (\ref{lambdaOfTau}) and
(\ref{epsilonOfTau}) for $\lambda(\tau)$ and $\varepsilon(\tau)$ then depends
analytically on $\lambda_{0}$ and $\varepsilon_{0}.$ It is then easily checked
that all quantities involved are even functions of $\varepsilon_{0}%
$---including $\varepsilon_{0}p_{\varepsilon,0}$---except for the leading
factor of $\varepsilon_{0}$ in the formula for $\varepsilon(\tau),$ making the
$\varepsilon$-component an odd function of $\varepsilon_{0}.$

We then compute the Jacobian of $\Psi_{\tau}$ at $(\lambda_{0},0)$, writing
$\lambda_{0}=x_{0}+iy_{0}$. If we differentiate (\ref{epsilonOfTau}) and then
evaluate at $\varepsilon_{0}=0,$ the leading factor of $\varepsilon_{0}$ in
the formula gives a simple result:%
\begin{equation}
\frac{\partial\varepsilon}{\partial x_{0}}=0;\quad\frac{\partial\varepsilon
}{\partial y_{0}}=0;\quad\frac{\partial\varepsilon}{\partial\varepsilon_{0}%
}=\exp\left\{  \operatorname{Re}\left[  \frac{\tau}{2}(2\lambda_{0}%
p_{\lambda,0}-1)\right]  \right\}  . \label{derivativesOfEpsilon}%
\end{equation}
Meanwhile, when $\varepsilon_{0}=0,$ we may use the formula
(\ref{lambdaWithE0}): $\lambda(\tau)=f_{s-\tau}(\chi_{s}(\lambda_{0})).$ Thus,
the Jacobian of $\Psi_{\tau}$ at $(\lambda_{0},0)$ has the form%
\[%
\begin{pmatrix}
J & \ast\\
0 & \frac{\partial\varepsilon}{\partial\varepsilon_{0}}%
\end{pmatrix}
,
\]
where $J$ is the $2\times2$ Jacobian of the map $f_{s-\tau}\circ\chi_{s}$ and
the value of $\ast$ is irrelevant. Since, by the identity $f_{s}(\chi
_{s}(z))=z$ and Theorem \ref{fInjective.thm}, $f_{s-\tau}\circ\chi_{s}$ is an
injective holomorphic map, its complex derivative cannot be zero and therefore
$J$ must be invertible. Since, also, $\partial\varepsilon/\partial
\varepsilon_{0}$ is positive by (\ref{derivativesOfEpsilon}), we find that the
Jacobian of $\Psi_{\tau}$ is invertible.
\end{proof}

\begin{proposition}
\label{SoutsideReg.prop}Fix $\tau$ with $\left\vert \tau-s\right\vert \leq s$
and a point $\tilde{\lambda}\in(\overline{\Sigma}_{s,\tau})^{c}.$ Then the map%
\[
(\lambda,\varepsilon)\mapsto S(s,\tau,\lambda,\varepsilon),
\]
initially defined for $\varepsilon>0,$ has an analytic extension defined for
$(\lambda,\varepsilon)$ in a neighborhood of $(\tilde{\lambda},0)$ and this
extension is even in $\varepsilon.$
\end{proposition}

\begin{proof}
We define a function $\operatorname{HJ}(\tau,\lambda_{0},\varepsilon_{0})$ by
the right-hand side of the first Hamilton--Jacobi formula (\ref{HJfirst}),
namely
\[
\operatorname{HJ}(s,\tau,\lambda_{0},\varepsilon_{0})=S(s,0,\lambda
_{0},\varepsilon_{0})+2\operatorname{Re}[\tau H_{0}]+\frac{1}{2}%
\operatorname{Re}\left[  \tau(\varepsilon_{0}p_{\varepsilon,0}+2\lambda
_{0}p_{\lambda,0})\right]  ,
\]
where it is understood that the initial momenta $p_{\lambda,0}$ and
$p_{\varepsilon,0}$ are always computed as functions of $\lambda_{0}$ and
$\varepsilon_{0}$ as in (\ref{pLambda0}) and (\ref{p0}). Fix $\tilde{\lambda}$
in $(\overline{\Sigma}_{s,\tau})^{c}.$ By Corollary \ref{lambda0lambda.cor},
we can find $\lambda_{0}$ in $\operatorname{supp}(\mu_{s})^{c}$ so that with
$\varepsilon_{0}=0,$ we have $\lambda(\tau)=\tilde{\lambda}.$ By Proposition
\ref{PhiT.prop} and the inverse function theorem, $\Psi_{\tau}$ has an
analytic inverse near $(\lambda_{0},0).$ Since, by (\ref{epsilonOfTau}),
$\varepsilon(\tau)$ always has the same sign as $\varepsilon_{0},$ we see that
the $\varepsilon_{0}$-component of $\Psi_{\tau}^{-1}(\lambda,\varepsilon)$ is
positive whenever $\varepsilon$ is positive. We then consider the function%
\[
\tilde{S}(s,\tau,\lambda,\varepsilon)=\operatorname{HJ}(s,\tau,\Psi_{\tau
}^{-1}(\lambda,\varepsilon)),
\]
which is an analytic function. Since the first Hamilton--Jacobi formula tells
us that $\tilde{S}$ agrees with $S$ when $\varepsilon>0,$ we see that
$\tilde{S}$ is the desired extension. Since the $\varepsilon$-component of
$\Psi_{\tau}$ is an odd function of $\varepsilon_{0},$ the $\varepsilon_{0}%
$-component of $\Psi_{\tau}^{-1}$ is an odd function of $\varepsilon.$ Since,
as is easily checked, $\mathrm{HJ}$ is an even function of $\varepsilon_{0},$
we see that the extended $S$ is even in $\varepsilon.$
\end{proof}

\begin{proof}
[Proof of Theorem \ref{outside.thm1}]We fix $\lambda$ in $(\overline{\Sigma
}_{s,\tau})^{c}$ and we use choose $\lambda_{0}\in\operatorname{supp}(\mu
_{s})^{c}$ so that with $\varepsilon_{0}=0,$ we get $\lambda(\tau)=\lambda$
and $\varepsilon(\tau)=0.$ (Such a $\lambda_{0}$ exists by Corollary
\ref{lambda0lambda.cor}.) We now multiply the second Hamilton--Jacobi formula
(\ref{HJsecond1}) by $\lambda(\tau),$ using this value of $\lambda_{0}$ and,
initially, $\varepsilon_{0}>0$:%
\[
\lambda(\tau)\frac{\partial S}{\partial\lambda}(s,\tau,\lambda(\tau
),\varepsilon(\tau))=\lambda(\tau)p_{\lambda}(\tau)=\lambda_{0}p_{\lambda,0},
\]
where the second equality is evident from the formulas (\ref{lambdaOfTau}) and
(\ref{pLambdaOfTau}) for $\lambda(\tau)$ and $p_{\lambda}(\tau).$ As
$\varepsilon_{0}$ tends to zero, $\lambda(\tau)$ tends to $\lambda$ and
$\varepsilon(\tau)$ tends to 0. Thus, by the regularity of $S$ established in
Proposition \ref{SoutsideReg.prop}, we may let $\varepsilon_{0}$ tend to zero
and obtain%
\[
\lambda\frac{\partial S_{0}}{\partial\lambda}(s,\tau,\lambda)=\lambda
_{0}p_{\lambda,0}.
\]

We now apply the formula (\ref{momentum2}), from which we can easily solve for
$\lambda_{0}p_{\lambda,0}$, to get%
\begin{equation}
\lambda\frac{\partial S_{0}}{\partial\lambda}(s,\tau,\lambda)=\chi_{s}%
(\lambda_{0})\int_{S^{1}}\frac{1}{\chi_{s}(\lambda_{0})-\xi}~d\mu_{0}(\xi).
\label{lambdaDSlambda2}%
\end{equation}
Finally, we recall that $\lambda=f_{s-\tau}(\chi_{s}(\lambda_{0}))$, so that
(recalling that $\chi_{s}$ is the inverse function to $f_{s}$), $\lambda
_{0}=f_{s}(\chi_{s-\tau}(\lambda)).$ Plugging this value for $\lambda_{0}$
into (\ref{lambdaDSlambda2}) gives the claimed formula (\ref{SderivOut}) for
$\lambda\partial S/\partial\lambda.$ If we divide (\ref{SderivOut}) by
$\lambda$ and note that $\chi_{s-\tau}(0)=0,$ we see that $\partial
S/\partial\lambda$ has a removable singularity at the origin. Thus, $\partial
S/\partial\lambda$ is a holomorphic function of $\lambda\in(\overline{\Sigma
}_{s,\tau})^{c},$ so that as in (\ref{BrownZero}), the Brown measure
$\mu_{s,\tau}$ is zero on $(\overline{\Sigma}_{s,\tau})^{c}.$
\end{proof}

\section{Inside the domain\label{inside.sec}}

We now briefly summarize the way the computation of the Brown measure
\textquotedblleft in the domain\textquotedblright\ (that is, where it is
\textit{not} zero) works in the present paper, since this is the main way our
paper differs from earlier ones. In the present paper, as in earlier works
such as \cite{DHKBrown}, \cite{HZ}, \cite{DemniHamdi}, and \cite{HHadditive},
one initially attempts to achieve the desired condition $\varepsilon(\tau)=0$
by taking $\varepsilon_{0}=0,$ leading to a determination of the region where
the Brown measure is zero. (Section \ref{outside.sec} of the present paper.)
But the way the idea of letting $\varepsilon_{0}$ tend to zero \textit{fails}
to work is different here than in the previous works. In the just-cited
papers, letting $\varepsilon_{0}$ tend to 0 fails at points $\lambda_{0}$
where, with $\varepsilon_{0}=0,$ the solution ceases to exist before one gets
to the $\tau$-value one is interested in. For each such $\lambda_{0}$, one
must look for a \textit{positive} value of $\varepsilon_{0}$ that gives
$\varepsilon(\tau)=0,$ and then compute the associated $\lambda(\tau).$

The situation in the present paper is different in that the lifetime of the
paths is always infinite, even if $\varepsilon_{0}$ tends to zero. (As a
result, we are able to avoid the involved analysis of the lifetime of
solutions in Sections 6.3 and 6.4 of \cite{DHKBrown} and Section 4.3 of
\cite{HZ}.) In this paper, letting $\varepsilon_{0}$ tend to 0 fails
\textit{at points }$\lambda_{0}$\textit{ for which the initial momenta
}$p_{\varepsilon,0}$\textit{ and }$p_{\lambda,0}$\textit{ become ill defined}
in this limit. The set of such points is one dimensional, consisting of the
support of the measure $\mu_{s}$ inside the unit circle.

For points $e^{i\phi}$ in $\operatorname{supp}(\mu_{s}),$ we will let
$\lambda_{0}$ approach $e^{i\phi}$ while \textit{simultaneously} letting
$\varepsilon_{0}$ tend to zero. A key point is that approaching $(e^{i\phi
},0)$ along different paths in $(\lambda_{0},\varepsilon_{0})$-space give
different limiting values of $\lambda(\tau).$ Thus, although
$\operatorname{supp}(\mu_{s})$ is one dimensional, the set of $\lambda(\tau
)$'s obtained will be two dimensional, consisting of the entire domain
$\Sigma_{s,\tau}.$

\subsection{Outline\label{insideOutline.sec}}

Recall the function $\phi^{s}$ introduced in Definition \ref{phiS.def} and the
function $\delta^{s,\tau}$ introduced in Theorem \ref{thm:Sigmastau}. For a
point $\lambda$ in $\Sigma_{s,\tau}$ with $\delta$-coordinate $\delta,$ we let
$\phi^{s,\tau}(\delta)$ be as in Notation \ref{phistau.notation}, which means
that
\[
\phi^{s,\tau}(\delta)=\phi^{s}(\theta^{s,\tau}(\delta)),
\]
where $\theta^{s,\tau}$ is the inverse function to $\delta^{s,\tau}.$ This
definition means simply that $\phi=\phi^{s,\tau}(\delta)$ is related to
$\delta$ as in Figures \ref{thetaphi.fig} and \ref{thetadelta.fig}.

In this section, we provide an outline of the computation of the Brown measure
$\mu_{s,\tau}$ inside the domain $\Sigma_{s,\tau},$ including a brief (but not
entirely rigorous) derivation of the formula (\ref{introFormula}) for the
Brown measure.

\subsubsection{Mapping into the domain\label{mapping.sec}}

Suppose we take $\varepsilon_{0}=0$ and $\lambda_{0}=e^{i\phi}$, where the
density of the measure $\mu_{s}$ at $e^{i\phi}$ is positive. Then the initial
momenta $p_{\lambda,0}$ and $p_{\varepsilon,0}$ become ill defined, because
the integrals (\ref{pLambda0}) and (\ref{p0}) defining them are not absolutely
convergent. Now, if we take $\varepsilon_{0}=0$ and we let $\lambda_{0}$
\textit{approach} such a point $e^{i\phi},$ either from inside or from outside
the unit circle, we get two \textit{different} limiting values of
$\lambda(\tau),$ both lying on the boundary of the domain $\Sigma_{s,\tau}.$
Specifically, the limiting values of $\chi_{s}(\lambda_{0})$ will be
$r_{s}(\theta)e^{i\theta}$ and $r_{s}(\theta)^{-1}e^{i\theta},$ which are the
two circled points on the left-hand side of Figure \ref{thetaphi.fig}; see
\cite[Proposition 2.3 and Proposition 3.7]{Zhong}. The limiting values of
$\lambda(\tau)=f_{s-\tau}(\chi_{s}(\lambda_{0}))$ will then be as on the
right-hand side of Figure \ref{thetadelta.fig}. To get limiting values of
$\lambda(\tau)$ in the \textit{interior} of $\Sigma_{s,\tau},$ we will choose
a family of paths in $(\lambda_{0},\varepsilon_{0})$-space approaching
$(e^{i\phi},0)$ so that the limiting value of $\varepsilon(\tau)$ is always
zero, but limiting value of $\lambda(\tau)$ depends on the path.

We compute the initial momenta in (untwisted) logarithmic coordinates
$\rho=\log\left\vert \lambda\right\vert $ and $\theta=\arg\lambda,$ so that
the initial momenta are $p_{\rho,0},$ $p_{\theta,0},$ and $p_{\varepsilon,0},$
that is, the derivatives of $S(s,0,\lambda_{0},\varepsilon_{0})$ with respect
to $\rho,$ $\theta,$ and $\varepsilon_{0}.$ We will find paths in
$(\lambda_{0},\varepsilon_{0})$-space approaching $(e^{i\phi},0)$ so that: (1)
$\varepsilon_{0}p_{\varepsilon,0}$ always approaches zero, (2) $p_{\theta,0}$
approaches a single number independent of the choice of path, and (3)
$p_{\rho,0}$ can approach a whole range of different values, depending on the
choice of path. Along all paths, the value of $\varepsilon(\tau)$ will
approach 0. We then compute the value of $\lambda(\tau)$ in twisted
logarithmic coordinates and we find that the range of possible values for
$p_{\rho,0}$ gives a range of possible values for the $v$-coordinate of
$\lambda(\tau),$ so that the values of $\lambda(\tau)$ trace out a segment of
an exponential spiral in $\Sigma_{s,\tau}.$ As $\phi$ varies over all points
where the density of $\phi$ is positive, these spiral segments will trace out
the whole domain $\Sigma_{s,\tau}.$ See Section \ref{surjectivity.sec} for
details of this argument.

\subsubsection{Computing the Brown measure\label{computeOutline.sec}}

In this section, we use the twisted logarithmic coordinates $v$ and $\delta$
introduced previously:%
\[
v=\frac{1}{\tau_{1}}\log\left\vert \lambda\right\vert ;\quad\delta=\arg
\lambda-\frac{\tau_{2}}{\tau_{1}}\log\left\vert \lambda\right\vert ,
\]
where $\tau=\tau_{1}+i\tau_{2}.$ We may easily verify the following formulas
for the derivatives of $S$ in those coordinates:%
\begin{align}
\frac{\partial S}{\partial v}  &  =2\operatorname{Re}\left[  \tau~\lambda
\frac{\partial S}{\partial\lambda}\right] \nonumber\\
\frac{\partial S}{\partial\delta}  &  =-2\operatorname{Im}\left[  \lambda
\frac{\partial S}{\partial\lambda}\right]  . \label{SderivTwisted}%
\end{align}
We recall the notation $S_{0}$ for the function%
\[
S_{0}(s,\tau,\lambda)=\lim_{\varepsilon\rightarrow0^{+}}S(s,\tau
,\lambda,\varepsilon).
\]

We now use the second Hamilton--Jacobi formula (\ref{HJsecond1}), which gives
a formula for $\partial S/\partial\lambda$ at the point $(s,\tau,\lambda
(\tau),\varepsilon(\tau)).$ As $(\lambda_{0},\varepsilon_{0})$ approaches
$(e^{i\phi},0)$ along our specially chosen paths, $\varepsilon(\tau)$ tends to
zero and $\lambda(\tau)$ tends to a point $\lambda$ in $\Sigma_{s,\tau}$.
\textit{Assuming} that $S$ has a $C^{1}$ extension to a neighborhood of
$(s,\tau,\lambda,0),$ we can multiply by $\lambda$ and let $\varepsilon_{0}$
tend to zero in the second Hamilton--Jacobi formula to obtain%
\begin{equation}
\lambda\frac{\partial S_{0}}{\partial\lambda}(s,\tau,\lambda)=\lambda
(\tau)p_{\lambda}(\tau)=\lambda_{0}p_{\lambda,0}, \label{logDerivOutline}%
\end{equation}
where the second equality follows directly from the formulas
(\ref{lambdaOfTau}) and (\ref{pLambdaOfTau}).

Meanwhile, let us take the formula (\ref{lambdaOfTau}) for $\lambda(\tau)$ and
evaluate the $v$-coordinates of this point in our limiting case, where
$\varepsilon_{0}=0$ and $\left\vert \lambda_{0}\right\vert =1.$ We obtain%
\begin{align}
v(\tau)  &  =\frac{1}{\tau_{1}}\operatorname{Re}\left[  \frac{\tau}%
{2}(2\lambda_{0}p_{\lambda,0}-1)\right] \nonumber\\
&  =\frac{1}{\tau_{1}}\operatorname{Re}\left[  \frac{\tau}{2}\left(
2\lambda\frac{\partial S_{0}}{\partial\lambda}-1\right)  \right] \nonumber\\
&  =\frac{1}{2\tau_{1}}\frac{\partial S_{0}}{\partial v}-\frac{1}{2},
\label{vOfTau}%
\end{align}
where we have used (\ref{logDerivOutline}) and (\ref{SderivTwisted}). A key
point here is that $v(\tau)$ does not depend on the value of $\lambda_{0},$
provided $\lambda_{0}$ is in the unit circle.

We can solve (\ref{vOfTau}) for $\partial S_{0}/\partial v$ as
\begin{equation}
\frac{\partial S_{0}}{\partial v}=2\tau_{1}v+\tau_{1}. \label{dSdVOutline}%
\end{equation}
There are two remarkable features of this formula: first, that it is very
simple and explicit, and second that (inside the domain) $\partial
S_{0}/\partial v$ depends only on $v$ and not on $\delta$. We can then take
another derivative and obtain
\[
\frac{\partial^{2}S_{0}}{\partial v^{2}}=2\tau_{1}.
\]

We now compute the $\delta$-coordinate of $\lambda(\tau)$ (with $\varepsilon
_{0}=0$ and $\left\vert \lambda_{0}\right\vert =1$). After computing that the
$\delta$-coordinate of $e^{i\theta}e^{\tau z}$ is $\theta+\frac{\left\vert
\tau\right\vert ^{2}}{\tau_{1}}\operatorname{Im}(z),$ we obtain%
\begin{align*}
\delta(\tau)  &  =\arg\lambda_{0}+\frac{\left\vert \tau\right\vert ^{2}}%
{\tau_{1}}\operatorname{Im}\left[  \lambda_{0}p_{\lambda,0}\right] \\
&  =\arg\lambda_{0}-\frac{1}{2}\frac{\left\vert \tau\right\vert ^{2}}{\tau
_{1}}\frac{\partial S_{0}}{\partial\delta},
\end{align*}
from which we obtain%
\begin{equation}
\frac{\partial S_{0}}{\partial\delta}=\frac{2\tau_{1}}{\left\vert
\tau\right\vert ^{2}}(\arg\lambda_{0}-\delta). \label{dSdDeltaOutline}%
\end{equation}
This formula is not as explicit as the formula (\ref{dSdVOutline}) for
$\partial S_{0}/\partial v$ because the result depends on the argument of the
point $\lambda_{0}$ in the unit circle.

Now, as discussed in Section \ref{mapping.sec}, the point $\lambda_{0}$ on the
unit circle is related to the point $\lambda$ in $\Sigma_{s,\tau}$ as in
Figures \ref{thetaphi.fig} and \ref{thetadelta.fig}, so that $\arg\lambda
_{0}=\phi^{s,\tau}(\delta).$ Thus, taking a second derivative, we obtain%
\[
\frac{\partial^{2}S_{0}}{\partial\delta^{2}}=\frac{2\tau_{1}}{\left\vert
\tau\right\vert ^{2}}\left(  \frac{d\phi^{s,\tau}(\delta)}{d\delta}-1\right)
.
\]
We then observe from (\ref{dSdVOutline}) that $\partial S_{0}/\partial v$ is
independent of $\delta,$ so that the mixed derivative $\partial^{2}%
S_{0}/\partial\delta\partial v$ is zero. Using the easily computed formula for
the Laplacian in $(v,\delta)$ coordinates, we obtain%
\begin{align*}
\Delta_{\lambda}S_{0}  &  =\frac{1}{\left\vert \lambda\right\vert ^{2}}\left(
\frac{1}{\tau_{1}^{2}}\frac{\partial^{2}S}{\partial v^{2}}-2\frac{\tau_{2}%
}{\tau_{1}^{2}}\frac{\partial^{2}S}{\partial\delta\partial v}+\frac{\left\vert
\tau\right\vert ^{2}}{\tau_{1}^{2}}\frac{\partial^{2}S}{\partial\delta^{2}%
}\right) \\
&  =\frac{2}{\tau_{1}}\frac{1}{\left\vert \lambda\right\vert ^{2}}\frac{d\phi
}{d\delta}.
\end{align*}

\begin{conclusion}
\label{brown.conclusion}If the preceding argument can be made rigorous, the
Brown measure $\mu_{s,\tau}$ in $\Sigma_{s,\tau}$ is given by%
\[
d\mu_{s,\tau}=\frac{1}{2\pi\tau_{1}}\frac{1}{\left\vert \lambda\right\vert
^{2}}\frac{d\phi}{d\delta}~dx~dy
\]
where $\lambda=x+iy.$ Here $\phi$ and $\delta$ are related as in Notation
\ref{phistau.notation} or Figures \ref{thetaphi.fig} and \ref{thetadelta.fig}.
We therefore recover the claimed formula (\ref{introFormula}) for the Brown measure.
\end{conclusion}

\subsection{Surjectivity\label{surjectivity.sec}}

In this section, we provide the details of the strategy outlined in Section
\ref{mapping.sec}, consisting of letting $\varepsilon_{0}$ approach zero while
\textit{simultaneously} letting $\lambda_{0}$ approach a point $e^{i\phi}$ in
the support of $\mu_{s}.$ By letting $(\lambda_{0},\varepsilon_{0})$ approach
$(e^{i\phi},0)$ along different paths, we will get different limiting values
of $\lambda(\tau),$ while still ensuring that $\varepsilon(\tau)$ approaches
0. The following theorem shows that for any $\lambda$ in $\Sigma_{s,\tau},$ we
can find a point $e^{i\phi}$ and a path in $(\lambda_{0},\varepsilon_{0}%
)$-space approaching $(e^{i\phi},0)$ so that the limiting value of
$\lambda(\tau)$ is $\lambda.$

We use the notation%
\[
\lambda(\tau;\lambda_{0},\varepsilon_{0});\quad\varepsilon(\tau;\lambda
_{0},\varepsilon_{0})
\]
to denote the curves $\lambda(\tau)$ and $\varepsilon(\tau)$ in
(\ref{lambdaOfTau}) and (\ref{epsilonOfTau}) with the indicated initial
conditions and with the initial momenta given by (\ref{pLambda0}) and
(\ref{p0}).

\begin{theorem}
\label{surjectivity.thm}Pick a point $e^{i\phi}$ in $S^{1}$ where the density
of the measure $\mu_{s}$ is positive and choose $\lambda_{0}$ to have the form%
\begin{equation}
\lambda_{0}=(1+c\varepsilon_{0})e^{i\phi},\quad c\in\mathbb{R}.
\label{lambdaOfEpsilon0}%
\end{equation}
Then for all $c\in\mathbb{R},$ we have%
\[
\lim_{\varepsilon_{0}\rightarrow0^{+}}\varepsilon(\tau;(1+c\varepsilon
_{0})e^{i\phi},\varepsilon_{0})=0.
\]
Furthermore, as $c$ varies over $\mathbb{R},$ the value of
\[
\lambda(\tau)=\lim_{\varepsilon_{0}\rightarrow0^{+}}\lambda(\tau
;(1+c\varepsilon_{0})e^{i\phi},\varepsilon_{0})
\]
traces out a curve in $\Sigma_{s,\tau}$ with $\delta$-coordinate fixed and
$v$-coordinate varying between $v_{1}^{s,\tau}(\delta)$ and $v_{2}^{s,\tau
}(\delta).$ As $e^{i\phi}$ varies over all points where the density of
$\mu_{s}$ is positive and $c$ varies over $\mathbb{R},$ the value of
$\lambda(\tau)$ fills out the entire domain $\Sigma_{s,\tau}.$
\end{theorem}

\begin{figure}[ptb]
\centering
\includegraphics[scale=0.55]{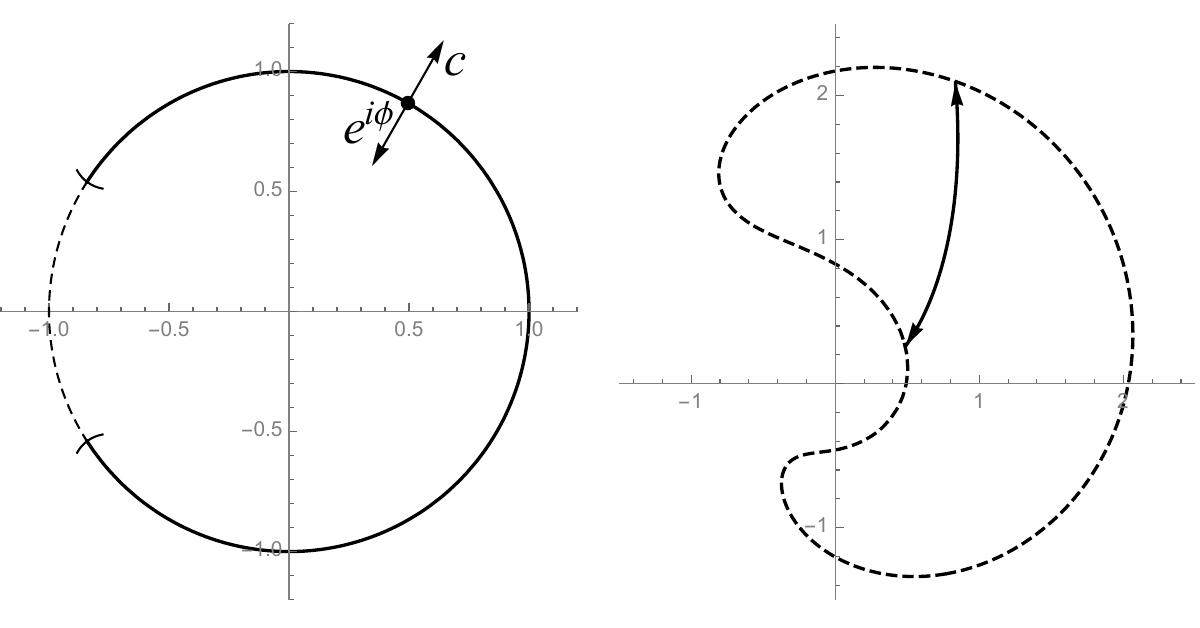}
\caption{A point $e^{i\phi}$ is \textquotedblleft blown up\textquotedblright
\ into an interval by approaching $e^{i\phi}$ along curves of the form
(\ref{lambdaOfEpsilon0}) with different values of $c$ (left). This interval
them maps to a portion of an exponential spiral in $\Sigma_{s,\tau}$ (right).
As $\phi$ varies over points where the density of $m_{s}$ is positive and $c$
varies over $\mathbb{R}$, the entire domain $\Sigma_{s,\tau}$ is swept out.}
\label{maptospiral.fig}
\end{figure}

Theorem \ref{surjectivity.thm} is illustrated in Figure \ref{maptospiral.fig}.%

The key to proving Theorem \ref{surjectivity.thm} is to understand how the
initial momenta---that is, the derivatives of $S(s,\tau,\lambda_{0}%
,\varepsilon_{0})$ at $\tau=0$---behave as $(\lambda_{0},\varepsilon_{0})$
approaches $(e^{i\phi},0)$ along a curve of the form (\ref{lambdaOfEpsilon0}%
).\ If we set $\tau=\varepsilon_{0}=0,$ the resulting function $S_{0}%
(s,0,\lambda_{0})$ is the logarithmic potential of the measure $\mu_{s}$. This
function has a fold-type singularity along the portion of the unit circle
where the density of $\mu_{s}$ is positive.

\begin{figure}[ptb]
\centering
\includegraphics[scale=0.6]{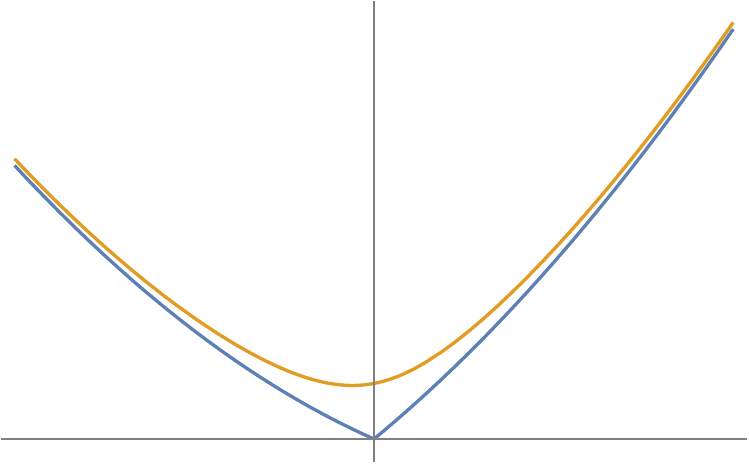}
\caption{Schematic plots of $S_{0}(s,0,e^{\rho+i\phi})$ (bottom) and
$S(s,0,e^{\rho+i\phi},\varepsilon_{0}),$ $\varepsilon_{0}>0$ (top), as
functions of $\rho$ with $\phi$ fixed.}
\label{smoothed.fig}
\end{figure}

The radial derivative of $S_{0}(s,0,\lambda_{0})$ therefore has a jump
discontinuity on portions of the unit circle, while the angular derivative of
$S_{0}(s,0,\lambda_{0})$ is continuous. But as soon as $\varepsilon_{0}$
becomes positive, $S(s,0,\lambda_{0},\varepsilon_{0})$ becomes smooth in
$\lambda_{0},$ so that the radial derivative takes on a continuous range of
values, as shown in Figure \ref{smoothed.fig}. It is then plausible that for
any number $p$ between the \textquotedblleft inner\textquotedblright\ and
\textquotedblleft outer\textquotedblright\ values of $\frac{\partial}%
{\partial\rho}S_{0}(s,0,e^{i\phi}),$ we can find a $\rho$ close to 0 and an
$\varepsilon_{0}$ close to zero for which $\frac{\partial}{\partial\rho
}S(s,0,e^{\rho+i\phi},\varepsilon_{0})$ is close to $p$; this claim is
verified in Proposition \ref{momentumLims.prop}.%

We start by recording what happens when $\varepsilon_{0}=0$ and $\lambda_{0}$
approaches the unit circle.

\begin{proposition}
\label{momentae0.prop}Fix a point $e^{i\phi}$ in the unit circle where the
density of $\mu_{s}$ is positive. Then the limits of $\partial S_{0}%
/\partial\rho(s,0,\lambda_{0})$ as $\lambda_{0}$ approaches $e^{i\phi}$ from
inside and outside the unit circle, respectively, exist and are finite. We
denote these limits as
\[
\frac{\partial S_{0}}{\partial\rho}(s,0,e^{i\phi})^{\operatorname{in}}\text{
and }\frac{\partial S_{0}}{\partial\rho}(s,0,e^{i\phi})^{\operatorname{out}}.
\]
Furthermore the limit of $\partial S_{0}/\partial\theta(s,0,\lambda_{0})$ as
$\lambda_{0}$ approaches $e^{i\phi}$ from points not in the unit circle
exists. We denote this limit as%
\[
\frac{\partial S_{0}}{\partial\theta}(s,0,e^{i\phi}).
\]

\end{proposition}

We will compute these limits and show that, in general, $\partial
S_{0}/\partial\rho(s,0,e^{i\phi})^{\operatorname{in}}$ and $\partial
S_{0}/\partial\rho(s,0,e^{i\phi})^{\operatorname{out}}$ are not equal. We now
turn to the limiting values of the momenta along the curves of the form
(\ref{lambdaOfEpsilon0}).

\begin{proposition}
\label{momentumLims.prop}Fix a point $e^{i\phi}$ in the unit circle where the
density of $\mu_{s}$ is positive. First, we have
\begin{equation}
\lim_{\varepsilon_{0}\rightarrow0^{+}}\frac{\partial S}{\partial\theta
}(s,0,(1+c\varepsilon_{0})e^{i\phi},\varepsilon_{0})=\frac{\partial S_{0}%
}{\partial\theta}(s,0,e^{i\phi}), \label{dSdthetaLim}%
\end{equation}
independent of $c.$ Second, we have%
\begin{equation}
\lim_{\varepsilon_{0}\rightarrow0}\varepsilon_{0}\frac{\partial S}%
{\partial\varepsilon_{0}}(s,0,(1+c\varepsilon_{0})e^{i\phi},\varepsilon_{0})=0
\label{e0dS}%
\end{equation}
independent of $c.$ Third, we have%
\begin{align}
\lim_{\varepsilon_{0}\rightarrow0^{+}}\frac{\partial S}{\partial\rho
}(s,0,(1+c\varepsilon_{0})e^{i\phi},\varepsilon_{0})  &  =\frac{1}{2}\left(
1-\frac{c}{\sqrt{1+c^{2}}}\right)  \frac{\partial S_{0}}{\partial\rho
}(s,0,e^{i\phi})^{\operatorname{in}}\nonumber\\
&  +\frac{1}{2}\left(  1+\frac{c}{\sqrt{1+c^{2}}}\right)  \frac{\partial
S_{0}}{\partial\rho}(s,0,e^{i\phi})^{\operatorname{out}}.
\label{dSdRhoInterpolate}%
\end{align}
Thus, as $c$ varies from $-\infty$ to $\infty,$ the limiting value of
$\partial S/\partial\rho$ varies continuously from $\partial S_{0}%
/\partial\rho(s,0,e^{i\phi})^{\operatorname{in}}$ to $\partial S_{0}%
/\partial\rho(s,0,e^{i\phi})^{\operatorname{out}}.$
\end{proposition}

It is convenient to relate the initial momenta in the preceding propositions
to the Herglotz integral of the measure $\mu_{s}$, namely the integral%
\begin{equation}
\int_{S^{1}}\frac{e^{i\phi^{\prime}}+qe^{i\phi}}{e^{i\phi^{\prime}}-qe^{i\phi
}}~d\mu_{s}(e^{i\phi^{\prime}}), \label{HergMus}%
\end{equation}
for $q\neq1.$ We will first use this idea to verify the statements in
Propositions \ref{momentae0.prop} and \ref{momentumLims.prop} about $\partial
S/\partial\theta.$

\begin{proof}
[Proof of results about $\partial S/\partial\theta$]We try to connect
$\partial S/\partial\theta$ to the imaginary part of the integral in
(\ref{HergMus}) by seeking $q\neq1$ such that%
\begin{equation}
\frac{\partial S}{\partial\theta}(s,0,r_{0}e^{i\phi},\varepsilon
_{0})=\operatorname{Im}\left[  \int_{S^{1}}\frac{e^{i\phi^{\prime}}+qe^{i\phi
}}{e^{i\phi^{\prime}}-qe^{i\phi}}~d\mu_{s}(e^{i\phi^{\prime}})\right]  .
\label{hergIdentity}%
\end{equation}
The left-hand side of (\ref{hergIdentity}) is computed in terms of the measure
$\mu_{s}$ using (\ref{theInitialCondIntegral}). After computing both sides,
(\ref{hergIdentity}) becomes%
\begin{align}
&  \int_{S^{1}}\frac{2\sin(\phi-\phi^{\prime})}{(r_{0}^{2}+1+\varepsilon
_{0}^{2})/r_{0}-2\cos(\phi-\phi^{\prime})}\,d\mu_{s}(e^{i\phi^{\prime}%
})\nonumber\\
&  =\int_{S^{1}}\frac{2\sin(\phi-\phi^{\prime})}{(q+1/q)-2\cos(\phi
-\phi^{\prime})}\,d\mu_{s}(e^{i\phi^{\prime}}), \label{hergIdentity2}%
\end{align}
which gives%
\begin{equation}
q+\frac{1}{q}=\frac{r_{0}^{2}+1+\varepsilon_{0}^{2}}{r_{0}}. \label{qEq}%
\end{equation}

We first consider the case $\varepsilon_{0}=0,$ which gives $q=r_{0}$ or
$q=1/r_{0},$ so that as $r_{0}$ tends to 1, $q$ also tends to 1. Still
assuming $\varepsilon_{0}=0,$ suppose $qe^{i\phi}$ approaches a point on the
unit circle where the density of $\mu_{s}$ is positive, from inside the
circle. Then since the density is analytic where it is positive
\cite[Proposition 3.6]{Zhong}, we can rewrite the integral on the right-hand
side of (\ref{hergIdentity}) as a contour integral and deform the contour to
see that the limit exists. The limit from outside the circle similarly exists
and the two limits are equal because the right-hand side of
(\ref{hergIdentity}), as computed on the right-hand side of
(\ref{hergIdentity2}), is invariant under $q\mapsto1/q.$ This establishes the
claim about $\partial S_{0}/\partial\theta$ in Proposition
\ref{momentae0.prop}.

We next consider the case $r_{0}=1+c\varepsilon_{0}.$ Since the left-hand side
of (\ref{qEq}) is invariant under $q\mapsto1/q,$ the two solutions are
reciprocals of one another, with the solution with $q<1$ simplifying to
\begin{equation}
q=1+\frac{1}{2}\frac{(1+c^{2})\varepsilon_{0}^{2}}{1+c\varepsilon_{0}}%
-\frac{1}{2}\varepsilon_{0}\frac{\sqrt{(1+c^{2})(4+\varepsilon_{0}%
^{2}+4c\varepsilon_{0}+c^{2}\varepsilon_{0}^{2})}}{1+c\varepsilon_{0}},
\label{qOfepsilon0}%
\end{equation}
provided $1+c\varepsilon_{0}>0.$ As $\varepsilon_{0}\rightarrow0^{+}$ with $c$
fixed, $q$ tends to 1. Thus,%
\[
\lim_{\varepsilon_{0}\rightarrow0^{+}}\frac{\partial S}{\partial\theta
}(s,0,(1+c\varepsilon_{0})e^{i\phi},\varepsilon_{0})=\lim_{q\rightarrow1^{-}%
}\operatorname{Im}\left[  \int_{S^{1}}\frac{e^{i\phi^{\prime}}+qe^{i\phi}%
}{e^{i\phi^{\prime}}-qe^{i\phi}}~d\mu_{s}(e^{i\phi^{\prime}})\right]  ,
\]
which is the same limit as we obtained with $\varepsilon_{0}=0$ and $r_{0}$
approaching 1, verifying (\ref{dSdthetaLim}).
\end{proof}

We now turn to $\partial S/\partial\rho.$

\begin{proof}
[Proof of results about $\partial S/\partial\rho$]Since $\partial/\partial
\rho=r\partial/\partial r$, we may compute that%
\begin{align*}
& \frac{\partial S}{\partial\rho}(s,0,r_{0}e^{i\phi},\varepsilon_{0})\\
& =1+\frac{r_{0}^{2}-1-\varepsilon_{0}^{2}}{r_{0}}\int_{S^{1}}\frac{1}%
{(r_{0}^{2}+1+\varepsilon_{0}^{2})/r_{0}-2\cos(\phi-\phi^{\prime})}~d\mu
_{s}(e^{i\phi^{\prime}}).
\end{align*}
Meanwhile, we compute that%
\begin{align}
& \operatorname{Re}\left[  \int_{S^{1}}\frac{e^{i\phi^{\prime}}+qe^{i\phi}%
}{e^{i\phi^{\prime}}-qe^{i\phi}}~d\mu_{s}(e^{i\phi^{\prime}})\right]
\nonumber\\
& =-\left(  q-\frac{1}{q}\right)  \int_{S^{1}}\frac{1}{q+1/q-2\cos(\phi
-\phi^{\prime})}~d\mu_{s}(e^{i\phi^{\prime}}).\label{ReHerg}%
\end{align}
We therefore find that if $q$ satisfies (\ref{qEq}), then
\begin{align}
&  \frac{\partial S}{\partial\rho}(s,0,r_{0}e^{i\phi},\varepsilon
_{0})\nonumber\\
&  =1-\frac{1}{q-1/q}\frac{r_{0}^{2}-1-\varepsilon_{0}^{2}}{r_{0}%
}\operatorname{Re}\left[  \int_{S^{1}}\frac{e^{i\phi^{\prime}}+qe^{i\phi}%
}{e^{i\phi^{\prime}}-qe^{i\phi}}~d\mu_{s}(e^{i\phi^{\prime}})\right]
.\label{HergSrho}%
\end{align}

We start with the case $\varepsilon_{0}=0,$ in which case, we may take
$q=r_{0}$ in (\ref{qEq}). Then (\ref{HergSrho}) reduces to
\begin{align}
&  \frac{\partial S}{\partial\rho}(s,0,r_{0}e^{i\phi},0)\nonumber\\
&  =1-\operatorname{Re}\left[  \int_{S^{1}}\frac{e^{i\phi^{\prime}}%
+r_{0}e^{i\phi}}{e^{i\phi^{\prime}}-r_{0}e^{i\phi}}~d\mu_{s}(e^{i\phi^{\prime
}})\right]  . \label{SrhoAtTau0}%
\end{align}
If we let $r_{0}e^{i\phi}$ approach a point on the unit circle where the
density of $\mu_{s}$ is positive, the limit from inside the disk exists,
giving%
\begin{align}
\frac{\partial S_{0}}{\partial\rho}(s,0,e^{i\phi})^{\operatorname{in}}  &
=1-\lim_{r_{0}\rightarrow1^{-}}\operatorname{Re}\left[  \int_{S^{1}}%
\frac{e^{i\phi^{\prime}}+r_{0}e^{i\phi}}{e^{i\phi^{\prime}}-r_{0}e^{i\phi}%
}~d\mu_{s}(e^{i\phi^{\prime}})\right] \nonumber\\
&  =1-m_{s}(\phi), \label{in}%
\end{align}
where $m_{s}$ is the density of the measure $\mu_{s}$ with respect to the
normalized Lebesgue measure on $S^{1}.$ (The second equality is obtained by
recognizing the right-hand side of (\ref{ReHerg}) as the Poisson integral of
$\mu_{s}$, evaluated at $qe^{i\phi}.$) The limit from outside the disk also
exists but is not equal to the limit from inside, because the right-hand side
of (\ref{ReHerg}) changes sign if $q$ is replaced by $1/q,$ giving%
\begin{align}
\frac{\partial S_{0}}{\partial\rho}(s,0,e^{i\phi})^{\operatorname{out}}  &
=1+\lim_{r_{0}\rightarrow1^{-}}\operatorname{Re}\left[  \int_{S^{1}}%
\frac{e^{i\phi^{\prime}}+r_{0}e^{i\phi}}{e^{i\phi^{\prime}}-r_{0}e^{i\phi}%
}~d\mu_{s}(e^{i\phi^{\prime}})\right] \nonumber\\
&  =1+m_{s}(\phi). \label{out}%
\end{align}

Meanwhile, if we take $r_{0}=1+c\varepsilon_{0}$ in (\ref{HergSrho}) and
choose $q$ as in (\ref{qOfepsilon0}), we may compute that%
\begin{equation}
\frac{1}{q-1/q}\frac{r_{0}^{2}-1-\varepsilon_{0}^{2}}{r_{0}}=\frac
{\varepsilon_{0}-c(2+c\varepsilon_{0})}{\sqrt{(1+c^{2})(4+\varepsilon_{0}%
^{2}+4c\varepsilon_{0}+c^{2}\varepsilon_{0}^{2})}}. \label{bigMess}%
\end{equation}
The right-hand side of this expression tends to $-c/\sqrt{1+c^{2}}$ as
$\varepsilon_{0}$ tends to zero. Meanwhile, as $\varepsilon_{0}$ tends to
zero, the quantity $q$ in (\ref{qOfepsilon0}) tends to 1 from below. Thus,
\begin{equation}
\lim_{\varepsilon\rightarrow0^{+}}\frac{\partial S}{\partial\rho
}(s,0,(1+c\varepsilon_{0})e^{i\phi},\varepsilon_{0})=1+\frac{c}{\sqrt{1+c^{2}%
}}m_{s}(\phi). \label{SrhoFinal}%
\end{equation}
Using this result and (\ref{in}) and (\ref{out}), we can easily check that the
left-hand side of (\ref{dSdRhoInterpolate}) agrees with the right-hand side.
\end{proof}

We finally verify (\ref{e0dS}).

\begin{proof}
We compute that
\[
\frac{\partial S}{\partial\varepsilon_{0}}(s,0,r_{0}e^{i\phi},\varepsilon
_{0})=2\varepsilon_{0}\int_{S^{1}}\frac{1}{(r_{0}^{2}+1+\varepsilon_{0}%
^{2})/r_{0}-2\cos(\phi-\phi^{\prime})}~d\mu_{s}(e^{i\phi^{\prime}}),
\]
so that, if $q$ satisfies (\ref{qEq}), we get%
\begin{equation}
\frac{\partial S}{\partial\varepsilon_{0}}(s,0,r_{0}e^{i\phi},\varepsilon
_{0})=-\frac{\varepsilon_{0}}{q-1/q}\operatorname{Re}\left[  \int_{S^{1}}%
\frac{e^{i\phi^{\prime}}+qe^{i\phi}}{e^{i\phi^{\prime}}-qe^{i\phi}}~d\mu
_{s}(e^{i\phi^{\prime}})\right]  . \label{e0Se0}%
\end{equation}
After checking that $q-1/q=\sqrt{\alpha^{2}+4\alpha},$ we can easily check
that the leading coefficient on the right-hand side of (\ref{e0Se0}) remains
finite as $\varepsilon_{0}$ tends to zero. Thus, after multiplying by
$\varepsilon_{0},$ the limit becomes zero.
\end{proof}

\begin{proof}
[Proof of Theorem \ref{surjectivity.thm}]According to Proposition
\ref{momentumLims.prop}, the quantity $\varepsilon_{0}p_{\varepsilon,0}$ tends
to zero along the path $\lambda_{0}=(1+c\varepsilon_{0})e^{i\phi}$. Meanwhile,
we can easily compute from (\ref{InitialMomentum1}) that
\begin{equation}
2\lambda_{0}p_{\lambda,0}-1=\left.  \left(  \frac{\partial S}{\partial\rho
}-i\frac{\partial S}{\partial\theta}\right)  \right\vert _{\tau=0}-1,
\label{2lambda0plambda0}%
\end{equation}
where the limiting behavior of $\partial S/\partial\rho$ and $\partial
S/\partial\theta$ was worked out in Proposition \ref{momentumLims.prop}. We
then use the formulas (\ref{epsilonOfTau}) and (\ref{lambdaOfTau}) for
$\varepsilon(\tau)$ and $\lambda(\tau).$ From (\ref{epsilonOfTau}), we see
that in the limit, since $\varepsilon_{0}$ is tending to zero and the initial
momenta remain finite, $\varepsilon(\tau)$ tends to zero. We then apply
(\ref{lambdaOfTau}), using (\ref{2lambda0plambda0}) and recalling that
$\varepsilon_{0}p_{\varepsilon,0}$ goes to zero, giving
\begin{align}
\lambda(\tau)  &  =\lambda_{0}\exp\left\{  \frac{\tau}{2}(2\lambda
_{0}p_{\lambda,0}-1)\right\} \nonumber\\
&  =\lambda_{0}e^{-\tau/2}\exp\left\{  -i\frac{\tau}{2}\left.  \frac{\partial
S}{\partial\theta}\right\vert _{\tau=0}\right\}  \exp\left\{  \frac{\tau}%
{2}\left.  \frac{\partial S}{\partial\rho}\right\vert _{\tau=0}\right\}  .
\label{lambdaLim}%
\end{align}
According to Proposition \ref{momentumLims.prop}, the limiting value of
$\partial S/\partial\theta$ is independent of the value of $c,$ while the
value of $\partial S/\partial\rho$ varies between $\partial S/\partial
\rho^{\operatorname{in}}$ and $\partial S/\partial\rho^{\operatorname{out}}.$

If we substitute the formulas (\ref{dSdthetaLim}) and (\ref{dSdRhoInterpolate}%
) from Proposition \ref{momentumLims.prop} into (\ref{lambdaLim}), we see
that
\begin{align}
\lim_{c\rightarrow-\infty}\lambda(\tau)  &  =\lambda_{0}\exp\left\{
\frac{\tau}{2}(2\lambda_{0}p_{\lambda,0}-1)\right\}  ^{\operatorname{in}%
}=f_{s-\tau}(\chi_{s}^{\operatorname{in}}(e^{i\phi}));\nonumber\\
\lim_{c\rightarrow+\infty}\lambda(\tau)  &  =\lambda_{0}\exp\left\{
\frac{\tau}{2}(2\lambda_{0}p_{\lambda,0}-1)\right\}  ^{\operatorname{out}%
}=f_{s-\tau}(\chi_{s}^{\operatorname{out}}(e^{i\phi})),
\label{lambdaEndpoints}%
\end{align}
where $\chi_{s}(e^{i\phi})^{\operatorname{in}}$ and $\chi_{s}(e^{i\phi
})^{\operatorname{out}}$ denote the limiting value of $\chi_{s}(\lambda_{0})$
as $\lambda_{0}$ approaches $e^{i\phi}$ from inside and outside the unit disk,
respectively. These two limiting values of $\lambda(\tau)$ then lie on the
boundary of $\Sigma_{s,\tau}$ and their $\delta$-coordinate is related to
$\phi$ as in Figures \ref{thetaphi.fig} and \ref{thetadelta.fig}. Meanwhile, a
curve of the form $ze^{\tau x}$, where $z$ is fixed and $x$ varies over
$\mathbb{R},$ is an exponential spiral along which the $\delta$-coordinate is
constant and the $v$-coordinate is varying. (Recall the definition
(\ref{coordinates}).) Thus, as $c$ varies with $\phi$ fixed, the value of
$\lambda(\tau)$ in (\ref{lambdaLim}) moves along a segment of an exponential
spiral connecting the two limiting values in (\ref{lambdaEndpoints}). (See
Figure \ref{maptospiral.fig}.) By \cite[Theorem 3.8]{Zhong}, as $\phi$ varies
over all points where the density of $\mu_{s}$ is positive, $\theta$ will vary
over all angles where $r_{s}(\theta)<1,$ so that (Lemma \ref{lem:u1u2}),
$\delta$ will vary over all angles for which $v_{1}^{s,\tau}(\delta
)<v_{2}^{s,\tau}(\delta).$ Thus, as $\phi$ varies over all points where the
density of $\mu_{s}$ is positive, the spiral segments will fill up the whole
of $\Sigma_{s,\tau}$.
\end{proof}

\subsection{Computation of the Brown measure inside the
domain\label{computeInside.sec}}

In this section, we justify the computation of the Brown measure in
$\Sigma_{s,\tau}$ outlined in Section \ref{computeOutline.sec}. The
computation given there is mostly rigorous, except that we need to show
$S(s,\tau,\lambda,\varepsilon)$ extends smoothly up to $\varepsilon=0.$

\subsubsection{Regularity of $S$ at $\varepsilon=0$\label{Sreg.sec}}

In this section, we establish the regularity for $S(s,\tau,\lambda
,\varepsilon)$ near $\varepsilon=0$ that is needed to rigorize the computation
of the Brown measure. We will show that the function%
\[
(\lambda,\varepsilon)\mapsto S(s,\tau,\lambda,\varepsilon),
\]
initially defined for $\varepsilon>0,$ has a $C^{1}$ extension to a
neighborhood of $(\tilde{\lambda},0)$, for any $\tilde{\lambda}$ in
$\Sigma_{s,\tau}.$ (See Theorem \ref{Sregularity.thm} below.)

We consider points $(\lambda_{0},\varepsilon_{0})$, where at first we require
$\varepsilon_{0}>0.$ We label such points using the coordinates
\[
\varepsilon_{0};\quad\phi=\arg(\lambda_{0});\quad c=\frac{\left\vert
\lambda_{0}\right\vert -1}{\varepsilon_{0}}.
\]
Now, if we let $\varepsilon_{0}$ tend to zero with a fixed value of $c,$ the
value $\left\vert \lambda_{0}\right\vert $ will approach 1 and $\lambda_{0}$
itself will approach $e^{i\phi}.$ In this limit, $\lambda_{0}$ is approaching
a single point $e^{i\phi},$ independent of $c$---and yet the limiting value of
$\lambda(\tau),$ computed in the previous subsection, does depend on $c$
(Theorem \ref{surjectivity.thm}). It is therefore natural to consider a
\textquotedblleft blown up\textquotedblright\ domain in which we consider all
non-negative values of $\varepsilon_{0}$ and all real values of $c.$ The point
here is that even when $\varepsilon_{0}=0,$ we count points with the same
value of $\phi$ but different values of $c$ as being different.

\begin{figure}[ptb]
\centering
\includegraphics[scale=0.7]{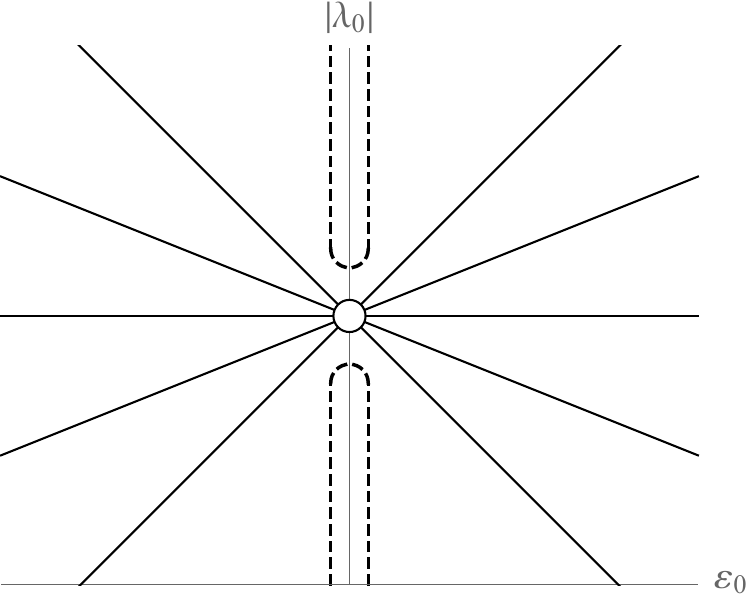}
\caption{For each $\phi,$ points in the $(\varepsilon_{0},\left\vert
\lambda_{0}\right\vert )$ plane obtained by approaching $(0,1)$ along lines
with dfferent slope $c$ are considered to be different.}
\label{blowup.fig}
\end{figure}

To establish the desired $C^{1}$ extension of $S,$ we need to allow the
variable $\varepsilon_{0}$ to be negative. Our blown up domain then consists
of triples $(\varepsilon_{0},c,\phi)$, where $\phi$ is an angle and
$\varepsilon_{0}$ and $c$ are real numbers with $1+c\varepsilon_{0}>0$. We may
then define $\lambda_{0}$ in terms of $(\varepsilon_{0},c,\phi)$ as
$\lambda_{0}=(1+c\varepsilon_{0})e^{i\phi}.$ Note that if $\varepsilon_{0}=0,$
then $\lambda_{0}=e^{i\phi}$; hence, $\left\vert \lambda_{0}\right\vert $ must
equal 1 when $\varepsilon_{0}=0.$ The blown up domain can be understood
geometrically as in Figure \ref{blowup.fig}. The variable $c$ is the slope of
a line in the $(\varepsilon_{0},\left\vert \lambda_{0}\right\vert )$ plane
passing through the point $(0,1)$ and we decree that we obtain different
points if we approach $(0,1)$ along lines lines of different slope.%

The first main result of this section is that the map sending $\lambda_{0}$
and $\varepsilon_{0}$ to $\lambda(\tau)$ and $\varepsilon(\tau),$ viewed as a
map on the blown up domain, has a $C^{1}$ extension with nonsingular Jacobian
up to (and a little beyond) $\varepsilon_{0}=0.$ Once this result is
established, we will use the first Hamilton--Jacobi formula to obtain the
desired $C^{1}$ extension of $S.$

\begin{theorem}
\label{mapRegularity.thm}Fix $\tau$ and consider the map%
\[
(\lambda_{0},\varepsilon_{0})\mapsto(\lambda(\tau),\varepsilon(\tau)),
\]
initially defined for $\lambda_{0}\in\mathbb{C}$ and $\varepsilon_{0}>0.$ We
compute using the coordinates $(\varepsilon_{0},\phi,c)$ on the domain side
and $(\varepsilon,\delta,v)$ on the range side, giving a map%
\[
\Psi_{\tau}(\varepsilon_{0},\phi,c)=(\varepsilon,\delta,v),
\]
initially defined for $\varepsilon_{0}>0,$ $\phi\in\lbrack-\pi,\pi),$ and
$c\in\mathbb{R}.$ Fix a real number $\tilde{c}$ and a point $\tilde{\phi}$
where the density $m_{s}$ of the measure $\nu_{s}$ is positive. Then (1)
$\Psi_{\tau}$ extends to a $C^{1}$ function defined for $(\varepsilon_{0}%
,\phi,c)$ in a neighborhood of $(0,\tilde{\phi},\tilde{c})$ inside
$\mathbb{R}^{3}$, and (2) the Jacobian of $\Psi_{\tau}$ at $(0,\tilde{\phi
},\tilde{c})$ is invertible.
\end{theorem}

We separate the proof into two parts, the proof of smoothness and the
computation of the Jacobian.

\begin{proof}
[Proof (smoothness)]The formulas for $\lambda(\tau)$ and $\varepsilon(\tau)$
in (\ref{lambdaOfTau}) and (\ref{epsilonOfTau}) are smooth functions of
$\lambda_{0},$ $\varepsilon_{0},$ and the quantities $\varepsilon
_{0}p_{\varepsilon,0}$ and $\lambda_{0}p_{\lambda,0}$. Since $\lambda
_{0}=(1+c\varepsilon_{0})e^{i\phi}$ is a smooth function of $(\varepsilon
_{0},\phi,c),$ it suffices to verify that $\varepsilon_{0}p_{\varepsilon,0}$
and $\lambda_{0}p_{\lambda,0}$ are $C^{1}$ functions of $(\varepsilon_{0}%
,\phi,c).$ In light of (\ref{2lambda0plambda0}), it then suffices to verify
smoothness of $\varepsilon_{0}p_{\varepsilon,0}$, $\partial S(s,0,\lambda
_{0},\varepsilon_{0})/\partial\rho$, and $\partial S(s,0,\lambda
_{0},0)/\partial\theta$ as functions of $\varepsilon_{0},$ $\phi,$ and $c.$ To
obtain the desired smoothness, we use the expressions (\ref{e0Se0}),
(\ref{HergSrho}) (together with (\ref{bigMess})), and (\ref{hergIdentity}). We
then use the previously noted \textquotedblleft deform the
contour\textquotedblright\ argument to show that all these expressions extend
smoothly as functions of $q$ past $q=1.$ Finally, $q$ itself, given as a
function of $c$ and $\varepsilon_{0}$ by (\ref{qOfepsilon0}), has a smooth
extension past $\varepsilon_{0}=0.$
\end{proof}

\begin{proof}
[Proof (Jacobian)]From (\ref{lambdaOfTau}) and (\ref{epsilonOfTau}), we have%
\begin{equation}
\varepsilon=\varepsilon_{0}\frac{\left\vert \lambda\right\vert }{\left\vert
\lambda_{0}\right\vert }=\varepsilon_{0}\frac{\left\vert \lambda\right\vert
}{1+c\varepsilon_{0}}. \label{zz0}%
\end{equation}
Since we have shown that $\left\vert \lambda\right\vert $ is a smooth function
of $(\varepsilon_{0},\phi,c),$ the leading factor of $\varepsilon_{0}$ in
(\ref{zz0}) easily gives that
\[
\left.  \frac{\partial\varepsilon}{\partial\varepsilon_{0}}\right\vert
_{\varepsilon_{0}=0}=\left\vert \lambda\right\vert ;\quad\left.
\frac{\partial\varepsilon}{\partial\phi}\right\vert _{\varepsilon_{0}%
=0}=\left.  \frac{\partial\varepsilon}{\partial c}\right\vert _{\varepsilon
_{0}=0}=0.
\]

Let us now evaluate $\lambda(\tau)$ at $\varepsilon_{0}=0$---that is, using
fixed values of $c$ and $\phi$ and letting $\varepsilon_{0}$ tend to zero as
in Proposition \ref{momentumLims.prop}. Let us use the superscript
\textquotedblleft$\operatorname{in}$\textquotedblright\ to mean that we
evaluate the relevant quantity with $\varepsilon_{0}=0$ in the limit as
$\lambda_{0}$ approaches $e^{i\phi}$ from \textit{inside} the unit circle.
Then, using (\ref{SrhoFinal}), we obtain
\[
\left.  \lambda(\tau)\right\vert _{\varepsilon_{0}=0}=\lambda_{0}\exp\left\{
\frac{\tau}{2}(2\lambda_{0}p_{\lambda,0}-1)^{\operatorname{in}}\right\}
\exp\left\{  \frac{\tau}{2}\left(  1+\frac{c}{\sqrt{1+c^{2}}}\right)
m_{s}(\phi)\right\}  ,
\]
where $m_{s}$ is the density of $\phi_{s},$ or, by Proposition
\ref{solutionse0.prop},%
\begin{equation}
\lambda(\tau)=f_{s-\tau}(\chi_{s}(\lambda_{0}))^{\operatorname{in}}%
\exp\left\{  \frac{\tau}{2}\left(  1+\frac{c}{\sqrt{1+c^{2}}}\right)
m_{s}(\phi)\right\}  . \label{lambdae0}%
\end{equation}

We see, then, that the factor of $f_{s-\tau}(\chi_{s}(\lambda_{0}%
))^{\operatorname{in}}$ depends only on $\phi$ and not on $c.$ Furthermore,
because the final exponential factor in (\ref{lambdae0}) has the form
$e^{a\tau}$ for some real number $a,$ this factor only affects the value of
$v$ and not the value of $\delta.$ Thus, the value of $\delta$ comes from
$f_{s-\tau}(\chi_{s}(\lambda_{0}))^{\operatorname{in}},$ which does does not
depend on $c$, and we see that $\partial\delta/\partial c=0$ and we find,
using (\ref{zz0}) that the Jacobian at $\varepsilon_{0}=0$ has the form%
\[
\left(
\begin{array}
[c]{ccc}%
\left\vert \lambda\right\vert  & 0 & 0\\
\ast & \frac{\partial\delta}{\partial\phi} & 0\\
\ast & \ast & \frac{\partial v}{\partial c}%
\end{array}
\right)  .
\]

Now, since the value of $\delta$ comes from $f_{s-\tau}(\chi_{s}(\lambda
_{0}))^{\operatorname{in}},$ we see that $\delta$ depends on $\phi$ as in
Figures \ref{thetaphi.fig} and \ref{thetadelta.fig}. Thus, by Proposition
\ref{prop:dpsidtheta} and the estimate \cite[Lemma 4.20]{HZ} for
$d\phi/d\theta,$ we see that $\partial\delta/\partial\phi$ is positive.
Furthermore, $\partial v/\partial c$ is easily computed (from the final
exponential factor in (\ref{lambdae0})) to be positive, provided we evaluate
at $\phi=\tilde{\phi},$ with $m_{s}(\tilde{\phi})>0.$ Finally, we recall that
the point $\lambda=0$ is never in $\Sigma_{s,\tau}$ and we can see that the
determinant of the Jacobian is positive.
\end{proof}

\begin{theorem}
\label{Sregularity.thm}Pick a point $\tilde{\lambda}\in\Sigma_{s,\tau}.$ Then
the map%
\[
(\lambda,\varepsilon)\mapsto S(s,\tau,\lambda,\varepsilon),
\]
initially defined for $\varepsilon>0,$ has a $C^{1}$ extension defined for
$(\lambda,\varepsilon)$ in a neighborhood of $(\tilde{\lambda},0).$
\end{theorem}

We emphasize that for $\varepsilon<0,$ the $C^{1}$ extension of $S(s,\tau
,\lambda,\varepsilon)$ is \textit{not} given by the formula (\ref{Sfunction}).
Indeed, since $\varepsilon(\tau)p_{\varepsilon}(\tau)$ is independent of
$\tau$ (from the formulas (\ref{epsilonOfTau}) and (\ref{pEpsilonOfTau})), we
can compute $p_{\varepsilon}(\tau)$ as $p_{\varepsilon,0}\varepsilon
_{0}/\varepsilon(\tau).$ If $\lambda_{0}=(1+c\varepsilon_{0})e^{i\phi}$, with
$m_{s}(\phi)>0,$ then as $\varepsilon_{0}$ tends to zero, both $p_{\varepsilon
,0}$ and $\varepsilon_{0}/\varepsilon(\tau)$ tend to nonzero limits. (For
$p_{\varepsilon,0}$, use (\ref{e0Se0}) and for $\varepsilon_{0}/\varepsilon
(\tau)$, use (\ref{epsilonOfTau}).) Using the second Hamilton--Jacobi formula
(\ref{HJsecond2}), we find that $\partial S/\partial\varepsilon$ evaluated at
$\varepsilon=0$ is nonzero, so that a $C^{1}$ extension of $S$ cannot be an
even function of $\varepsilon.$ To put it differently, if we \textit{did}
define $S$ by the same formula (\ref{Sfunction}) even for $\varepsilon<0,$
then $S$ would \textit{not} be $C^{1}$ but would have a singularity of
absolute-value type as a function of $\varepsilon$ near $\varepsilon=0.$

\begin{corollary}
\label{S0pde.cor}Consider the function $S_{0}$ given by $S_{0}(s,\tau
,\lambda)=\lim_{\varepsilon\rightarrow0^{+}}S(s,\tau,\lambda,\varepsilon).$
Then for $\tau$ satisfying $\left\vert \tau-s\right\vert <s$ and $\lambda$ in
$\Sigma_{s,\tau},$ the function $S_{0}$ satisfies the PDE obtained by formally
setting $\varepsilon=0$ in (\ref{thePDE}), namely%
\begin{equation}
\frac{\partial S_{0}}{\partial\tau}=-\frac{1}{2}\left(  \lambda^{2}\left(
\frac{\partial S_{0}}{\partial\lambda}\right)  ^{2}-\lambda\frac{\partial
S_{0}}{\partial\lambda}\right)  . \label{S0PDE}%
\end{equation}
Furthermore, suppose we define curves $\lambda(\tau)$ and $p_{\lambda}(\tau)$
by taking the curves in (\ref{lambdaOfTau})\ and (\ref{pLambdaOfTau}) with
$\lambda_{0}=(1+c\varepsilon_{0})e^{i\phi}$ as in (\ref{lambdaOfEpsilon0}%
)\ and letting $\varepsilon_{0}$ tend to zero. Then $S_{0}$ satisfies the
first Hamilton--Jacobi formula%
\begin{equation}
S_{0}(s,\tau,\lambda(\tau))=S_{0}(s,0,\lambda_{0})+2\operatorname{Re}[\tau
H_{0}]+\operatorname{Re}\left[  \tau\lambda_{0}p_{\lambda,0}\right]
\label{S0HJfirst}%
\end{equation}
and the second Hamilton--Jacobi formula%
\begin{equation}
\frac{\partial S_{0}}{\partial\lambda}(s,\tau,\lambda(\tau))=p_{\lambda}%
(\tau), \label{S0HJsecond}%
\end{equation}
where $H_{0}$ is the value of the Hamiltonian $H$ in (\ref{compexHamiltonian})
with $\varepsilon p_{\varepsilon}=0$ and $\lambda p_{\lambda}$ replaced by
$\lambda_{0}p_{\lambda,0}.$
\end{corollary}

\begin{proof}
Theorem \ref{Sregularity.thm} justifies formally setting $\varepsilon=0$ in
the PDE. Using Theorem \ref{Sregularity.thm} again along with Proposition
\ref{momentumLims.prop} allow us to let $\varepsilon_{0}$ tend to zero in the
Hamilton--Jacobi formulas (\ref{HJfirst}), (\ref{HJsecond1}), and
(\ref{HJsecond2}).
\end{proof}

\begin{proof}
[Proof of Theorem \ref{Sregularity.thm}]We define a function
$\operatorname{HJ}(s,\tau,\lambda_{0},\varepsilon_{0})$ by the right-hand side
of the first Hamilton--Jacobi formula, namely%
\[
\operatorname{HJ}(s,\tau,\lambda_{0},\varepsilon_{0})=S(s,0,\lambda
_{0},\varepsilon_{0})+2\operatorname{Re}[\tau H_{0}]+\frac{1}{2}%
\operatorname{Re}\left[  \tau(\varepsilon_{0}p_{\varepsilon,0}+2\lambda
_{0}p_{\lambda,0})\right]
\]
where it is understood that the initial momenta $p_{\varepsilon,0}$ and
$p_{\lambda,0}$ are always computed as functions of $\lambda_{0}$ and
$\varepsilon_{0}$ as in (\ref{pLambda0}) and (\ref{p0}). We express
$(\lambda_{0},\varepsilon_{0})$ in terms of the coordinates $(\varepsilon
_{0},c,\phi)$ and work in the blown up domain. Then the proof of Theorem
\ref{mapRegularity.thm} shows that $\varepsilon_{0}p_{\varepsilon,0}$ and
$\lambda_{0}p_{\lambda,0}$ extend to $C^{1}$ functions of these variables
defined near $(0,\tilde{c},\tilde{\phi}).$ A similar argument shows that
$S(s,0,\lambda_{0},\varepsilon_{0}),$ which is the regularized log potential
of the measure $\mu_{s},$ also has a $C^{1}$ extension.

We now appeal to Theorem \ref{mapRegularity.thm} and the inverse function
theorem to construct an inverse to the map $\Psi_{\tau}$ near $(0,\delta,v).$
Since the extended map $\Psi_{\tau}$ uses the same formula (\ref{epsilonOfTau}%
) for $\varepsilon(\tau),$ but with $C^{1}$ extensions of the momenta, we see
that the sign of $\varepsilon$ is always the same as the sign of
$\varepsilon_{0}.$ Thus, the $\varepsilon_{0}$-component of $\Psi_{\tau}%
^{-1}(\varepsilon,\delta,v)$ is positive when $\varepsilon$ is positive. We
then consider%
\[
\tilde{S}(s,\tau,\lambda,\varepsilon)=\operatorname{HJ}(s,\tau,\Psi_{\tau
}^{-1}(\varepsilon,\delta,v)),
\]
which is a $C^{1}$ function. The first Hamilton--Jacobi formula tells us that
$\tilde{S}$ agrees with $S$ when $\varepsilon>0.$ Thus, $\tilde{S}$ is the
desired extension.
\end{proof}

\subsubsection{The Brown measure\label{theBrownMeasure.sec}}

We now rigorize the argument outlined in Section \ref{computeOutline.sec}.

\begin{proof}
[Proof of Conclusion \ref{brown.conclusion}]Choose an arbitrary point
$\lambda$ in $\Sigma_{s,\tau}.$ By Theorem \ref{surjectivity.thm}, we can find
$\phi$ and $c$ so that if $\lambda_{0}=(1+c\varepsilon_{0})e^{i\phi},$ then as
$\varepsilon_{0}$ tends to zero, $\varepsilon(\tau)$ tends to zero and
$\lambda(\tau)$ tends to $\lambda.$ As shown in the proof of Theorem
\ref{surjectivity.thm}, the $\delta$-coordinate of $\lambda$ is then related
to $\phi$ as in Figures \ref{thetaphi.fig} and \ref{thetadelta.fig}. Thus, as
in Section \ref{computeOutline.sec}, the limiting value of $\arg\lambda_{0}$
is $\phi^{s,\tau}(\delta).$

We then apply the second Hamilton--Jacobi formula (\ref{HJsecond1}) with this
choice of $\lambda_{0},$ where we initially require $\varepsilon_{0}>0$. After
multiplying by $\lambda,$ we have%
\begin{equation}
\lambda\frac{\partial S}{\partial\lambda}(s,\tau,\lambda(\tau),\varepsilon
(\tau))=\lambda(\tau)p_{\lambda}(\tau)=\lambda_{0}p_{\lambda,0},
\label{SuLimit}%
\end{equation}
where the second equality follows from the formulas (\ref{lambdaOfTau}) and
(\ref{pLambdaOfTau}) for $\lambda(\tau)$ and $p_{\lambda}(\tau).$ Theorem
\ref{Sregularity.thm} allows us to let $\varepsilon_{0}$ tend to zero on the
left-hand side of (\ref{SuLimit}) to obtain $\frac{\partial}{\partial\lambda
}S_{0}(s,\tau,\lambda).$ Meanwhile, Proposition \ref{momentumLims.prop} tells
us that as $\varepsilon_{0}$ tends to zero, $\lambda_{0}p_{\lambda,0}$
approaches some finite number, which we continue to call $\lambda
_{0}p_{\lambda,0},$ and we obtain the claimed formula (\ref{logDerivOutline}).
The computation of $\partial S_{0}/\partial v$ and $\partial S_{0}%
/\partial\delta$ then proceeds as in Section \ref{computeOutline.sec}, leading
to the formulas (\ref{dSdVOutline}) and (\ref{dSdDeltaOutline}) for $\partial
S_{0}/\partial v$ and $\partial S_{0}/\partial\delta.$

From these formulas, we can see that the partial derivatives of $S_{0}%
(s,\tau,\lambda)$ are themselves continuously differentiable. We may therefore
compute the Laplacian of $S_{0}(s,\tau,\lambda)$ in the pointwise sense, and
the rest of the derivation in Section \ref{computeOutline.sec} may be applied
without change.
\end{proof}

We have already shown (Theorem \ref{outside.thm1}) that the Brown measure
$\mu_{s,\tau}$ is zero outside the closure of $\Sigma_{s,\tau}$ and we have
computed $\mu_{s,\tau}$ on the open set $\Sigma_{s,\tau}.$ It remains only to
show that $\mu_{s,\tau}$ does not give mass to the boundary of $\Sigma
_{s,\tau}.$

\begin{theorem}
\label{fullMass.thm}The Brown measure $\mu_{s,\tau}$ assigns measure 1 to the
open set $\Sigma_{s,\tau}.$ Thus, $\mu_{s,\tau}$ has no mass on the boundary
of $\Sigma_{s,\tau}$.
\end{theorem}

Intuitively, the mass of the boundary should be zero if the limiting values of
$S_{0}(s,\tau,\lambda)$ and its first derivatives are the same when
approaching $\partial\Sigma_{s,\tau}$ from inside $\Sigma_{s,\tau}$ as when
approaching $\partial\Sigma_{s,\tau}$ from outside $\overline{\Sigma}_{s,\tau
}.$ And, indeed, these limiting values are the same. This claim can be
verified by using the first and second Hamilton--Jacobi formulas
(\ref{HJfirst}), (\ref{HJsecond1}), and (\ref{HJsecond2}) and checking that
$\lambda_{0}$ and $p_{\lambda,0}$ can be computed as continuous functions of
$\lambda$ over all of $\mathbb{C}$, even when $\lambda$ is in the boundary of
$\Sigma_{s,\tau}.$ (Check that the \textquotedblleft inside\textquotedblright%
\ values of $\lambda_{0}$ and $p_{\lambda,0}$ in Section
\ref{surjectivity.sec} agree with the \textquotedblleft
outside\textquotedblright\ values of $\lambda_{0}$ and $p_{\lambda,0}$ in
Section \ref{outside.sec} when $\lambda$ approaches a point on the boundary of
$\Sigma_{s,\tau}.$)

If we then integrate $S_{0}(s,\tau,\lambda)$ against $\Delta\chi(\lambda)$ for
some test function $\lambda,$ we may split the integral into integrals over
$\Sigma_{s,\tau}$ and its complement. We may then hope to integrate by parts
using Green's second identity. Since $S_{0}$ and its first derivatives are
continuous up to the boundary, the boundary terms from inside and from outside
should cancel, giving the integral of $(\Delta S_{0})\chi$ over $\Sigma
_{s,\tau}$ and the integral of $(\Delta S_{0})\chi$ over $\Sigma_{s,\tau}%
^{c}.$ Since the \textquotedblleft outside\textquotedblright\ integral is zero
by Theorem \ref{outside.thm1}, we would be left with the integral of $(\Delta
S_{0})\chi$ over $\Sigma_{s,\tau},$ showing that the Brown measure is
supported on $\Sigma_{s,\tau}.$ This argument is carried out the case $\tau=s$
and $\mu_{0}=\delta_{1}$ in \cite[Proposition 7.13]{DHKBrown}.

The difficulty with carrying this argument out in general is that we would
have to justify the use of Green's second identity, using, say, Federer's
theorem \cite[Section 12.2]{morgan}. But to use Federer's theorem, we would
need to know that the boundary $\Sigma_{s,\tau}$ is a rectifiable curve, which
we do not know if $\mu_{0}$ is a completely arbitrary probability measure on
the unit circle. (By contrast, if $\mu_{0}=\delta_{1},$ as in \cite{DHKBrown},
then the boundary of $\Sigma_{s,\tau}$ is smooth, except possibly at a single point.)

We will therefore use an indirect approach, similar to the one developed in
\cite[Lemma 4.25]{HZ} and used also in \cite[Theorem 7.9]{HHadditive}. In
\cite[Proposition 4.31]{HZ}, it is shown that when $\tau=s,$ the Brown measure
$\mu_{s,s}$ gives mass one to the domain $\Sigma_{s,s}=\Sigma_{s}.$ In Section
\ref{push.sec}, we will use a push-forward result to prove that $\mu_{s,\tau}$
assigns mass 1 to $\Sigma_{s,\tau},$ for all nonzero $\tau$ with $\left\vert
\tau-s\right\vert \leq s.$ See Section \ref{fullMass.sec} for the details of
this argument.

\section{Relating different values of $\tau$\label{push.sec}}

The goal of this section is to prove Theorem \ref{pushforward.thm}, which
states that the Brown measure $\mu_{s,\tau}$ is the push-forward of $\mu
_{s,s}$ under a certain map $\Phi_{s,\tau}:\overline{\Sigma}_{s}%
\rightarrow\overline{\Sigma}_{s,\tau}.$ Recall the notation%
\[
S_{0}(s,\tau,\lambda)=\lim_{\varepsilon\rightarrow0^{+}}S(s,\tau
,\lambda,\varepsilon),
\]
so that the Brown measure $\mu_{s,\tau}$ is the distributional Laplacian of
$S_{0}(s,\tau,\lambda)$ with respect to $\lambda.$

\subsection{The maps}

Throughout this section, we assume that $\mu_{0}$ (the law of the unitary
element $u$ in the expression $ub_{s,\tau}$) is fixed. Recall the
\textquotedblleft characteristic curves\textquotedblright\ $\lambda(\tau)$
defined in (\ref{lambdaOfTau}). We will define the map $\Phi_{s,\tau}$ in
terms of these curves and then compute the map more explicitly in Proposition
\ref{domainMap.prop}.

\begin{definition}
\label{PhiMap.def}Fix a positive real number $s$ and a nonzero complex number
$\tau$ such that $\left\vert \tau-s\right\vert \leq s$. For each $\lambda
\in\overline{\Sigma}_{s},$ define $\lambda_{0}$ in the unit circle by%
\[
\lambda_{0}=\exp\left\{  i\phi^{s}(\theta)\right\}  ,
\]
where $\theta=\arg\lambda$ and $\phi^{s}$ is as in Definition \ref{phiS.def}.
Let $p_{\theta,0}=\partial S_{0}/\partial\theta(0,0,\lambda_{0})$ be computed
from $\lambda_{0}$ as in Proposition \ref{momentae0.prop}. By Theorem
\ref{surjectivity.thm}, we can choose $p_{\rho,0}$ between $\partial
S_{0}/\partial\rho^{\operatorname{in}}$ and $\partial S_{0}/\partial
\rho^{\operatorname{out}}$ so that with initial conditions $(\lambda
_{0},p_{\theta,0},p_{\rho,0})$---and with $\varepsilon_{0}p_{\varepsilon,0}%
=0$---we have%
\begin{equation}
\lambda(s)=\lambda. \label{lambdaTau1}%
\end{equation}
Then define a map $\Phi_{s,\tau}:\overline{\Sigma}_{s}\rightarrow
\overline{\Sigma}_{s,\tau}$ by setting
\[
\Phi_{s,\tau}(\lambda)=\lambda(\tau),
\]
where $\lambda(\tau)$ is computed using \textbf{the same initial conditions}
$(\lambda_{0},p_{\theta,0},p_{\rho,0})$ (with $\varepsilon_{0}p_{\varepsilon
,0}=0$) that give (\ref{lambdaTau1}).
\end{definition}

\begin{figure}[ptb]
\centering
\includegraphics[scale=0.9]{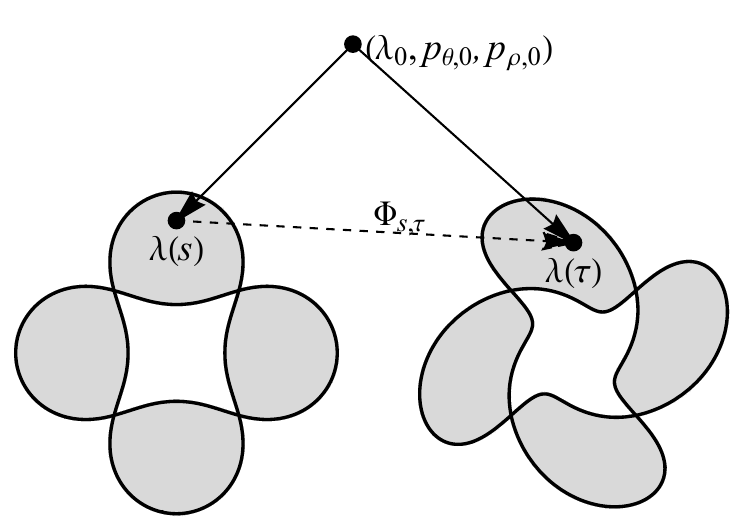}
\caption{To compute $\Phi_{s,\tau}(\lambda),$ we choose initial conditions
$(\lambda_{0},p_{\theta,0},p_{\rho,0})$ so that $\lambda(s)=\lambda$ and then
set $\Phi_{s,\tau}(\lambda)$ equal to $\lambda(\tau),$ computed with the same
initial conditions.}
\label{mapdef.fig}
\end{figure}

See Figure \ref{mapdef.fig}. The following is the main result about the maps
$\Phi_{s,\tau}.$

\begin{theorem}
\label{pushGeneral.thm}Fix a positive real number $s$ and a nonzero complex
number $\tau$ such that $\left\vert \tau-s\right\vert \leq s$. Then
$\Phi_{s,\tau}$ is a homeomorphism of $\overline{\Sigma}_{s}$ onto
$\overline{\Sigma}_{s,\tau}$ and a diffeomorphism of $\Sigma_{s}$ onto
$\Sigma_{s,\tau}.$ Furthermore, the Brown measure $\mu_{s,\tau}$ is the
push-forward of the Brown measure $\mu_{s,s}$ under the map $\Phi_{s,\tau}.$
\end{theorem}

In light of Definition \ref{PhiMap.def}, the theorem can be interpreted as
follows: as $\tau$ varies, the Brown measure $\mu_{s,\tau}$ pushes forward
along the characteristic curves of the PDE for its log potential.

Note that we have already established (Theorem \ref{outside.thm1}) that the
Brown measure $\mu_{s,s}$ is zero outside $\overline{\Sigma}_{s}%
=\overline{\Sigma}_{s,s}.$ It is therefore meaningful to speak of pushing
forward $\mu_{s,s}$ by a map that is only defined on $\overline{\Sigma}_{s}.$
The homeomorphism and diffeomorphism properties of $\Phi_{s,\tau}$ are
established in Proposition \ref{domainMap.prop} in this section. The proof of
the push-forward property is given in Section \ref{pushProof.sec}.

In \cite[Proposition 2.5]{DHKBrown} and \cite[Corollary 4.30]{HZ}, a map
$\Phi_{s}$ from $\Sigma_{s}$ into the unit circle is constructed with the
property that the push-forward of the Brown measure $\mu_{s,s}$ is the measure
$\mu_{s},$ the law of $uu_{s},$ where $u_{s}$ is the free unitary Brownian
motion and $u$ is the unitary element in the expression $ub_{s,\tau}$, taken
to be freely independent of $u_{s}.$ Now, $uu_{s}$ is just the limit in $\ast
$-distribution of $ub_{s,\tau}$ as $\tau$ tends to zero (see (\ref{AtTauZero})
and Appendix \ref{continuousDependence.sec}), and we will show in Proposition
\ref{domainMap.prop} that $\Phi_{s}$ is the the limit of $\Phi_{s,\tau}$ as
$\tau$ tends to zero. Theorem \ref{pushGeneral.thm} thus shows that the
relationship between $\mu_{s,s}$ and $\mu_{s}$ in \cite{DHKBrown,HZ} is a
limiting case of a whole family of relationships among the measures
$\mu_{s,\tau}$ with $s$ fixed and $\tau$ varying.

Of course, it follows from Theorem \ref{pushGeneral.thm} that the push-forward
of $\mu_{s,\tau_{1}}$ under the map $\Phi_{s,\tau_{2}}\circ\Phi_{s,\tau_{1}%
}^{-1}$ is $\mu_{s,\tau_{2}},$ for any two nonzero complex numbers $\tau_{1}$
and $\tau_{2}$ with $\left\vert \tau_{1}-s\right\vert \leq s$ and $\left\vert
\tau_{2}-s\right\vert \leq s.$

We now compute the map $\Phi_{s,\tau}$ more explicitly.

\begin{proposition}
\label{domainMap.prop}Fix a positive number $s$ and a nonzero complex number
$\tau$ satisfying $\left\vert \tau-s\right\vert \leq s.$ Then for all
$\lambda=re^{i\theta}$ in $\overline{\Sigma}_{s,}$ we may compute
$\Phi_{s,\tau}(\lambda)$ as%
\begin{equation}
\Phi_{s,\tau}(re^{i\theta})=\left(  \frac{r}{r_{s}(\theta)}\right)  ^{\tau
/s}f_{s-\tau}(r_{s}(\theta)e^{i\theta}). \label{PhiStau}%
\end{equation}
In particular, $\Phi_{s,\tau}$ agrees with $f_{s-\tau}$ on the boundary of
$\Sigma_{s}$ and maps each radial segment to a portion of an exponential
spiral (a curve of the form $u\mapsto ce^{u\tau}$ for some $c\in\mathbb{C}$)
in $\Sigma_{s,\tau}.$ From these formulas, we can verify that $\Phi_{s,\tau}$
is a homeomorphism of $\overline{\Sigma}_{s}$ onto $\overline{\Sigma}_{s,\tau
}$ and a diffeomorphism of $\Sigma_{s}$ onto $\Sigma_{s,\tau}.$ Furthermore,
we have%
\begin{equation}
\lim_{\tau\rightarrow0}\Phi_{s,\tau}(\lambda)=\Phi_{s}(\lambda),
\label{PhiLim}%
\end{equation}
where, as in \cite{DHKBrown,HZ}, $\Phi_{s}:\overline{\Sigma}_{s}\rightarrow
S^{1}$ is given by%
\begin{equation}
\Phi_{s}(re^{i\theta})=f_{s}(r_{s}(\theta)e^{i\theta}). \label{phiSdef}%
\end{equation}

\end{proposition}

See Figure \ref{pushmap0.fig}. We note that $\Phi_{s,\tau}(\lambda)$ depends
holomorphically on $\tau$ with $\lambda$ and $s$ fixed; this property of
$\Phi_{s,\tau}$ will be used in Section \ref{pushFromTauEqualsS.sec}. On the
other hand, $\Phi_{s,\tau}(\lambda)$ \textit{does not} depend holomorphically
on $\lambda$ with $s$ and $\tau$ fixed.

\begin{proof}
Suppose we take $\varepsilon_{0}p_{\varepsilon,0}=0$, fix a point $\lambda
_{0}=e^{i\phi}$ in the support of $\mu_{s}$ in the unit circle, and take
$p_{\theta,0}$ as in Proposition \ref{momentae0.prop}. Then the proof of
Theorem \ref{surjectivity.thm} shows that the point $\lambda(s)$ will have the
form $\lambda_{0}e^{i\alpha}r,$ where $r>0$ and the value of $\alpha$ is
independent of the limiting value of $p_{\rho,0}.$ As $p_{\rho,0}$ varies
between $\partial S/\partial\rho^{\operatorname{in}}$ and $\partial
S/\partial\rho^{\operatorname{out}},$ the value of $r$ will vary, so that
$\lambda(s)$ will range over a radial segment inside $\Sigma_{s}$ connecting
the inner and outer boundaries. Then by (\ref{lambdaOfTau}), $\lambda(\tau)$
will have the form $\lambda_{0}e^{i\tau\alpha/s}r^{\tau/s}.$ Thus,
$\Phi_{s,\tau}$ will map the radial segment in $\Sigma_{s}$ in a log-linear
fashion to a spiral segment in $\Sigma_{s,\tau}.$

We then check that the formula in (\ref{PhiStau}) captures this behavior. We
can verify that the right-hand side of (\ref{PhiStau}) agrees with $f_{s-\tau
}(re^{i\theta})$ both on the inner boundary at $r=r_{s}(\theta)$ (trivially)
and on the outer boundary at $r=r_{s}(\theta)^{-1}$ (using Lemma
\ref{lem:u1u2}). It follows that the endpoints of the radial segment map to
the endpoints of the spiral segment. As $r$ varies, the right-hand side of
(\ref{PhiStau}) traces out an exponential spiral of the form $cr^{\tau/s}$,
along which $\delta$ is constant. Thus, (\ref{PhiStau}) agrees with the map in
the previous paragraph.

That $\Phi_{s,\tau}$ is a homeomorphism of $\overline{\Sigma}_{s}$ onto
$\overline{\Sigma}_{s,\tau}$ follows from the bijective nature of the
correspondence between $\theta$ and $\delta$ (Point \ref{DomainDelta.point} of
Theorem \ref{thm:Sigmastau}), together with Lemma \ref{lem:u1u2}, which tells
us that a radial segment in $\Sigma_{s}$ cannot collapse to a single point in
$\Sigma_{s,\tau}.$ That $\Phi_{s,\tau}$ is a diffeomorphism on $\Sigma_{s}$
follows by using polar coordinates $(r,\theta)$ in $\Sigma_{s}$ and twisted
logarithmic coordinates $(v,\delta)$ in $\Sigma_{s,\tau}.$ Then $\delta$
depends only on $\theta$ and, by Proposition \ref{prop:dpsidtheta},
$d\delta/d\theta>0$. Furthermore, $\partial v/\partial r>0$ from the explicit
formula in (\ref{PhiStau}). Thus, the Jacobian of $\Phi_{s,\tau}$ in these
coordinates is upper triangular with positive entries on the diagonal.

Finally, the expression on the right-hand side of (\ref{PhiStau}) is
continuous up to $\tau=0$ and
\[
\lim_{\tau\rightarrow0}\Phi_{s,\tau}(re^{i\theta})=f_{s}(r_{s}(\theta
)e^{i\theta}),
\]
as claimed in (\ref{PhiLim}) and (\ref{phiSdef}).
\end{proof}

We now give formulas for the composite map $\Phi_{s,\tau^{\prime}}\circ
\Phi_{s,\tau}^{-1}$ and the $\bar{\lambda}$-derivative of its logarithm.
Equation (\ref{eq:logPhi}) shows that $\Phi_{s,\tau^{\prime}}\circ\Phi
_{s,\tau}^{-1}$ can be computed from the Brown measure $\mu_{s,\tau}$ (and its
log potential). Equation (\ref{eq:DLogPhi}), meanwhile, shows precisely how
$\Phi_{s,\tau^{\prime}}\circ\Phi_{s,\tau}^{-1}$ (or, equivalently, its
logarithm) fails to be conformal on $\Sigma_{s,\tau}.$

\begin{proposition}
\label{PhiComposite.prop}Let $\tau$ and $\tau^{\prime}$ be two nonzero complex
numbers such that $\left\vert s-\tau\right\vert \leq s$ and $\left\vert
s-\tau^{\prime}\right\vert \leq s$. Then the map $\Phi_{s,\tau^{\prime}}%
\circ\Phi_{s,\tau}^{-1}$ can be computed as follows:
\begin{equation}
\Phi_{s,\tau^{\prime}}(\Phi_{s,\tau}^{-1}(\lambda))=\lambda\exp\left\{
\frac{(\tau^{\prime}-\tau)}{2}\left(  2\lambda\frac{\partial S_{0}}%
{\partial\lambda}(s,\tau,\lambda)-1\right)  \right\}  . \label{eq:logPhi}%
\end{equation}
Therefore, the $\bar{\lambda}$-derivative of $\log\Phi_{s,\tau^{\prime}}%
(\Phi_{s,\tau}^{-1}(\lambda))$ has a direct relation to the density
$W(s,\tau,\lambda)$ of $\mu_{s,\tau},$ as follows:
\begin{equation}
\bar{\lambda}\frac{\partial}{\partial\bar{\lambda}}\log\Phi_{s,\tau^{\prime}%
}(\Phi_{s,\tau}^{-1}(\lambda))=\pi(\tau^{\prime}-\tau)\left\vert
\lambda\right\vert ^{2}W(s,\tau,\lambda),\quad\lambda\in\Sigma_{s,\tau}.
\label{eq:DLogPhi}%
\end{equation}

\end{proposition}

\begin{proof}
Fix $\lambda\in\overline{\Sigma}_{s,\tau}$ and choose initial conditions
$\lambda_{0}$ and---through the limiting process in Definition
\ref{PhiMap.def}---$p_{\lambda,0}$ such that, with $\varepsilon_{0}%
p_{\varepsilon_{0}}=0$, we have
\[
\lambda(\tau)=\lambda.
\]
Then by Definition \ref{PhiMap.def}, we will have $\Phi_{s,\tau}^{-1}%
(\lambda)=\lambda(s),$ where $\lambda(s)$ is computed with the same initial
conditions as $\lambda(\tau).$ Then since $\lambda(s)=\Phi_{s,\tau}%
^{-1}(\lambda)$, we may continue to use the same initial conditions and
compute
\[
\Phi_{s,\tau^{\prime}}(\Phi_{s,\tau}^{-1}(\lambda))=\lambda(\tau^{\prime}).
\]

By Corollary \ref{S0pde.cor}, $\lambda(\tau)$ and $\lambda(\tau^{\prime})$ can
be computed using (\ref{lambdaOfTau}) with $\varepsilon_{0}p_{\varepsilon
,0}=0$. Using initial conditions $\lambda_{0}$ and $p_{\lambda,0}$ as in the
previous paragraph, we then have
\[
\lambda=\lambda(\tau)=\lambda_{0}\exp\left(  \frac{\tau}{2}(2\lambda
_{0}p_{\lambda,0}-1)\right)
\]
and
\[
\Phi_{s,\tau^{\prime}}(\Phi_{s,\tau}^{-1}(\lambda))=\lambda(\tau^{\prime
})=\lambda_{0}\exp\left(  \frac{\tau^{\prime}}{2}(2\lambda_{0}p_{\lambda
,0}-1)\right)  .
\]
Thus,
\begin{align*}
\Phi_{s,\tau^{\prime}}(\Phi_{s,\tau}(\lambda))  &  =\lambda\exp\left(
\frac{(\tau^{\prime}-\tau)}{2}(2\lambda_{0}p_{\lambda,0}-1)\right) \\
&  =\lambda\exp\left(  \frac{(\tau^{\prime}-\tau)}{2}(2\lambda(\tau
)p_{\lambda}(\tau)-1)\right)  ,
\end{align*}
where, in the second equality, we have used that $\lambda(\tau)p_{\lambda
}(\tau)$ is independent of $\tau$.

Since $\lambda(\tau)=\lambda$, we may apply the second Hamilton--Jacobi
formula (Corollary \ref{S0pde.cor}) to obtain (\ref{eq:logPhi}). Finally,
(\ref{eq:DLogPhi}) follows by applying $\partial/\partial\bar{\lambda}$ to the
logarithm of (\ref{eq:logPhi}), recalling that $W=\frac{1}{4\pi}%
\Delta_{\lambda}S_{0}.$
\end{proof}

We now begin working toward the proof of Theorem \ref{pushGeneral.thm}. In
Section \ref{derivDensity.sec}, we will compute how the Brown measure
$\mu_{s,\tau}$ varies as $\tau$ varies. Then in Section
\ref{pushFromTauEqualsS.sec}, we will compute how the push-forward of
$\mu_{s,s}$ by $\Phi_{s,\tau}$ varies as $\tau$ varies. The proof (Section
\ref{pushProof.sec}) will be based on a similarity between these two
computations. The proof does not rely on the \textit{formula} for the Brown
measure obtained in Section \ref{inside.sec} but instead uses only that the
log potential of the Brown measure satisfies the $\varepsilon=0$ PDE in
Corollary \ref{S0pde.cor}.

\subsection{The derivative of the density of the Brown
measure\label{derivDensity.sec}}

We let $W(s,\tau,\lambda)$ denote the density (in $\Sigma_{s,\tau}$) of the
Brown measure $\mu_{s,\tau}.$ We begin by computing how $W(s,\tau,\lambda)$
varies as $\tau$ varies with $s$ fixed. Recall that $S_{0}(s,\tau,\lambda)$ is
the limit of $S(s,\tau,\lambda,\varepsilon)$ as $\varepsilon$ tends to zero
from above and that the density \thinspace$W$ of the Brown measure is
$\frac{1}{4\pi}$ times the Laplacian of $S_{0}$ with respect to $\lambda.$

\begin{theorem}
\label{dWdTau.thm}For all $s>0,$ the function $W$ satisfies the equation%
\begin{equation}
\frac{\partial W}{\partial\tau}=-\frac{\partial}{\partial\lambda}\left[
\lambda\left(  \lambda\frac{\partial S_{0}}{\partial\lambda}-\frac{1}%
{2}\right)  W\right]  \label{wPDE}%
\end{equation}
for all $\tau$ with $\left\vert \tau-s\right\vert <s$ and all $\lambda$ in
$\Sigma_{s,\tau}.$
\end{theorem}

Note that (\ref{wPDE}) is not a self-contained PDE\ for $W,$ since the
right-hand side involves $S_{0}$ as well as $W.$

\begin{proof}
According to Corollary \ref{S0pde.cor}, $S_{0}$ satisfies the $\varepsilon=0$
version of the PDE (\ref{thePDE}) inside $\Sigma_{s,\tau}$, as given in
(\ref{S0PDE}). We now apply the operator%
\[
\frac{1}{4\pi}\Delta_{\lambda}=\frac{1}{\pi}\frac{\partial^{2}}{\partial
\lambda\partial\bar{\lambda}}%
\]
to both sides of (\ref{S0PDE}), leaving the $\lambda$-derivative on the
right-hand side unevaluated. This gives%
\begin{equation}
\frac{1}{\pi}\frac{\partial^{2}}{\partial\lambda\partial\bar{\lambda}}%
\frac{\partial S_{0}}{\partial\tau}=-\frac{1}{2}\frac{1}{\pi}\frac{\partial
}{\partial\lambda}\left(  2\lambda^{2}\frac{\partial S_{0}}{\partial\lambda
}\frac{\partial^{2}S_{0}}{\partial\lambda\partial\bar{\lambda}}-\lambda
\frac{\partial^{2}S_{0}}{\partial\lambda\partial\bar{\lambda}}\right)  .
\label{lapS0}%
\end{equation}
After interchanging the $\tau$-derivative with the $\lambda$- and
$\bar{\lambda}$-derivatives on the left-hand side of (\ref{lapS0}), we get the
claimed result.
\end{proof}

\subsection{Pushing forward the Brown measure from $\tau=s$%
\label{pushFromTauEqualsS.sec}}

It is shown in Lemma 4.25 and Theorem 4.28 of \cite{HZ} that the Brown measure
$\mu_{s,s}$ has full mass on the open set $\Sigma_{s}.$ (This statement is
stronger than what we have proved in Theorem \ref{outside.thm1}, which tells
us only that $\mu_{s,s}$ has full mass on $\overline{\Sigma}_{s}.$) In this
section, we start with the Brown measure $\mu_{s,s}$ on $\Sigma_{s}$ and push
it forward to each domain of the form $\Sigma_{s,\tau}$, with $\tau\neq0$ and
$\left\vert \tau-s\right\vert \leq s.$ We then compute how the density of the
pushed-forward measure varies with $\tau.$ Thus, we define%
\[
\hat{\mu}_{s,\tau}=\text{ push-forward of }\mu_{s,s}\text{ by }\Phi_{s,\tau}%
\]
and we let $\hat{W}(s,\tau,\lambda)$ denote the density (in $\Sigma_{s,\tau}$)
of $\hat{\mu}_{s,\tau}.$

\begin{theorem}
\label{WhatEq.thm}For each fixed $s>0,$ the function $\hat{W}$ satisfies the
same equation as the function $W$ in Theorem \ref{dWdTau.thm}, namely
\[
\frac{\partial\hat{W}}{\partial\tau}=-\frac{\partial}{\partial\lambda}\left[
\lambda\left(  \lambda\frac{\partial S_{0}}{\partial\lambda}-\frac{1}%
{2}\right)  \hat{W}\right]  ,
\]
for all $\tau$ with $\left\vert \tau-s\right\vert <s.$
\end{theorem}

We begin by computing $\partial\hat{W}/\partial\tau$ in a general way. The
following result does not rely on the specific form of the map $\Phi_{s,\tau}$
but only on one crucial property of it---that $\Phi_{s,\tau}(\lambda)$ depends
holomorphically on $\tau$ with $\lambda$ fixed.

\begin{proposition}
\label{whatGen.prop}For each fixed $s>0,$ and all $\tau$ with $\left\vert
\tau-s\right\vert <s,$ we have
\[
\frac{\partial\hat{W}}{\partial\tau}=-\frac{\partial}{\partial\lambda}\left[
Q\hat{W}\right]  ,
\]
where%
\begin{equation}
Q(s,\tau,\lambda)=\left.  \frac{\partial}{\partial\tau^{\prime}}\Phi
_{s,\tau^{\prime}}(\Phi_{s,\tau}^{-1}(\lambda))\right\vert _{\tau^{\prime
}=\tau}. \label{qq}%
\end{equation}

\end{proposition}

Before proving the proposition, we establish a useful lemma.

\begin{lemma}
\label{Bruce.lem}Suppose $f_{r}$ is a family of diffeomorphisms such that
$f_{r}$ and $f_{r}^{-1}$ depend smoothly on $r$, $\alpha$ is a fixed
differential form, and we define $\beta_{r}$ as the pull-back of $\alpha$ by
$f_{r}^{-1}.$ Then%
\begin{equation}
\frac{d\beta_{r}}{dr}=-\mathcal{L}_{X_{r}}\beta_{r}, \label{LieDeriv}%
\end{equation}
where $X_{r}$ is the vector field given by
\[
X_{r}(m)=\left.  \frac{d}{dr^{\prime}}f_{r^{\prime}}(f_{r}^{-1}(m))\right\vert
_{r^{\prime}=r}.
\]

\end{lemma}

\begin{proof}
We write $\alpha$ locally as $\alpha=g~dx_{1}\wedge\cdots\wedge dx_{n},$ so
that
\[
\beta_{r}=(g\circ f_{r}^{-1})d[x_{1}\circ f_{r}^{-1}]\wedge\cdots\wedge
d[x_{n}\circ f_{r}^{-1}],
\]
from which an elementary calculation shows%
\begin{align*}
\frac{d\beta_{r}}{dr}  &  =Y_{r}(g\circ f_{r}^{-1})d[x_{1}\circ f_{r}%
^{-1}]\wedge\cdots\wedge d[x_{n}\circ f_{r}^{-1}]\\
&  +\sum_{j}(g\circ f_{r}^{-1})d[x_{1}\circ f_{r}^{-1}]\wedge\cdots\wedge
d[Y_{r}(x_{j}\circ f_{r}^{-1})]\wedge\cdots\wedge d[x_{n}\circ f_{r}^{-1}]\\
&  =\mathcal{L}_{Y_{r}}\beta_{r},
\end{align*}
where%
\[
Y_{r}(m)=(f_{r})_{\ast}\left(  \frac{d}{dr}f_{r}^{-1}(m)\right)  .
\]
We then differentiate the identity $f_{r}(f_{r}^{-1}(m))=m$ with respect to
$r$ to show that $Y_{r}=-X_{r}.$
\end{proof}

\begin{proof}
[Proof of Proposition \ref{whatGen.prop}]It will be convenient to use the
symbol \textquotedblleft$\operatorname{vect}$\textquotedblright\ to convert a
complex number to a vector in the plane; that is,
\begin{equation}
\operatorname{vect}(\alpha+i\beta)=(\alpha,\beta). \label{vectDef}%
\end{equation}
Consider a curve of the form%
\begin{equation}
\tau(r)=\tau+r\tau^{\prime} \label{TauR}%
\end{equation}
and define a vector field $V_{r}$ by
\begin{equation}
V_{r}(\lambda)=\operatorname{vect}\left(  \left.  \frac{d}{dr^{\prime}}%
\Phi_{s,\tau(r^{\prime})}(\Phi_{s,\tau(r)}^{-1}(\lambda))\right\vert
_{r^{\prime}=r}\right)  . \label{VrFirst}%
\end{equation}
Since $\Phi_{s,\tau}(\lambda)$ depends holomorphically on $\tau$ with
$\lambda$ fixed, we may compute
\begin{align}
\left.  \frac{d}{dr^{\prime}}\Phi_{s,\tau(r^{\prime})}(\Phi_{s,\tau(r)}%
^{-1}(\lambda))\right\vert _{r^{\prime}=r}  &  =\frac{\partial\Phi_{s,\tau}%
}{\partial\tau}(\Phi_{s,\tau(r)}^{-1}(\lambda))\frac{d\tau}{dr}\nonumber\\
&  =\tau^{\prime}\left.  \frac{\partial\Phi_{s,\tau}}{\partial\tau}%
(\Phi_{s,\tau(r)}^{-1}(\lambda))\right\vert _{\tau=\tau(r)}. \label{phiDeriv}%
\end{align}

Now, since $\hat{W}(s,\tau(r),\cdot)$ is the density of the push-forward of
$\mu_{s,s}$ by $\Phi_{s,\tau(r)},$ the associated differential form $\hat
{W}(s,\tau(r),\lambda)~dx\wedge dy$ is the pull-back of $W(s,s,\lambda
)~dx\wedge dy$ by the inverse map $\Phi_{s,\tau(r)}^{-1}.$ By Lemma
\ref{Bruce.lem}, we obtain%
\begin{equation}
\frac{d}{dr}\hat{W}(s,\tau(r),\lambda)~dx\wedge dy=-\mathcal{L}_{V_{r}}%
[\hat{W}(s,\tau(r),\lambda)~dx\wedge dy]. \label{ref1}%
\end{equation}
We now compute $\mathcal{L}_{V_{r}}$ using Cartan's formula, one term of which
is zero because $\hat{W}~dx\wedge dy$ is a top-degree form. Now, using
(\ref{phiDeriv}), we find that the vector field $V_{r}$ in (\ref{VrFirst}) is
equal to $\operatorname{vect}(\tau^{\prime}Q(s,\tau(r),\lambda)),$ where $Q$
is as in (\ref{qq}). We may then easily compute that
\begin{equation}
\mathcal{L}_{V_{r}}[\hat{W}~dx\wedge dy]=\left(  \tau^{\prime}\frac{\partial
}{\partial\lambda}(Q\hat{W})+\overline{\tau^{\prime}}\frac{\partial}%
{\partial\bar{\lambda}}(\bar{Q}\hat{W})\right)  ~dx\wedge dy. \label{ref2}%
\end{equation}

We now argue as in the proof of Theorem \ref{thm:dSdtau}. Recalling that
$\tau(r)=\tau+r\tau^{\prime},$ we find that%
\[
\left.  \frac{\partial}{\partial r}\hat{W}(s,\tau(r),\lambda)\right\vert
_{r=0}=\frac{\partial\hat{W}}{\partial\tau}(s,\tau,\lambda)\tau^{\prime}%
+\frac{\partial\hat{W}}{\partial\bar{\tau}}(s,\tau,\lambda)\overline
{\tau^{\prime}}.
\]
Thus, using (\ref{ref1}) and (\ref{ref2}), we obtain
\begin{align*}
\frac{\partial\hat{W}}{\partial\tau}(s,\tau,\lambda)  &  =\frac{\partial
}{\partial\tau^{\prime}}\left(  \left.  \frac{\partial}{\partial r}\hat
{W}(s,\tau(r),\lambda)\right\vert _{r=0}\right) \\
&  =-\left.  \frac{\partial}{\partial\lambda}(Q(s,\tau(r),\lambda)\hat
{W})\right\vert _{r=0},
\end{align*}
as claimed.
\end{proof}

To prove Theorem \ref{WhatEq.thm}, we simply need to compute the function $Q$
in Proposition \ref{whatGen.prop}.

\begin{proof}
[Proof of Theorem \ref{WhatEq.thm}]To compute $\Phi_{s,\tau}^{-1}(\lambda),$
using Definition \ref{PhiMap.def}, we choose our initial conditions so that
$\lambda(\tau)=\lambda,$ and then $\Phi_{s,\tau}^{-1}(\lambda)$ is equal to
$\lambda(s).$ Then to compute $\Phi_{s,\tau^{\prime}}(\Phi_{s,\tau}%
^{-1}(\lambda)),$ we continue to use the same initial conditions---since they
give $\lambda(s)=\Phi_{s,\tau}^{-1}(\lambda)$---and then $\Phi_{s,\tau
^{\prime}}(\Phi_{s,\tau}^{-1}(\lambda))$ is equal to $\lambda(\tau^{\prime}).$
Thus, the values $\Phi_{s,\tau^{\prime}}(\Phi_{s,\tau}^{-1}(\lambda))$ have
the form $\lambda(\tau^{\prime}),$ where all curves are computed using a fixed
set of initial conditions. We then differentiate the formula
(\ref{lambdaOfTau}) for $\lambda(\tau)$ with $\varepsilon_{0}p_{\varepsilon
,0}=0,$ and find that
\[
\frac{d\lambda(\tau)}{d\tau}=\lambda(\tau)\left(  \lambda_{0}p_{\lambda
,0}-\frac{1}{2}\right)  .
\]
Since $\lambda(\tau)p_{\lambda}(\tau)$ is independent of $\tau$ (from the
formulas (\ref{lambdaOfTau}) and (\ref{epsilonOfTau})) and using the second
Hamilton--Jacobi formula (\ref{S0HJsecond}) for $S_{0}$, we get%
\begin{equation}
\frac{d\lambda(\tau)}{d\tau}=\lambda(\tau)\left(  \lambda(\tau)\frac{\partial
S_{0}}{\partial\lambda}(s,\tau,\lambda(\tau))-\frac{1}{2}\right)  .
\label{derivZtau}%
\end{equation}
Thus,
\[
Q(s,\tau,\lambda)=\lambda\left(  \lambda\frac{\partial S_{0}}{\partial\lambda
}(s,\tau,\lambda)-\frac{1}{2}\right)
\]
and Theorem \ref{WhatEq.thm} follows from Proposition \ref{whatGen.prop}.
\end{proof}

\subsection{Proof of the push-forward property\label{pushProof.sec}}

We now prove our main push-forward result, Theorem \ref{pushGeneral.thm}. The
proof is based on the fact that $W$ (the density of the measure $\mu_{s,\tau}%
$) and $\hat{W}$ (the density of the push-forward of $\mu_{s,s}$ by
$\Phi_{s,\tau}$) satisfy the same PDE with the same initial condition at
$r=0.$ (In both equations, recorded in Theorems \ref{dWdTau.thm} and
\ref{WhatEq.thm}, we regard $S_{0}$ as a known quantity and try to solve for
$W$ or $\hat{W}.$) The agreement between the PDEs, in turn, relies on an
agreement between the expression for $d\lambda/d\tau$ in (\ref{derivZtau}) and
the coefficient of $W$ inside the square brackets on the right-hand side of
(\ref{wPDE}).

Now, the quantity in square brackets in (\ref{wPDE}), comes from applying
$\frac{1}{\pi}\partial/\partial\bar{\lambda}$ to the right-hand side of the
PDE (\ref{S0PDE}), while the formula for $\partial\lambda/\partial\tau$ comes
from the Hamilton--Jacobi analysis of the same PDE. While we will refrain from
stating a precise result, one could develop a general theory whenever a
real-valued function $S_{0}$ satisfies a PDE of the form%
\[
\frac{\partial S_{0}}{\partial\tau}=F\left(  \lambda,\frac{\partial S_{0}%
}{\partial\lambda}\right)  ,
\]
where $F$ is a holomorphic of two complex variables.

\begin{proof}
[Proof of Theorem \ref{pushGeneral.thm}]We will prove that the function $W$ in
Section \ref{derivDensity.sec} and the function $\hat{W}$ in Section
\ref{pushFromTauEqualsS.sec} are equal, showing that the push-forward of
$\mu_{s,s}$ by $\Phi_{s,\tau}$ agrees with $\mu_{s,\tau(r)}$ on the open set
$\Sigma_{s,\tau(r)}$. Then since the set $\Sigma_{s}$ has measure 1 with
respect to $\mu_{s,s}$ \cite[Lemma 4.25 and Theorem 4.28]{HZ}, the set
$\Sigma_{s,\tau(r)}$ must have measure 1 with respect to $\mu_{s,\tau(r)}.$
Thus, the push-forward of $\mu_{s,s}$ by $\Phi_{s,\tau}$ is actually equal to
$\mu_{s,\tau}.$

We treat the function $S_{0}$ as a known quantity, and note that this function
occurs both in the PDE for $W$ (Theorem \ref{dWdTau.thm}) and in the identical
PDE for $\hat{W}$ (Theorem \ref{WhatEq.thm}). We can then solve this PDE for
$W$ (or $\hat{W}$) by differentiating along the characteristic curves
$\lambda(\tau).$ We compute that%
\[
\frac{\partial}{\partial\tau}W(s,\tau,\lambda(\tau))=\frac{\partial
W}{\partial\tau}(s,\tau,\lambda(\tau))+\frac{\partial W}{\partial\lambda
}(s,\tau,\lambda(\tau))\frac{\partial\lambda}{\partial\tau}.
\]
Then using the PDE (\ref{wPDE}) and the formula (\ref{derivZtau}) for
$\partial\lambda/\partial\tau,$ we obtain, after noting a cancellation and
dividing both sides by $W$:%
\[
\frac{\partial}{\partial\tau}\log W(s,\tau,\lambda(\tau))=-\frac{\partial
}{\partial\lambda}\left[  \lambda\left(  \lambda\frac{\partial S_{0}}%
{\partial\lambda}-\frac{1}{2}\right)  \right]  .
\]

We now consider a nonzero $\tau$ with $\left\vert \tau-s\right\vert \leq s$
and consider the curve $\tau(r)=s+r(\tau-s),$ which is simply the curve in
(\ref{TauR}) with $\tau=s$ and $\tau^{\prime}=(\tau-s).$ We then wish to show
that
\begin{equation}
W(s,\tau(r),\lambda)=\hat{W}(s,\tau(r),\lambda) \label{WEqualsW}%
\end{equation}
for all $0\leq r\leq1,$ where equality at $r=1$ shows that $W=\hat{W}.$ We
then compute that for $0\leq r<1,$ we have%
\begin{align}
\frac{d}{dr}\log W(s,\tau(r),\lambda(\tau(r)))  &  =\frac{\partial}%
{\partial\tau}\log W(s,\tau(r),\lambda(\tau(r)))\frac{d\tau}{dr}\nonumber\\
&  +\frac{\partial}{\partial\bar{\tau}}\log W(s,\tau(r),\lambda(\tau
(r)))\frac{d\bar{\tau}}{dr}\nonumber\\
&  =-2\operatorname{Re}\left[  (\tau-s)H_{r}(\lambda(\tau(r)))\right]
\label{dLogW}%
\end{align}
where%
\[
H_{r}(\lambda)=-\frac{\partial}{\partial\lambda}\left[  \lambda\left(
\lambda\frac{\partial S_{0}}{\partial\lambda}(\tau(r),\lambda)-\frac{1}%
{2}\right)  \right]  .
\]
Of course, $\hat{W}$ satisfies the same equation with $W$ changed to $\hat{W}$
everywhere. We impose the restriction $r<1$ because in the borderline case
$\left\vert \tau-s\right\vert =s,$ the PDEs for $W$ and $\hat{W}$ are not
known to hold at $\tau.$

At $r=0,$ both functions equal the density of the Brown measure $\mu_{s,s}.$
We may then integrate in (\ref{dLogW}) to obtain%
\begin{align*}
\log W(s,\tau,\lambda(\tau(r)))  &  =\log\hat{W}(s,\tau,\lambda(\tau(r)))\\
&  =\log W(s,s,0)-2\int_{0}^{r}\operatorname{Re}\left[  (\tau-s)H_{r^{\prime}%
}(\lambda(\tau(r^{\prime})))\right]  ~dr^{\prime},
\end{align*}
for $r<1.$ We now wish to let $r$ tend to 1. This is only a problem if
$\left\vert \tau-s\right\vert =s.$ (If $\left\vert \tau-s\right\vert <s,$ the
PDEs for $\hat{W}$ and $W$ are applicable even a little beyond $r=1.$) As long
as $\tau\neq0,$ all the estimates in Section \ref{sect:fInjBdry} concerning
the boundary behavior of $f_{s-\tau}$ remain applicable even when $\left\vert
\tau-s\right\vert =s,$ from which good behavior of the map $\Phi_{s,\tau}$
follows. Thus, we can see that $\hat{W}(s,\tau,\lambda(\tau(r)))$ is
continuous up to $r=1.$ The analysis in Section \ref{Sreg.sec} of $S$ at
$\varepsilon=0$ is similarly applicable for $\tau\neq0$ and $\left\vert
\tau-s\right\vert =s,$ from which we can see that $W(s,\tau,\lambda(\tau(r)))$
is continuous up to $r=1.$ Thus, (\ref{WEqualsW}) holds at $r=1,$ showing that
$W=\hat{W}.$
\end{proof}

\subsection{The proof of Theorem \ref{fullMass.thm}\label{fullMass.sec}}

The theorem in question says that the Brown measure $\mu_{s,\tau}$ gives full
mass to the open set $\Sigma_{s,\tau}.$ This is a strengthening of Theorem
\ref{outside.thm1}, which says that $\mu_{s,\tau}$ gives full mass to the
closed set $\overline{\Sigma}_{s,\tau}.$ As discussed after the statement of
Theorem \ref{fullMass.thm} in Section \ref{theBrownMeasure.sec}, if the
boundary of $\Sigma_{s,\tau}$ is a rectifiable curve, Theorem
\ref{fullMass.thm} follows from the agreement between the first derivatives of
$S_{0}$ as we approach the boundary from inside and from outside the
domain---but we do not know that rectifiability of the boundary holds in general.

We actually have already proved Theorem \ref{fullMass.thm} in general as part
of the proof of Theorem \ref{pushGeneral.thm} in the previous subsection. The
proof relies on a result of Ho and Zhong \cite[Lemma 4.25 and Theorem
4.28]{HZ}, that the Brown measure of $\mu_{s,s}$ has full mass on $\Sigma
_{s}.$ Then the push-forward result on the open domains then tells us that%
\[
\mu_{s,\tau}(\Sigma_{s,\tau})=\mu_{s,s}(\Phi_{s,\tau}^{-1}(\Sigma_{s,\tau
}))=\mu_{s,s}(\Sigma_{s})=1.
\]

Now, the proof of the just-cited result of Ho and Zhong also relies on a
push-forward result. Recall that $\mu_{s}$ denotes the law of the unitary
element $uu_{s},$ which is known \cite[Proposition 3.6]{Zhong} to have a
continuous density with respect to Lebesgue measure on the unit circle. The
push-forward result of Ho and Zhong (implicit in their Lemma 4.25 and stated
explicitly in Proposition 4.31) asserts that the push-forward of the
restriction of $\mu_{s,s}$ to $\Sigma_{s}$ coincides with the restriction of
$\mu_{s}$ to a certain open set $V_{s}$ in the unit circle. By Theorems 2.11
and 4.10 in \cite{HZ}, the set $V_{s}$ is precisely the set where the density
of $\mu_{s}$ is positive, so that $\mu_{s}(V_{s})$ equals 1. It then follows
that $\mu_{s,s}(\Sigma_{s})$ equals 1 as well.

Intuitively, the push-forward result in \cite{HZ} can be understood as the
$\tau\rightarrow0$ limit of Theorem \ref{pushGeneral.thm} in the present
paper, but their result does not actually follow from our result, without
additional effort.

\appendix{}

\section{A factorization of $b_{s,\tau}$\label{factorization.app}}

In this section, we establish the factorization result in Theorem
\ref{thm:factorization} for elements of the form $b_{s,\tau}$ (Notation
\ref{bstau.notation}). This result is used in the derivation of the PDE\ for
$S$ in Section \ref{pde.sec}.

One could prove a similar result in the finite-$N$ setting of
\cite{DHKcomplex}, replacing, say, $b_{s,\tau}(r)$ by the corresponding
Brownian motion in $GL(N;\mathbb{C}).$ According to the equation between Eq.
(1.13) and Eq. (1.14) in \cite{DHKcomplex}, the generator of the $(s,\tau
)$-Brownian motion in $GL(N;\mathbb{C})$ is the operator $s\Delta-\tau
\partial^{2}-\bar{\tau}\bar{\partial}^{2},$ where $\Delta$, $\partial^{2},$
and $\bar{\partial}^{2}$ are certain left-invariant differential operators on
$GL(N;\mathbb{C}).$ In this setting, the factorization result holds because
all three of these operators commute with one another \cite[Corollary
5.7]{DHKcomplex}.

One obvious strategy for trying to prove Theorem \ref{thm:factorization} in
the free setting is to try to prove that $b_{s,\tau}(r)b_{s^{\prime}%
,\tau^{\prime}}^{\prime}(r)$ and $b_{s+s^{\prime},\tau+\tau^{\prime}}(r)$ have
the same $\ast$-distribution for all $r.$ We may then try to show that both
sets of $\ast$-moments satisfy the same differential equation with respect to
$r$, with the same trivial initial condition at $r=0.$ This approach does not,
however, work in general. Actually, we will implement precisely this strategy
in Section \ref{unitaryFact.sec} \textit{in the special case where}
$\tau^{\prime}=0.$ In the $\tau^{\prime}=0$ case, we get $b_{s^{\prime}%
,0}(r)=u_{rs^{\prime}},$ where $u_{\cdot}$ is the free unitary Brownian
motion, and the unitarity of $u_{rs^{\prime}}$ is used in an essential way in
the analysis. (See, specifically, the derivation of (\ref{2dAr3}) in the proof
of Proposition \ref{prop:bfactor}.)

To show that the $\ast$-moments of $b_{s,\tau}b_{s^{\prime},\tau^{\prime}%
}^{\prime}$ and $b_{s+s^{\prime},\tau+\tau^{\prime}}$ are the same, we follow
a three-step process. In Step 1, we will prove this result (as we have just
discussed) in the case $\tau^{\prime}=0.$ In Step 2, we will determine how the
$\ast$-moments of $b_{s,\tau}$ (multiplied on the left and right by freely
independent elements) change as we\textit{ vary }$\tau$\textit{ with }%
$s$\textit{ fixed}. To determine this, we will need to use the factorization
result in Step 1. Finally, in Step 3, we will use Step 2 to show that the
$\ast$-moments of
\[
b_{s,t\tau}b_{s^{\prime},t\tau^{\prime}}\text{ and }b_{s+s^{\prime}%
,t(\tau+\tau^{\prime})}%
\]
satisfy the same differential equation in $t,$ with the same value at $t=0.$
We will then conclude the $\ast$-moments are equal for all $t$; setting $t=1$
gives the factorization theorem.

\subsection{Factoring with a free unitary Brownian
motion\label{unitaryFact.sec}}

In this section, we carry out Step 1 in the proof of Theorem
\ref{thm:factorization}, corresponding to the special case $\tau^{\prime}=0.$
In this case, the free multiplicative Brownian motion $b_{s^{\prime}%
,0}^{\prime}(r)$ is computable as%
\[
b_{s^{\prime},0}^{\prime}(r)=u_{rs^{\prime}},
\]
where $u_{\cdot}$ is a free unitary Brownian motion freely independent of
$b_{s,\tau}(r).$ The unitarity of $u_{\cdot}$ gives us the identities
\begin{align}
b_{s,\tau}(r)b_{s,\tau}(r)^{\ast}  &  =(b_{s,\tau}(r)u_{rs^{\prime}%
})(b_{s,\tau}(r)u_{rs^{\prime}})^{\ast}\nonumber\\
u_{rs^{\prime}}^{\ast}u_{rs^{\prime}}  &  =1, \label{eq:bb*Combine}%
\end{align}
which of course do not hold if $u_{rs^{\prime}}$ is replaced by a generic
element of the form $b_{s^{\prime},\tau^{\prime}}^{\prime}(r).$ These
identities will play an essential role in the proof of the following result.

\begin{proposition}
\label{prop:bfactor} Choose $s>0,$ $s^{\prime}>0$, and $\tau\in\mathbb{C}$ so
that $\left\vert \tau-s\right\vert \leq s$. Let $u_{\cdot}$ be a free unitary
Brownian motion that is freely independent of $b_{s,\tau}(r).$ Then for any
$r>0$, the random variables $b_{s+s^{\prime},\tau}(r)$ and
\begin{equation}
a_{s,s^{\prime},\tau}(r):=b_{s,\tau}(r)u_{rs^{\prime}} \label{astauDef}%
\end{equation}
have the same $\ast$-distribution.
\end{proposition}

We will prove Proposition~\ref{prop:bfactor} using the free It\^{o} rules. If
we think of $u_{rs^{\prime}}$ as a free stochastic process with time-variable
$r,$ we obtain the scaled version of the free SDE for $u_{t}$, namely
\[
du_{rs^{\prime}}=u_{rs^{\prime}}\left(  i\sqrt{s^{\prime}}~dx_{s}%
-\frac{s^{\prime}}{2}dr\right)  .
\]
Using this result, the SDE (\ref{bstauSDE}) for $b_{s,\tau}(r),$ and the
stochastic product rule (\ref{eq:Itoprod}), we find that the process
$a_{s,s^{\prime},\tau}(r)$ satisfies the following free SDE:%
\begin{align}
da_{s,s^{\prime},\tau}(r)  &  =i\sqrt{s^{\prime}}a_{s,s^{\prime},\tau
}(r)~dx_{r}\nonumber\\
&  +ib_{s,\tau}(r)~dw_{s,\tau}(r)~u_{rs^{\prime}}-\frac{1}{2}(s+s^{\prime
}-\tau)a_{s,s^{\prime},\tau}(r)~dr, \label{eq:aSDE}%
\end{align}
where $x_{r}$ is a semicircular Brownian motion that is freely independent of
$w_{s,\tau}(r).$

\begin{proof}
[Proof of Proposition~\ref{prop:bfactor}]In this proof, we assume $s_{0},$
$s,$ and $\tau$ are fixed and use the notations%
\begin{align*}
b_{r}  &  =b_{s+s^{\prime},\tau}(r);\quad w_{r}=w_{s+s^{\prime},\tau}(r);\\
\tilde{b}_{r}  &  =b_{s,\tau}(r);\quad a_{r}=a_{s,s^{\prime},\tau}%
(r);\quad\tilde{w}_{r}=w_{s,\tau}(r).\quad.
\end{align*}
We prove the result by induction on the length of the $\ast$-moment. It is
clear that $b_{r}$ and $a_{r}$ have the same zeroth-order $\ast$-moments.
Suppose that all the $\ast$-moments of order less than $n$ of $b_{r}$ and
$a_{r}$ are equal and consider a $\ast$-moment of length $n,$ which we write
as $\operatorname{tr}\left[  b_{r}^{\varepsilon_{1}}\cdots b_{r}%
^{\varepsilon_{n}}\right]  $ or $\operatorname{tr}[a_{r}^{\varepsilon_{1}%
}\cdots a_{r}^{\varepsilon_{n}}],$ with $\varepsilon_{j}\in\{1,\ast\}.$

In the case of the $\ast$-moment of $b_{r},$ the stochastic product rule tells
us that
\[
d\operatorname{tr}\left[  b_{r}^{\varepsilon_{1}}\cdots b_{r}^{\varepsilon
_{n}}\right]  =\sum_{j}\operatorname{tr}\left[  b_{r}^{\varepsilon_{1}}\cdots
db_{r}^{\varepsilon_{j}}\cdots b_{r}^{\varepsilon_{n}}\right]  +\sum
_{j<k}\operatorname{tr}\left[  b_{r}^{\varepsilon_{1}}\cdots db_{r}%
^{\varepsilon_{j}}\cdots db_{r}^{\varepsilon_{k}}\cdots b_{r}^{\varepsilon
_{n}}\right]  ,
\]
where the sum over $j<k$ is empty (having a value of 0) if $n=1.$ Using the
SDE (\ref{bstauSDE}) of $b_{r}$ and the It\^{o} rules (\ref{eq:3Ito}), we may
easily compute that
\begin{align}
&  \frac{d}{dr}\operatorname{tr}\left[  b_{r}^{\varepsilon_{1}}\cdots
b_{r}^{\varepsilon_{n}}\right]  \nonumber\\
&  =-\frac{1}{2}\operatorname{tr}\left[  b_{r}^{\varepsilon_{1}}\cdots
b_{r}^{\varepsilon_{n}}\right]  \sum_{j}(s+s^{\prime}-\tau^{\varepsilon_{j}%
})\nonumber\\
&  -\sum_{\varepsilon_{j}=\varepsilon_{k}=1}(s+s^{\prime}-\tau
)\operatorname{tr}\left[  b_{r}^{\varepsilon_{1}}\cdots b_{r}^{\varepsilon
_{j}}b_{r}^{\varepsilon_{k+1}}\cdots b_{r}^{\varepsilon_{n}}\right]
\operatorname{tr}\left[  b_{r}^{\varepsilon_{j+1}}\cdots b_{r}^{\varepsilon
_{k}}\right]  \nonumber\\
&  \quad-\sum_{\varepsilon_{j}=\varepsilon_{k}=\ast}(s+s^{\prime}-\bar{\tau
})\operatorname{tr}\left[  b_{r}^{\varepsilon_{1}}\cdots b_{r}^{\varepsilon
_{j}}b_{r}^{\varepsilon_{k+1}}\cdots b_{r}^{\varepsilon_{n}}\right]
\operatorname{tr}\left[  b_{r}^{\varepsilon_{j+1}}\cdots b_{r}^{\varepsilon
_{k}}\right]  \nonumber\\
&  \quad+\sum_{\varepsilon_{j}=1,~\varepsilon_{k}=\ast}(s+s^{\prime
})\,\operatorname{tr}\left[  b_{r}^{\varepsilon_{1}}\cdots b_{r}%
^{\varepsilon_{j}}b_{r}^{\varepsilon_{k}}\cdots b_{r}^{\varepsilon_{n}%
}\right]  \operatorname{tr}\left[  b_{r}^{\varepsilon_{j+1}}\cdots
b_{r}^{\varepsilon_{k-1}}\right]  \nonumber\\
&  \quad+\sum_{\varepsilon_{j}=\ast,~\varepsilon_{k}=1}(s+s^{\prime
})\,\operatorname{tr}\left[  b_{r}^{\varepsilon_{1}}\cdots b_{r}%
^{\varepsilon_{j-1}}b_{r}^{\varepsilon_{k+1}}\cdots b_{r}^{\varepsilon_{n}%
}\right]  \operatorname{tr}\left[  b_{r}^{\varepsilon_{j}}\cdots
b_{r}^{\varepsilon_{k}}\right]  ,\label{momentDiff}%
\end{align}
where in the last four lines, each sum is over all pairs $(j,k)$ with $j<k$
and $\varepsilon_{j}$ and $\varepsilon_{k}$ satisfying the indicated
conditions. In the first line, we take $\tau^{\varepsilon_{j}}=\bar{\tau}$ if
$\varepsilon_{j}=\ast.$

The first term on the right-hand side of (\ref{momentDiff}) is a multiple of
the $\ast$-moment on the left-hand side of the equation. All the remaining
terms on the right-hand side of (\ref{momentDiff}) are either multiples of the
left-hand side (e.g., the term $j=k-1$ in the second sum on the right-hand
side) or products of $\ast$-moments of $b$ of lower degree. Thus,
(\ref{momentDiff}) has the form%
\begin{equation}
\frac{d}{dr}\operatorname{tr}\left[  b_{r}^{\varepsilon_{1}}\cdots
b_{r}^{\varepsilon_{n}}\right]  =C\operatorname{tr}\left[  b_{r}%
^{\varepsilon_{1}}\cdots b_{r}^{\varepsilon_{n}}\right]  +f(r),
\label{momentDiff2}%
\end{equation}
where $f(r)$ may be considered, inductively, as a known quantity. The equation
(\ref{momentDiff2}) may then be solved by the method of integrating factors.

We now show that $\operatorname{tr}\left[  a_{r}^{\varepsilon_{1}}\cdots
a_{r}^{\varepsilon_{n}}\right]  $ satisfies the same ODE in (\ref{momentDiff2}%
). By the stochastic product rule,
\begin{align}
d\operatorname{tr}\left[  a_{r}^{\varepsilon_{1}}\cdots a_{r}^{\varepsilon
_{n}}\right]    & =\sum_{j}\operatorname{tr}\left[  a_{r}^{\varepsilon_{1}%
}\cdots da_{r}^{\varepsilon_{j}}\cdots a_{r}^{\varepsilon_{n}}\right]
\nonumber\\
& +\sum_{j<k}\operatorname{tr}\left[  a_{r}^{\varepsilon_{1}}\cdots
da_{r}^{\varepsilon_{j}}\cdots da_{r}^{\varepsilon_{k}}\cdots a_{r}%
^{\varepsilon_{n}}\right]  .\label{eq:amomentIto}%
\end{align}
By the SDE (\ref{eq:aSDE}) of $a$ and the free It\^{o} rules (\ref{eq:scIto})
and (\ref{eq:3Ito}), the first term on the right-hand side of
(\ref{eq:amomentIto}) is simply
\begin{equation}
-\frac{1}{2}\operatorname{tr}\left[  a_{r}^{\varepsilon_{1}}\cdots
a_{r}^{\varepsilon_{n}}\right]  \sum_{j}(s+s^{\prime}-\tau^{\varepsilon_{j}%
}).\label{eq:aFirstSum}%
\end{equation}
This term is the same as the first term in the right-hand side of
(\ref{momentDiff}), after replacing the $b_{r}$ by $a_{r}.$

For the second sum on the right-hand side of (\ref{eq:amomentIto}), we
separate the computation into four cases of $(\varepsilon_{j},\varepsilon
_{k})\in\{1,\ast\}^{2}$. If $(\varepsilon_{j},\varepsilon_{k})=(1,1)$, then,
after dropping terms that are equal to zero, we get
\begin{align}
&  \operatorname{tr}\left[  a_{r}^{\varepsilon_{1}}\cdots da_{r}%
^{\varepsilon_{j}}\cdots da_{r}^{\varepsilon_{k}}\cdots a_{r}^{\varepsilon
_{n}}\right]  \nonumber\\
&  =-s^{\prime}\operatorname{tr}\left[  a_{r}^{\varepsilon_{1}}\cdots
a_{r}^{\varepsilon_{j}}dx_{r}\cdots a_{r}^{\varepsilon_{k}}dx_{r}\cdots
a_{r}^{\varepsilon_{n}}\right]  \nonumber\\
&  +\operatorname{tr}\left[  a_{r}^{\varepsilon_{1}}\cdots(i\tilde{b}%
_{r}\,d\tilde{w}_{r}\,u_{rs^{\prime}})\cdots(i\tilde{b}_{r}\,d\tilde{w}%
_{r}\,u_{rs^{\prime}})\cdots a^{\varepsilon_{n}}\right]  .\label{2dAr1}%
\end{align}
We then apply the It\^{o} rules (\ref{eq:scIto}) and (\ref{eq:3Ito}). The
second term on the right-hand side of (\ref{2dAr1}) is the crucial one: after
using (\ref{eq:3Ito}), the $\tilde{b}_{r}$ to the left of the first
$d\tilde{w}_{r}$ and the $u_{rs^{\prime}}$ to the right of the second
$d\tilde{w}_{r}$ combine to form $a_{r}=a_{r}^{\varepsilon_{j}}.$ Meanwhile,
after a cyclic permutation inside the trace, the $u_{rs^{\prime}}$ to the
right of first $d\tilde{w}_{r}$ and the $\tilde{b}_{r}$ to the left of second
$d\tilde{w}_{r}$ also combine to form $a_{r}=a_{r}^{\varepsilon_{k}}.$ After
combining terms, we end up with%
\begin{align*}
&  \operatorname{tr}\left[  a_{r}^{\varepsilon_{1}}\cdots da_{r}%
^{\varepsilon_{j}}\cdots da_{r}^{\varepsilon_{k}}\cdots a_{r}^{\varepsilon
_{n}}\right]  \\
&  =-(s+s^{\prime}-\tau)\operatorname{tr}\left[  a_{r}^{\varepsilon_{1}}\cdots
a_{r}^{\varepsilon_{j}}a_{r}^{\varepsilon_{k+1}}\cdots a_{r}^{\varepsilon_{n}%
}\right]  \operatorname{tr}\left[  a_{r}^{\varepsilon_{j+1}}\cdots
a_{r}^{\varepsilon_{k}}\right]  .
\end{align*}
This result matches the $\varepsilon_{j}=\varepsilon_{k}=1$ term in
(\ref{momentDiff}). The case where $(\varepsilon_{j},\varepsilon_{k}%
)=(\ast,\ast)$ similarly matches the corresponding term in (\ref{momentDiff}).

We now turn to the case $(\varepsilon_{j},\varepsilon_{k})\in(1,\ast)$. We
compute
\begin{align*}
\operatorname{tr}\left[  a_{r}^{\varepsilon_{1}}\cdots da_{r}^{\varepsilon
_{j}}\cdots da_{r}^{\varepsilon_{k}}\cdots a_{r}^{\varepsilon_{n}}\right]   &
=s^{\prime}\operatorname{tr}\left[  a_{r}^{\varepsilon_{1}}\cdots
a_{r}^{\varepsilon_{j}}dx^{\prime}\cdots dx^{\prime}a_{r}^{\varepsilon_{k}%
}\cdots a_{r}^{\varepsilon_{n}}\right] \\
&  +\operatorname{tr}\left[  a^{\varepsilon_{1}}\cdots\tilde{b}_{r}%
\,d\tilde{w}_{r}\,u_{rs^{\prime}}\cdots u_{rs^{\prime}}^{\ast}\,d\tilde{w}%
_{r}^{\ast}\,\tilde{b}_{r}^{\ast}\cdots a^{\varepsilon_{n}}\right]  .
\end{align*}
Using the It\^{o} rules (\ref{eq:scIto}) and (\ref{eq:3Ito}), we obtain%
\begin{align}
&  \operatorname{tr}\left[  a_{r}^{\varepsilon_{1}}\cdots da_{r}%
^{\varepsilon_{j}}\cdots da_{r}^{\varepsilon_{k}}\cdots a_{r}^{\varepsilon
_{n}}\right] \nonumber\\
&  =s^{\prime}\operatorname{tr}\left[  a_{r}^{\varepsilon_{1}}\cdots
a_{r}^{\varepsilon_{j}}a_{r}^{\varepsilon_{k}}\cdots a_{r}^{\varepsilon_{n}%
}\right]  \operatorname{tr}\left[  a_{r}^{\varepsilon_{j+1}}\cdots
a_{r}^{\varepsilon_{k-1}}\right] \nonumber\\
&  +s\operatorname{tr}\left[  a_{r}^{\varepsilon_{1}}\cdots\tilde{b}_{r}%
\tilde{b}_{r}^{\ast}\cdots a_{r}^{\varepsilon_{n}}\right]  \operatorname{tr}%
\left[  u_{rs^{\prime}}a_{r}^{\varepsilon_{j+1}}\cdots a_{r}^{\varepsilon
_{k-1}}u_{rs^{\prime}}^{\ast}\right] \nonumber\\
&  =(s+s^{\prime})\operatorname{tr}\left[  a_{r}^{\varepsilon_{1}}\cdots
a_{r}^{\varepsilon_{j}}a_{r}^{\varepsilon_{k}}\cdots a_{r}^{\varepsilon_{n}%
}\right]  \operatorname{tr}\left[  a_{r}^{\varepsilon_{j+1}}\cdots
a_{r}^{\varepsilon_{k-1}}\right]  , \label{2dAr3}%
\end{align}
where, in the last equality, we have used both lines of (\ref{eq:bb*Combine}).
Note that the unitarity of $u_{rs^{\prime}}$ is used in an essential way here.
The result of (\ref{2dAr3}) matches the $\varepsilon_{j}=1,$ $\varepsilon
_{k}=\ast$ term in (\ref{momentDiff}). The case $(\varepsilon_{j}%
,\varepsilon_{k})=(\ast,1)$ is similar. Summing up all the four cases and
(\ref{eq:aFirstSum}), we see that $\operatorname{tr}\left[  a_{r}%
^{\varepsilon_{1}}\cdots da_{r}^{\varepsilon_{j}}\cdots da_{r}^{\varepsilon
_{k}}\cdots a_{r}^{\varepsilon_{n}}\right]  $ satisfies the same equation as
in (\ref{momentDiff}), after changing $b_{r}$ to $a_{r}.$

By the induction hypothesis, all the lower-order $\ast$-moments of $a_{r}$ and
$b_{r}$ are equal. Thus, $\operatorname{tr}\left[  a_{r}^{\varepsilon_{1}%
}\cdots a_{r}^{\varepsilon_{n}}\right]  $ satisfies the same ODE as
$\operatorname{tr}\left[  b_{r}^{\varepsilon_{1}}\cdots b_{r}^{\varepsilon
_{n}}\right]  $ as in (\ref{momentDiff2}), with the same constant $C$ and the
same function $f$. Since, also, both functions equal $1$ at $r=0,$ they are
equal for all $r.$
\end{proof}

\subsection{Varying $\tau$ with $s$ fixed}

In this section, we carry out Step 2 in the proof of
Theorem~\ref{thm:factorization}, by looking at how the $\ast$-moments of
$b_{s,\tau}$---multiplied on the left and right by freely independent
elements---vary as $\tau$ varies with $s$ fixed. When doing the calculation,
we will obtain a convenient cancellation that eliminates certain problematic
terms, a cancellation that would not occur if we were varying both $s$ and
$\tau$ simultaneously. This computation supports the overall philosophy of
this paper, that we obtain the nicest results when $\tau$ is varied while
keeping $s$ fixed.

\begin{lemma}
\label{lem:tmoment} Fix $s$ and $\tau$ with $\left\vert \tau-s\right\vert \leq
s$ and fix a sequence $(\varepsilon_{1},\ldots,\varepsilon_{n})$ with values
in $\{1,\ast\}$. Suppose $a_{1},a_{2}\in\mathcal{A}$ and $b_{s,t\tau}(r)$ are
all freely independent. Define
\[
B=a_{1}b_{s,t\tau}a_{2}%
\]
and
\[
f(s,t)=\operatorname{tr}[B^{\varepsilon_{1}}\cdots B^{\varepsilon_{k}}].
\]
Then for all $t<1$, we have
\begin{align}
\frac{\partial f}{\partial t}  & =\frac{f(s,t)}{2}\sum_{j}\tau^{\varepsilon
_{j}}\nonumber\\
& -\sum_{\varepsilon_{j}=\varepsilon_{k}}\tau^{\varepsilon_{j}}%
\operatorname{tr}[B^{\varepsilon_{1}}\cdots B^{\varepsilon_{j}}B^{\varepsilon
_{k+1}}\cdots B^{\varepsilon_{n}}]\operatorname{tr}[B^{\varepsilon_{j+1}%
}\cdots B^{\varepsilon_{k}}],\label{dfdt}%
\end{align}
where $\tau^{\varepsilon_{j}}=\bar{\tau}$ if $\varepsilon_{j}=\ast$ and where
the second sum on the right-hand side of (\ref{dfdt}) is over all pairs $j<k$
with $\varepsilon_{j}=\varepsilon_{k}.$
\end{lemma}

We emphasize that the second sum on the right-hand side of (\ref{dfdt})
contains only terms where $\varepsilon_{j}=\varepsilon_{k}$; the terms with
$\varepsilon_{j}\neq\varepsilon_{k}$ cancel out in the process of computing
the derivative with respect to $t$.

\begin{proof}
Let $B_{r}=a_{1}b_{s,t\tau}(r)a_{2},$ so that $B=B_{1}$. Since $b_{s,t\tau
}(r)$ has the same $\ast$-distribution as $b_{rs,rt\tau},$ we have
\begin{equation}
f(rs,rt)=\operatorname{tr}[B_{r}^{\varepsilon_{1}}\cdots B_{r}^{\varepsilon
_{k}}]. \label{fScaled}%
\end{equation}
Then by the chain rule, we have%
\[
\frac{d}{dr}f(rs,rt)=s\frac{\partial f}{\partial s}(rs,rt)+t\frac{\partial
f}{\partial t}(rs,rt).
\]
Solving this for the $\partial f/\partial t$ term and setting $r=1$ gives%
\begin{equation}
\frac{\partial f}{\partial t}(s,t)=\frac{1}{t}\left.  \frac{d}{dr}%
f(rs,rt)\right\vert _{r=1}-\frac{s}{t}\frac{\partial f}{\partial s}(s,t).
\label{chain}%
\end{equation}
We now proceed to compute each of the two terms on the right-hand side of
(\ref{chain}).

We first differentiate $f(rs,rt)$ with respect to $r$, using (\ref{fScaled}).
Using the stochastic product rule and the It\^{o} formulas, we easily obtain
\begin{align}
&  \frac{d}{dr}f(rs,rt)\nonumber\\
&  =\sum_{j}\frac{\operatorname{tr}[B_{r}^{\varepsilon_{1}}\cdots
dB_{r}^{\varepsilon_{j}}\cdots B_{r}^{\varepsilon n}]}{dr}+\sum_{j<k}%
\frac{\operatorname{tr}[B_{r}^{\varepsilon_{1}}\cdots dB_{r}^{\varepsilon_{j}%
}\cdots dB_{r}^{\varepsilon_{k}}\cdots B_{r}^{\varepsilon_{n}}]}%
{dr}\nonumber\\
&  =-f(rs,rt)\sum_{j}\frac{s-t\tau^{\varepsilon_{j}}}{2}+\sum_{j<k}%
\frac{\operatorname{tr}[B_{r}^{\varepsilon_{1}}\cdots dB_{r}^{\varepsilon_{j}%
}\cdots dB_{r}^{\varepsilon_{k}}\cdots B_{r}^{\varepsilon_{n}}]}%
{dr}.\label{eq:dgdr}%
\end{align}
We will analyze the last sum in (\ref{eq:dgdr}) in a moment.

We next differentiate $f(s,t)$ with respect to $s.$ Since $\left\vert
t\tau-s\right\vert <s$ for any $t<1$, we can choose $s_{0}<s$ such that
$\left\vert t\tau-s_{0}\right\vert <s.$ Then by Proposition~\ref{prop:bfactor}%
, the element
\[
U_{s}:=a_{1}b_{s_{0},t\tau}u_{s-s_{0}}a_{2}%
\]
has the same $\ast$-distribution as $a_{1}b_{s,t\tau}a_{2}$. We then
calculate
\begin{align}
\frac{\partial f}{\partial s}(s,t)  &  =\sum_{j}\frac{\operatorname{tr}%
[U_{s}^{\varepsilon_{1}}\cdots dU_{s}^{\varepsilon_{j}}\cdots U_{s}%
^{\varepsilon_{n}}]}{ds}+\sum_{j<k}\frac{\operatorname{tr}[U_{s}%
^{\varepsilon_{1}}\cdots dU_{s}^{\varepsilon_{j}}\cdots dU_{s}^{\varepsilon
_{k}}\cdots U_{s}^{\varepsilon_{n}}]}{ds}\nonumber\\
&  =-\frac{n}{2}f(s,t)+\sum_{j<k}\frac{\operatorname{tr}[U_{s}^{\varepsilon
_{1}}\cdots dU_{s}^{\varepsilon_{j}}\cdots dU_{s}^{\varepsilon_{k}}\cdots
U_{s}^{\varepsilon_{n}}]}{ds}, \label{eq:dgds}%
\end{align}
where we have computed the first sum as
\[
\sum_{j=1}^{n}\frac{\operatorname{tr}[U_{s}^{\varepsilon_{1}}\cdots
dU_{s}^{\varepsilon_{j}}\cdots U_{s}^{\varepsilon_{n}}]}{ds}=-\frac{1}{2}%
\sum_{j=1}^{n}\operatorname{tr}[U_{s}^{\varepsilon_{1}}\cdots U_{s}%
^{\varepsilon_{n}}]=-\frac{n}{2}f(s,t).
\]
We will analyze the last sum in (\ref{eq:dgds}) in a moment.

Using (\ref{chain}) together with (\ref{eq:dgdr}) and (\ref{eq:dgds}), we
obtain
\begin{align}
&  \frac{\partial f}{\partial t}-\frac{f(s,t)}{2}\sum_{j}\tau^{\varepsilon
_{j}}\nonumber\\
&  =\sum_{j<k}\left(  \frac{1}{t}\left.  \frac{\operatorname{tr}%
[B_{r}^{\varepsilon_{1}}\cdots dB_{r}^{\varepsilon_{j}}\cdots dB_{r}%
^{\varepsilon_{k}}\cdots B_{r}^{\varepsilon_{n}}]}{dr}\right\vert
_{r=1}\right. \nonumber\\
&  -\left.  \frac{s}{t}\frac{\operatorname{tr}[U_{s}^{\varepsilon_{1}}\cdots
dU_{s}^{\varepsilon_{j}}\cdots dU_{s}^{\varepsilon_{k}}\cdots B_{r}%
^{\varepsilon_{n}}]}{ds}\right)  . \label{df2}%
\end{align}
We then analyze the sum on the right-hand side of (\ref{df2}) case by case,
starting with the critical cases where $(\varepsilon_{j},\varepsilon_{k})$
equals $(1,\ast)$ or $(\ast,1).$ If $(\varepsilon_{j},\varepsilon_{k}%
)=(1,\ast)$, then
\begin{align*}
&  \operatorname{tr}[B_{r}^{\varepsilon_{1}}\cdots dB_{r}^{\varepsilon_{j}%
}\cdots dB_{r}^{\varepsilon_{k}}\cdots B_{r}^{\varepsilon_{n}}]\\
&  =\operatorname{tr}[B_{r}^{\varepsilon_{1}}\cdots(a_{1}b_{s,t\tau
}(r)\,dw_{r}~a_{2})\cdots(a_{2}^{\ast}~dw_{r}^{\ast}~b_{s,t\tau}(r)^{\ast
}a_{1}^{\ast})\cdots B_{r}^{\varepsilon_{n}}]\\
&  =s\operatorname{tr}[B_{r}^{\varepsilon_{1}}\cdots a_{1}b_{s,t\tau
}(r)b_{s,t\tau}(r)^{\ast}a_{1}^{\ast}\cdots B_{r}^{\varepsilon_{n}%
}]\operatorname{tr}[a_{2}B_{r}^{\varepsilon_{j+1}}\cdots B_{r}^{\varepsilon
_{k-1}}a_{2}^{\ast}]~dr
\end{align*}
where $\operatorname{tr}[a_{2}B_{r}^{\varepsilon_{j+1}}\cdots B_{r}%
^{\varepsilon_{k-1}}a_{2}]$ means $\operatorname{tr}[a_{2}a_{2}^{\ast}]$ if
$k=j+1$. Similarly,
\begin{align}
&  \operatorname{tr}[U_{s}^{\varepsilon_{1}}\cdots dU_{s}^{\varepsilon_{j}%
}\cdots dU_{s}^{\varepsilon_{k}}\cdots U_{s}^{\varepsilon_{n}}]\nonumber\\
&  =\operatorname{tr}[U_{s}^{\varepsilon_{1}}\cdots(a_{1}b_{s_{0},t\tau
}u_{s-s_{0}}\,dx_{s}a_{2})\cdots(a_{2}^{\ast}dx_{s}u_{s-s_{0}}^{\ast}%
b_{s_{0},t\tau}^{\ast}a_{1}^{\ast})\cdots U_{s}^{\varepsilon_{n}}]\nonumber\\
&  =\operatorname{tr}[U_{s}^{\varepsilon_{1}}\cdots(a_{1}b_{s_{0},t\tau
}b_{s_{0},t\tau}^{\ast}a_{1}^{\ast})\cdots U_{s}^{\varepsilon_{n}%
}]\operatorname{tr}[a_{2}U_{s}^{\varepsilon_{j+1}}\cdots U_{s}^{\varepsilon
_{k-1}}a_{2}^{\ast}]~ds. \label{lotsOfUs}%
\end{align}

In the last line of (\ref{lotsOfUs}), we note that $b_{s_{0},t\tau}%
b_{s_{0},t\tau}^{\ast}=(b_{s_{0},t\tau}u_{s-s_{0}})(b_{s_{0},t\tau}u_{s-s_{0}%
})^{\ast}$. Then after recalling the factorization result in Proposition
\ref{prop:bfactor}, we get a cancellation:
\[
\frac{1}{t}\left.  \frac{\operatorname{tr}[B_{r}^{\varepsilon_{1}}\cdots
dB_{r}^{\varepsilon_{j}}\cdots dB_{r}^{\varepsilon_{k}}\cdots B_{r}%
^{\varepsilon_{n}}]}{dr}\right\vert _{r=1}-\frac{s}{t}\frac{\operatorname{tr}%
[U_{s}^{\varepsilon_{1}}\cdots dU_{s}^{\varepsilon_{j}}\cdots dU_{s}%
^{\varepsilon_{k}}\cdots U_{s}^{\varepsilon_{n}}]}{ds}=0.
\]
The above equality also holds for the case $(\varepsilon_{j},\varepsilon
_{l})=(\ast,1)$ by a similar computation that again uses the unitarity of
$u_{s-s_{0}}.$

We now consider the case $\varepsilon_{j}=\varepsilon_{l}=1$. In this case, we
can directly verify that
\begin{align*}
&  \frac{1}{t}\left.  \frac{\operatorname{tr}[B_{r}^{\varepsilon_{1}}\cdots
dB_{r}^{\varepsilon_{j}}\cdots dB_{r}^{\varepsilon_{k}}\cdots B_{r}%
^{\varepsilon_{n}}]}{dr}\right\vert _{r=1}\\
&  =\frac{s-t\tau}{t}\operatorname{tr}[B^{\varepsilon_{1}}\cdots
B^{\varepsilon_{j}}B^{\varepsilon_{k+1}}\cdots B_{r}^{\varepsilon_{n}%
}]\operatorname{tr}[B^{\varepsilon_{j+1}}\cdots B^{\varepsilon_{k}}]
\end{align*}
and
\begin{align*}
&  \frac{s}{t}\frac{\operatorname{tr}[U_{s}^{\varepsilon_{1}}\cdots
dU_{s}^{\varepsilon_{j}}\cdots dU_{s}^{\varepsilon_{k}}\cdots U_{s}%
^{\varepsilon_{n}}]}{ds}\\
&  =\frac{s}{t}\operatorname{tr}[U_{s}^{\varepsilon_{1}}\cdots U_{s}%
^{\varepsilon_{j}}U_{s}^{\varepsilon_{k+1}}\cdots U_{s}^{\varepsilon_{n}%
}]\operatorname{tr}[U_{s}^{\varepsilon_{j+1}}\cdots U_{s}^{\varepsilon_{k}}]\\
&  =\frac{s}{t}\operatorname{tr}[B^{\varepsilon_{1}}\cdots B^{\varepsilon_{j}%
}B^{\varepsilon_{k+1}}\cdots B^{\varepsilon_{n}}]\operatorname{tr}%
[B^{\varepsilon_{j+1}}\cdots B^{\varepsilon_{k}}].
\end{align*}
Subtracting these equations gives a cancellation of the terms involving $s/t.$
The above equations hold in the case $\varepsilon_{j}=\varepsilon_{l}=\ast$,
except with $s-t\tau$ replaced by $s-t\bar{\tau}$ in the first line. Combining
the four cases in (\ref{df2}), we obtain the result claimed in the proposition.
\end{proof}

\subsection{The proof of the main result}

We are now ready for the proof of our main factorization result.

\begin{proof}
[Proof of Theorem~\ref{thm:factorization}]Write%
\begin{align*}
A_{t,t^{\prime}}  &  =b_{s,t\tau}b_{s^{\prime},t^{\prime}\tau^{\prime}}\\
B_{t}  &  =b_{s+s^{\prime},t(\tau+\tau^{\prime})}.
\end{align*}
We want show that $A_{t,t}$ and $B_{t}$ have the same $\ast$-moments by
mathematical induction on the order $n$ of $\ast$-moments, starting from the
trivial case $n=0.$ In the induction step, we will obtain differential
equations for the $\ast$-moments of $A_{t,t}$ and $B_{t}$ with respect to $t.$

Assume all the $\ast$-moments of length less than $n$ are equal, consider a
sequence $(\varepsilon_{1},\ldots,\varepsilon_{n})$ taking values in
$\{1,\ast\},$ and consider the functions
\[
f(t)=\operatorname{tr}[A_{t,t}^{\varepsilon_{1}}\cdots A_{t,t}^{\varepsilon
_{k}}]
\]
and
\[
g(t)=\operatorname{tr}[B_{t}^{\varepsilon_{1}}\cdots B_{t}^{\varepsilon_{k}%
}].
\]
Since $u_{s}u_{s^{\prime}}^{\prime}$ has the same $\ast$-distribution as
$u_{s+s^{\prime}}$ where $u$ and $u^{\prime}$ are freely independent free
unitary Brownian motions, $f(0)=g(0)$. Our goal is to show that $f(t)=g(t)$
for $0\leq t\leq1$.

Let
\[
h(t,t^{\prime})=\operatorname{tr}[A_{t,t^{\prime}}^{\varepsilon_{1}}\cdots
A_{t,t^{\prime}}^{\varepsilon}].
\]
Applying Lemma~\ref{lem:tmoment} with $a_{1}=1$, $a_{2}=b_{s,t\tau}$, we have
\begin{align*}
\frac{\partial h}{dt}(t,t^{\prime})  & =\frac{h(t,t^{\prime})}{2}\sum_{j}%
\tau^{\varepsilon_{j}}\\
& -\sum_{\varepsilon_{j}=\varepsilon_{k}}\tau^{\varepsilon_{j}}%
\operatorname{tr}[A_{t,t^{\prime}}^{\varepsilon_{1}}\cdots A_{t,t^{\prime}%
}^{\varepsilon_{j}}A_{t,t^{\prime}}^{\varepsilon_{k+1}}\cdots A_{t,t^{\prime}%
}^{\varepsilon_{n}}]\operatorname{tr}[A_{t,t^{\prime}}^{\varepsilon_{j+1}%
}\cdots A_{t,t^{\prime}}^{\varepsilon_{k}}].
\end{align*}
Similarly, by applying Lemma~\ref{lem:tmoment} with $a_{1}=b_{s,t\tau}$,
$a_{2}=1$, and with $s^{\prime},\tau^{\prime}$ in place of $s,\tau$, we have
\begin{align*}
\frac{\partial h}{dt^{\prime}}  & =\frac{h(t,t^{\prime})}{2}\sum_{j}%
(\tau^{\prime})^{\varepsilon_{j}}\\
& -\sum_{\varepsilon_{j}=\varepsilon_{k}}(\tau^{\prime})^{\varepsilon_{j}%
}\operatorname{tr}[A_{t,t^{\prime}}^{\varepsilon_{1}}\cdots A_{t,t^{\prime}%
}^{\varepsilon_{j}}A_{t,t^{\prime}}^{\varepsilon_{k+1}}\cdots A_{t,t^{\prime}%
}^{\varepsilon_{n}}]\operatorname{tr}[A_{t,t^{\prime}}^{\varepsilon_{j+1}%
}\cdots A_{t,t^{\prime}}^{\varepsilon_{k}}].
\end{align*}
Since $f(t)=h(t,t)$, the chain rule tells us that
\begin{align}
\frac{df}{dt} &  =\frac{\partial h}{\partial t}(t,t)+\frac{\partial
h}{\partial t^{\prime}}(t,t)\nonumber\\
&  =\frac{f(t)}{2}\sum_{j}(\tau+\tau^{\prime})^{\varepsilon_{j}}\nonumber\\
&  \qquad-\sum_{\varepsilon_{j}=\varepsilon_{k}}(\tau+\tau^{_{^{\prime}}%
})^{\varepsilon_{j}}\operatorname{tr}[A_{t,t}^{\varepsilon_{1}}\cdots
A_{t,t}^{\varepsilon_{j}}A_{t,t}^{\varepsilon_{k+1}}\cdots A_{t,t}%
^{\varepsilon_{n}}]\operatorname{tr}[A_{t,t}^{\varepsilon_{j+1}}\cdots
A_{t,t}^{\varepsilon_{k}}].\label{eq:fmomentODE}%
\end{align}

Meanwhile, by directly applying Lemma~\ref{lem:tmoment}, we find that the
function $g$ satisfies the same formula as in (\ref{eq:fmomentODE}), but with
$A_{t,t}$ replaced by $B_{t}$ everywhere. Now, by the induction hypothesis,
all the lower moments of $A_{t,t}$ and $B_{t}$ are equal. We therefore
conclude that $f$ and $g$ satisfy the same first-order linear ODE. Since we
also have $f(0)=g(0),$ as we have noted, we conclude that $f(t)=g(t)$ for all
$t<1$. Equality also holds when $t=1$ by the continuous dependence of
$b_{s,\tau}$ on $\tau$ (Appendix \ref{continuousDependence.sec}). This
completes the proof using mathematical induction, and we conclude that
$A_{t,t}$ and $B_{t}$ have the same $\ast$-distribution.
\end{proof}

\section{Continuous dependence of $b_{s,\tau}$ on $\tau$%
\label{continuousDependence.sec}}

We consider a free SDE of the form%
\[
db_{r}=b_{r}(\alpha~dx_{r}+\beta~d\tilde{x}_{r}+\gamma~dr)
\]
with initial condition $b_{0},$ where $\alpha,$ $\beta,$ and $\gamma$ are
fixed complex numbers. In integral form, this reads as%
\begin{equation}
b_{r}=b_{0}+\alpha\int_{0}^{r}b_{u}~dx_{u}+\beta\int_{0}^{r}b_{u}~d\tilde
{x}_{u}+\gamma\int_{0}^{r}b_{u}~du. \label{integralForm}%
\end{equation}
We will show that the solution depends continuously on the constants $\alpha,$
$\beta,$ and $\gamma.$ We first give an a priori bound on the size of the solutions.

\begin{proposition}
\label{brEst.prop}Given any $r_{0}>0$, we have the inequality
\[
\Vert b_{r}\Vert^{2}\leq\Vert b_{0}\Vert^{2}\left(  1+(2\sqrt{2}%
(|\alpha|+\left\vert \beta\right\vert )+\sqrt{r_{0}}|\gamma|)^{2}\right)
e^{\left(  1+(2\sqrt{2}(|\alpha|+\left\vert \beta\right\vert )+\sqrt{r_{0}%
}|\gamma|)^{2}\right)  r}%
\]
for all $r\leq r_{0}$.
\end{proposition}

\begin{proof}
We take the norm of both sides of (\ref{integralForm}) and apply the free
Burkholder--Davis--Gundy inequality \cite[Theorem 3.2.1]{BS1}, giving
\begin{equation}
\Vert b_{r}\Vert\leq\Vert b_{0}\Vert+2\sqrt{2}(|\alpha|+\left\vert
\beta\right\vert )\left(  \int_{0}^{r}\Vert b_{u}\Vert^{2}\,du\right)
^{1/2}+|\gamma|\int_{0}^{r}\Vert b_{u}\Vert\,du. \label{intEst}%
\end{equation}
Using the Cauchy--Schwarz inequality twice (once for the last integral in
(\ref{intEst}) and once for numbers) we have
\begin{align}
\Vert b_{r}\Vert &  \leq\Vert b_{0}\Vert+(2\sqrt{2}(|\alpha|+\left\vert
\beta\right\vert )+|\gamma|\sqrt{r_{0}})\left(  \int_{0}^{r}\Vert b_{u}%
\Vert^{2}du\right)  ^{1/2}\nonumber\\
&  \leq\left(  \Vert b_{0}\Vert^{2}+\int_{0}^{r}\Vert b_{u}\Vert^{2}du\right)
^{1/2}\left(  1+(2\sqrt{2}(|\alpha|+\left\vert \beta\right\vert )+\sqrt{r_{0}%
}|\gamma|)^{2}\right)  ^{1/2}. \label{intEst2}%
\end{align}
If we square both sides of (\ref{intEst2}) and apply Gronwall's inequality
\cite[Theorem 1.2.2]{pach}, we obtain
\[
\Vert b_{r}\Vert^{2}\leq\Vert b_{0}\Vert^{2}\left(  1+(2\sqrt{2}%
(|\alpha|+\left\vert \beta\right\vert )+\sqrt{r_{0}}|\gamma|)^{2}\right)
e^{\left(  1+(2\sqrt{2}(|\alpha|+\left\vert \beta\right\vert )+\sqrt{r_{0}%
}|\gamma|)^{2}\right)  r}%
\]
for all $r\leq r_{0}$, as desired.
\end{proof}

Now consider $b_{r}$ and $b_{r}^{\prime}$ satisfying (\ref{integralForm}) with
constants $\alpha,$ $\beta,$ $\gamma$ and $\alpha^{\prime},$ $\beta^{\prime},$
$\gamma^{\prime},$ respectively, and define
\[
v(r)=\Vert b_{r}-b_{r}^{\prime}\Vert.
\]
Then we have the following estimate, showing that $b_{r}$ is close to
$b_{r}^{\prime}$ when $(\alpha,\beta,\gamma)$ is close to $(\alpha^{\prime
},\beta^{\prime},\gamma^{\prime}).$

\begin{theorem}
\label{continuousDep.thm}Fix $r_{0}$ and assume $\alpha,\beta,\gamma
,\alpha^{\prime},\beta^{\prime},\gamma^{\prime}$ lie in a disk of some fixed
radius. Then there exists a constant $C>0$ such that
\[
v(r)^{2}\leq C\left(  2\sqrt{2}(\left\vert \alpha-\alpha^{\prime}\right\vert
+\left\vert \beta-\beta^{\prime}\right\vert )+\sqrt{r_{0}}\left\vert
\gamma-\gamma^{\prime}\right\vert \right)  ^{2}e^{Cr},\quad r\leq r_{0}.
\]
In particular, $v(r)\rightarrow0$ uniformly in $r\in\lbrack0,r_{0}]$ as
$(\alpha,\beta)\rightarrow(\alpha^{\prime},\beta^{\prime})$.
\end{theorem}

\begin{proof}
We compute
\begin{align*}
b_{t}-b_{t}^{\prime}  &  =\int_{0}^{t}[(\alpha-\alpha^{\prime})b_{r}%
+\alpha^{\prime}(b_{r}-b_{r}^{\prime})]\,dx_{r}+\int_{0}^{t}[(\beta
-\beta^{\prime})b_{r}+\beta^{\prime}(b_{r}-b_{r}^{\prime})]\,d\tilde{x}_{r}\\
&  +\int_{0}^{t}[(\gamma-\gamma^{\prime})b_{r}+\gamma^{\prime}(b_{r}%
-b_{r}^{\prime})]\,dr.
\end{align*}
By the free Burkholder--Davis--Gundy inequality,
\begin{align*}
v(r)  &  \leq2\sqrt{2}\left(  \int_{0}^{r}[(\left\vert \alpha-\alpha^{\prime
}\right\vert +\left\vert \beta-\beta^{\prime}\right\vert )\Vert b_{u}%
\Vert+(\left\vert \alpha^{\prime}\right\vert +\left\vert \beta^{\prime
}\right\vert )v(u)]^{2}\,du\right)  ^{1/2}\\
&  +\int_{0}^{r}(\left\vert \gamma-\gamma^{\prime}\right\vert \Vert b_{u}%
\Vert+\left\vert \gamma^{\prime}\right\vert v(u))\,du.
\end{align*}
We now apply the Minkowski inequality to the first integral and
Cauchy--Schwarz inequalities to the second integral and collect terms to
obtain%
\begin{align*}
v(r)  &  \leq\left(  2\sqrt{2}(\left\vert \alpha-\alpha^{\prime}\right\vert
+\left\vert \beta-\beta^{\prime}\right\vert )+\sqrt{r_{0}}\left\vert
\gamma-\gamma^{\prime}\right\vert \right)  \left(  \int_{0}^{r}\Vert
b_{u}\Vert^{2}du\right)  ^{1/2}\\
&  +\left(  2\sqrt{2}(\left\vert \alpha^{\prime}\right\vert +\left\vert
\beta\right\vert )+\sqrt{r_{0}}\left\vert \gamma^{\prime}\right\vert \right)
\left(  \int_{0}^{r}v(u)^{2}\,du\right)  ^{1/2}.
\end{align*}

Applying the two-dimensional Cauchy--Schwarz inequality to the above estimate,
we have
\begin{align}
v(r)^{2}  &  \leq\left[  \left(  2\sqrt{2}(\left\vert \alpha-\alpha^{\prime
}\right\vert +\left\vert \beta-\beta^{\prime}\right\vert )+\sqrt{r_{0}%
}\left\vert \gamma-\gamma^{\prime}\right\vert \right)  ^{2}+\int_{0}%
^{r}v(u)^{2}du\right] \nonumber\\
&  \times\left[  \int_{0}^{r}\Vert b_{u}\Vert^{2}du+\left(  2\sqrt
{2}(\left\vert \alpha^{\prime}\right\vert +\left\vert \beta^{\prime
}\right\vert )+\sqrt{r_{0}}\left\vert \gamma^{\prime}\right\vert \right)
^{2}\right]  . \label{vtEst}%
\end{align}
We may now bound the second factor on the right-hand side of (\ref{vtEst}) by
a constant $C$ using Proposition \ref{brEst.prop}. Applying Gronwall's
inequality gives the claimed estimate.
\end{proof}

\subsection*{Acknowledgments}

The authors thank Hari Bercovici, Guillaume C\'{e}bron, Bruce Driver, Todd
Kemp, Vaki Nikitopoulos, and Ping Zhong for useful discussions. We especially
thank Bruce Driver for showing us the proof of Lemma \ref{Bruce.lem}. We also
thank the referee for reading the paper carefully and making useful suggestions.

\subsection*{Data availability}

This article has no associated data.

\end{document}